\documentclass[paper=a4,headings=normal,DIV=16]{scrartcl}
\usepackage[T1]{fontenc}
\usepackage[utf8]{inputenc}
\usepackage[english]{babel}
\usepackage{
  amsmath,
  amssymb,
  amsfonts,
  amsxtra,
  amsthm,
  array,
  authblk,
  babelbib,
  booktabs,
  graphicx,
  mathrsfs,
  microtype,
  float,
  tikz,
  xcolor
}

\usetikzlibrary{arrows.meta, calc, decorations.pathreplacing}

\usepackage{hyperref}
\hypersetup{
  colorlinks=true,
  linkcolor=blue!80!black,
  citecolor=green!50!black
}

\numberwithin{equation}{section}
\theoremstyle{plain}
\newtheorem{thm}{Theorem}[section]
\newtheorem{lem}[thm]{Lemma}
\newtheorem{prop}[thm]{Proposition}
\newtheorem{cor}[thm]{Corollary}
\theoremstyle{definition}
\newtheorem{defi}[thm]{Definition}
\newtheorem{ex}[thm]{Example}
\theoremstyle{remark}
\newtheorem{rem}[thm]{Remark}

\DeclareMathOperator{\Curl}{curl}
\DeclareMathOperator{\Div}{div}
\DeclareMathOperator{\Grad}{grad}
\DeclareMathOperator{\A}{\mathcal{A}}
\DeclareMathOperator{\dom}{dom}
\DeclareMathOperator{\ran}{ran}

\newcommand{\mult}{\mathord{\cdot}}
\DeclareMathOperator{\supp}{supp}
\newcommand{\N}{\mathbb{N}}

\newcommand{\R}{\mathbb{R}}
\newcommand{\C}{\mathbb{C}}
\newcommand{\dd}[1]{\mathop{\mathrm{d}#1}}
\newcommand{\coleq}{\mathrel{\mathop:}=}
\newcommand{\eqcol}{=\mathrel{\mathop:}}
\newcommand{\ie}{\mbox{i.e.}}

\renewcommand{\Re}{\operatorname{Re}}

\renewcommand{\setminus}{\smallsetminus}
\newcommand{\BLO}{\mathcal{B}} 
\newcommand{\Hw}{L^2} 
\newcommand{\Hs}{\mathcal{H}} 
\newcommand{\HX}{X} 
\newcommand{\spatH}[1]{\mathcal{H}^{#1}} 
\newcommand{\Lip}{\mathrm{Lip}}
\newcommand{\snl}{q}
\newcommand{\ph}{N}

\title{\Large Well-Posedness and Exponential Stability of Nonlinear Maxwell Equations for Dispersive Materials with Interface}
\author[1]{Tomá\v{s} Dohnal}
\author[1]{Mathias Ionescu-Tira\footnote{mathias.ionescu-tira@mathematik.uni-halle.de}}
\author[2]{Marcus Waurick}
\affil[1]{\large Institute for Mathematics, Martin-Luther University Halle-Wittenberg}
\affil[2]{\large Institute for Applied Analysis, TU Bergakademie Freiberg}

\date{\today}

\begin{document}
\maketitle

\begin{abstract}
  \noindent  \textbf{\abstractname.}
  \noindent In this paper we consider an abstract Cauchy problem for a Maxwell system modelling electromagnetic fields in the presence of an interface between optical media.
  The electric polarization is in general time-delayed and nonlinear, turning the macroscopic Maxwell equations into a system of nonlinear integro-differential equations.
  Within the framework of evolutionary equations, we obtain well-posedness in function spaces exponentially weighted in time and of different spatial regularity and formulate various conditions on the material functions, leading to exponential stability on a bounded domain.
  \bigskip




\end{abstract}

\section{Introduction}

The macroscopic Maxwell equations governing electromagnetic fields $E = E(t,x)$, $H = H(t,x)$ ($t \in \R$, $x\in\Omega\subseteq\R^3$) in matter are given by
\begin{equation}
    \begin{aligned}
      \partial_t {D} - \Curl {H} &= -J &\qquad \Div D &= \rho \\
      \partial_t {B} + \Curl {E} &= 0 &\qquad \Div B &= 0,
    \end{aligned}
    \label{eq:Maxwell_diff}
\end{equation}
where $J$ and $\rho$ are the current and charge densities.
The auxiliary fields (the material response) $D, B$ are induced by $E,H$ through the constitutive relations
\begin{equation*}
  D(E) = \epsilon_0 E + P(E), \qquad B(H) = \mu H,
\end{equation*}
with a real constant $\epsilon_0 > 0$ and a symmetric matrix-valued function $\mu\colon \Omega \to \R^{3\times 3}$.
In reality, most material response is time-delayed, but the delay is often short enough to be modeled by instantaneous coefficients.
If the material includes metals (e.g. in plasmonics), the memory effect (and the implied dispersion, \ie, frequency dependence of the material functions in frequency domain) is significant enough to warrant a detailed discussion, see~\cite{MaierPlasmonics}.
Materials with memory also occur in nonlinear optics;
usually (see \cite{Boyd}), an analytic expansion
\begin{equation*}
  P(E) = {P}^{(1)}(E) + {P}^{(2)}(E) + {P}^{(3)}(E) + \cdots
\end{equation*}
is assumed, where each $P^{(n)}$ is given by the time-delayed action of a tensor with rank $(n+1)$.
We thus consider throughout variants of the model
\begin{equation}
  P(E)(t) = \int_\R\cdots\int_\R
  \chi(t - s_1, \ldots, t - s_k)\,
  \snl(E(s_1), \ldots, E(s_k)) \dd{s_1}\cdots\dd{s_k},
  \label{eq:Polarization}
\end{equation}
with $k\in \N$, where the kernel $\chi$ is causal in the sense that $\chi(s_1,\ldots,s_k) = 0$ whenever $s_j < 0$ for some $j \in \{1,\dots,k\}$.
Candidates for the vector function $\snl$ are subjected to various conditions involving Lipschitz-continuity.
Due to the convolution in \eqref{eq:Polarization} the dielectric constant, obtained from the Fourier transform (in time) of $\chi$, is frequency dependent. In other words, the material is dispersive. In the linear case, \ie, $k = 1$ and $q(E) = E$, a standard model of $\chi$ is the Drude--Lorentz model, see Appendix~\ref{sec:app1}.

\begin{figure}[htbp]
  \centering
  \begin{tikzpicture}[line join=round,scale=1]
    \fill[color=olive,opacity=.2] (0,0) --++ (6,0) --++ (0,1) --++ (-6,0) -- cycle; 
    \fill[color=olive,opacity=.3] (6,1) --++ (0,-1) --++ (40:1.4) --++ (0,1) -- cycle; 
    \fill[color=olive,opacity=.1] (0,1) --++ (0,-1) --++ (40:1.4) --++ (0,1) -- cycle; 
    \fill[color=olive,opacity=.15] (40:1.4) --++ (6,0) --++ (0,1) --++ (-6,0) -- cycle; 
    \fill[color=olive,opacity=.15] (0,1) --++ (6,0) --++ (40:1.4) --++ (-6,0) -- cycle; 

    \fill[green,opacity=.2] (0,0) --++ (6,0) --++ (40:1.4) --++ (-6,0) -- cycle; 

    \fill[color=blue,opacity=.2] (0,-1) --++ (6,0) --++ (0,1) --++ (-6,0) -- cycle; 
    \fill[color=blue,opacity=.3] (6,0) --++ (0,-1) --++ (40:1.4) --++ (0,1) -- cycle; 
    \fill[color=blue,opacity=.1] (0,0) --++ (0,-1) --++ (40:1.4) --++ (0,1) -- cycle; 
    \fill[color=blue,opacity=.15] (40:1.4) ++ (0,-1) --++ (6,0) --++ (0,1) --++ (-6,0) -- cycle; 
    \fill[color=blue,opacity=.2] (0,-1) --++ (6,0) --++ (40:1.4) --++ (-6,0) -- cycle; 

    \path (6,0) ++ (40:1.4) node[right] {$\Gamma$};

    \node[right] at (3.3,1.4) {$\Omega_1$};
    \node[right] at (3.3,-.7) {$\Omega_2$};

    \node[color=purple] at (2,.6) {\tiny $+\,+$};
    \node[color=purple] at (2,.4) {\tiny $+$};

    \node[color=violet] at (2.8,.6) {\tiny $-$};
    \node[color=violet] at (2.8,.4) {\tiny $-\,-$};

    \node[color=purple] at (3.6,.6) {\tiny $+\,+$};
    \node[color=purple] at (3.6,.4) {\tiny $+$};

    \node[color=violet] at (4.4,.6) {\tiny $-$};
    \node[color=violet] at (4.4,.4) {\tiny $-\,-$};

    \node[color=purple] at (5.2,.6) {\tiny $+\,+$};
    \node[color=purple] at (5.2,.4) {\tiny $+$};
  \end{tikzpicture}
  \quad
  \begin{tikzpicture}
    \node[right] at (0,.8)  {$\Omega_1$};
    \node[right] at (0,-.8) {$\Omega_2$};
    \fill[color=blue,opacity=.2] (0,0) rectangle (7,-1.4);
    \fill[color=olive,opacity=.2] (0,0) rectangle (7,1.4);
    \draw[gray,very thick,->] (6,0) --++ (0,-1) node[right] {$n$};
    \draw[very thick] (0,0) node[left] {$\Gamma$} -- (7,0);
    \draw[thick,->] (2.1,.2) sin +(.4,1) cos (2.9,.2);
    \draw[thick,<-] (3,.2)   sin +(.4,1) cos (3.8,.2);
    \draw[thick,->] (3.9,.2) sin +(.4,1) cos (4.7,.2);
    \draw[thick,->] (2.1,-.2) sin +(.4,-.4) cos (2.9,-.2);
    \draw[thick,<-] (3,-.2)   sin +(.4,-.4) cos (3.8,-.2);
    \draw[thick,->] (3.9,-.2) sin +(.4,-.4) cos (4.7,-.2);
    \draw[gray,very thick,->] (2.5,-1) --++ (2,0);
  \end{tikzpicture}
  \caption{Schematic depiction of the electric field of a travelling surface wave induced by charge oscillations at an interface $\Gamma$ between two media $\Omega_1$ and $\Omega_2$.}
  \label{fig:SPP_schematic}
\end{figure}
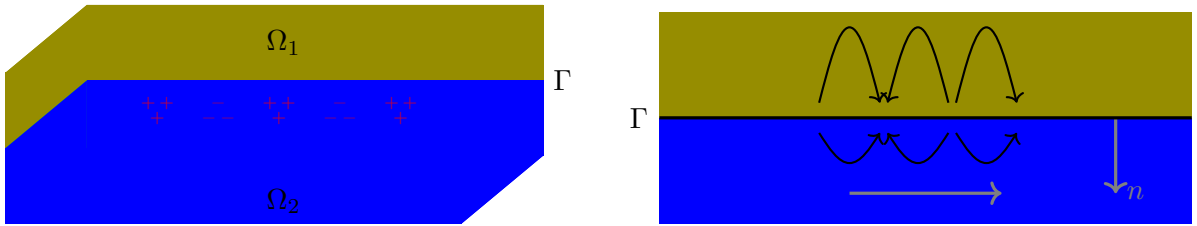

Our aim is to study Maxwell's equations in nonlinear optics for the interface problem depicted in Figure \ref{fig:SPP_schematic},
where $\Omega = \Omega_1 \sqcup \Gamma \sqcup \Omega_2$ ($\sqcup$ denotes disjoint union) consists of two domains $\Omega_1$ and $\Omega_2$, each with their own material response, separated by the interface $\Gamma$.
In this case, the two systems resulting from \eqref{eq:Maxwell_diff} (on $\Omega_i$ respectively) have to be supplemented by transmission conditions:
If $\Gamma$ is a 2-dimensional $C^1$-manifold with a normal vector field $n$, these are given by (see~\cite{DautrayLions1})
\begin{equation}
  \left[ {n} \times {E} \right]_\Gamma = \left[ {n} \times {H} \right]_\Gamma = 0, \quad \left[ {n} \cdot {D} \right]_\Gamma = \left[ {n} \cdot {B} \right]_\Gamma = 0
    \label{eq:interface_cond}
\end{equation}
and are due to the absence of surface charges and currents.
Here $\left[ f \right]_\Gamma \coleq \left.(f_1 - f_2)\right|_{\Gamma}$ denotes the jump of $f$ across $\Gamma$ (in the sense of traces), where $f_i = f\rvert_{\Omega_i}$.
Condition \eqref{eq:interface_cond} is equivalent to the continuity of the tangential components of $E,H$ and transversal components of $D,B$.

A possible application of this interface setting lies in the modelling of surface waves such as surface plasmon polaritons (SPPs).
These are evanescent electromagnetic fields resulting from charge excitation (achieved through coupling mechanisms) at an interface; see also Figure \ref{fig:SPP_schematic}.
SPPs exist at metal-dielectric interfaces: in the linear case, one can obtain the existence of travelling surface waves which satisfy a highly nonlinear dispersion relation due to frequency-dependent material response (see \cite{raether}).
In addition, as metals are intrinsically lossy, SPPs experience exponentially fast damping in time. In the theory this is reflected by the so-called exponential stability.

An important nonlinear effect at interfaces is {surface second-harmonic generation}, which originates from quadratic nonlinearities, see \cite{Shen1989}.
In fact, the quadratic nonlinearity is typically of leading order at interfaces where the inversion symmetry of the material is broken.
Analytical and numerical studies of this effect have been performed, for instance, in \cite{BaoDobsonSHG, AmmariBaoHamdache}.
Our results cover a wide class of nonlinearities, including the quadratic one.
\medskip

For nonlinear Maxwell systems with interface, a local well-posedness theory is available for materials
without memory, see~\cite{SchnaubeltSpitz2022} and the references therein.
An application of the latter to surface waves can be found in \cite{DST}.
The Maxwell problem with memory can be reformulated as an instantaneous system if the susceptibilities  satisfy certain assumptions, see \cite{SchneiderUecker2003}.
Up to our knowledge the only rigorous well-posedness analysis of nonlinear Maxwell equations with memory in the literature is the recent \cite{SchnaubeltBresch}. 
Here a semigroup approach is used and the results are limited to the local well-posedness on bounded domains. Our well-posedness results are in some cases global and include exponential stability.
\medskip

There are several aims of this article.
First, we provide a well-posedness theory (global in time) for Maxwell systems with interface and with a nonlinear material response given by nonlocal models.
The Maxwell system is formulated within the framework of evolutionary equations in the sense of Picard \cite{MilaniPicard, PicardStructObserv} (we also refer to \cite{STW22} as a general reference). As this theory works in Bochner spaces in space-time, memory effects can be treated more naturally (see also \cite{SuessWaurick} for similar nonlocal models in SPDEs).

Second, the formulation in spaces of higher spatial regularity allows for a wider class of nonlinearities and complements previous work on spatial regularity for evolutionary equations, see \cite{PTW_MaxReg} and also \cite{TW_NonautoMaxReg}.

Third, conditions for exponential stability of linear and nonlinear systems are provided; this is based on work in \cite{TrostorffLinExpStab, Trostorff_habil}.

Finally, treating the paradigmatic Maxwell case may open up similar strategies for general evolutionary equations.

\bigskip

This article is structured as follows.
Section \ref{sec:wp} is concerned with the well-posedness of the Maxwell system in the functional analytic framework of evolutionary equations. Since here we work exclusively with (weighted) $L^2$-spaces, no regularity of the interface or the boundary is needed.
Subsection \ref{sec:higher_spatial_regularity} deals with higher regularity in space, which requires some regularity of the interface and the boundary as well as the boundedness of $\Omega$ (whereas the time regularity is covered by the general theory and is not Maxwell-specific).

In Section \ref{sec:exp_stab} we examine exponential stability for the electric and magnetic field.
The result for the linear second-order formulation (wave equation for the electric field) is Theorem \ref{thm:Maxwell_ExpStab_SecondOrder}. The proof uses the theory in \cite{Trostorff_habil} and the Picard--Weber--Weck selection theorem (\ref{thm:PWW_selection_thm}) as a key ingredient.
The latter restricts $\Omega$ to a bounded (weak Lipschitz) domain.
Imposing suitable Lipschitz-continuity on the nonlinearity, a corresponding result is obtained for the nonlinear system by a fixed-point argument.

Using similar methods as in Theorem~\ref{thm:Maxwell_ExpStab_SecondOrder}, we establish exponential stability for the full Maxwell system, for materials of a different class, in Theorem~\ref{thm:Maxwell_ExpStab_FirstOrder} for fields in $L^2(\Omega)^3$, and in Theorem~\ref{thm:Maxwell_ExpStab_H2} for fields in $H^2(\Omega_1)^3\oplus H^2(\Omega_2)^3$. Theorem~\ref{thm:Maxwell_ExpStab_H2} allows for a wider class of nonlinearities, which are considered in Theorem~\ref{thm:NonlinFixH2}.

The theory is accompanied by several examples. An overview of these applications is found in Appendix~\ref{sec:app_examples_overview}.

\section{The Maxwell system as an evolutionary equation}
\label{sec:wp}

\subsection{Maxwell operator and boundary conditions}

We introduce the functional analytic setup in which we treat the Maxwell system \eqref{eq:Maxwell_diff} together with the transmission conditions.
Let again $\Omega = \Omega_1 \sqcup \Gamma \sqcup \Omega_2$ and set
\begin{equation*}
  \Hs \coleq L^2(\Omega)^3.
\end{equation*}
We denote by $\Curl_0$ the closure of the operator
\begin{equation*}
  C_c^\infty(\Omega)^{3} \ni \varphi =
  \begin{pmatrix}
    \varphi_1 \\ \varphi_2 \\ \varphi_3
  \end{pmatrix}
  \mapsto
  \nabla\times\varphi =
  \begin{pmatrix}
    \partial_2 \varphi_{3} - \partial_{3}\varphi_{2} \\
    \partial_3 \varphi_{1} - \partial_{1}\varphi_{3} \\
    \partial_1 \varphi_{2} - \partial_{2}\varphi_{1}
  \end{pmatrix} \in C_c^\infty(\Omega)^{3}
\end{equation*}
in $\Hs$, with
\begin{equation*}
  \dom(\Curl_0) = H_0(\Curl,\Omega) \coleq
  \overline{C_c^\infty(\Omega)^3}^{\lVert \cdot \rVert_{H(\Curl)}},
\end{equation*}
where $\lVert \varphi \rVert_{H(\Curl)} = \bigl( \lVert \varphi \rVert_{L^2}^2 + \lVert \nabla\times \varphi \rVert_{L^2}^2 \bigr)^{1/2}$.
We further set $\Curl \coleq \Curl_0^*$.
It is then easy to derive
\begin{equation*}
  \dom(\Curl) = H(\Curl,\Omega) \coleq \{ u \in \Hs : \Curl u \in \Hs \}.
\end{equation*}
Consider the operator
\begin{equation*}
  \A \coleq
  \begin{pmatrix}
    0 & -\Curl \\ \Curl_0 & 0
  \end{pmatrix}.
\end{equation*}
Then, by construction (since $\Curl^* = (\Curl_0^*)^* = \Curl_0$),
\begin{equation*}
  \A^* = 
  \begin{pmatrix}
    0 & -\Curl \\ \Curl_0 & 0
  \end{pmatrix}^* =
  \begin{pmatrix}
    0 & \Curl_0^* \\ -\Curl^* & 0
  \end{pmatrix} =
  \begin{pmatrix}
    0 & \Curl \\ -\Curl_0 & 0
  \end{pmatrix} = -\A,
\end{equation*}
\ie, $\A \colon H_0(\Curl,\Omega)\times H(\Curl,\Omega) \subseteq L^2(\Omega)^3 \to L^2(\Omega)^3$ is \emph{skew-selfadjoint}.

\begin{rem}
  For the application of the well-posedness in Sections \ref{sec:well_posedness} and \ref{sec:initial_values} the domains $\Omega_1$, $\Omega_2$ may be quite general; in particular no regularity of the boundary is needed, and $\Omega_1$, $\Omega_2$ may be unbounded (if $\Omega = \R^3$, then $\Curl = \Curl_0$). However, to deal with higher regularity and exponential stability in Sections \ref{sec:higher_spatial_regularity} and \ref{sec:exp_stab} they need to be bounded with more regularity of the boundary.
\end{rem}

The above choice of 
\begin{equation*}
  \dom(\A) = H_0(\Curl,\Omega)\times H(\Curl,\Omega)
\end{equation*}
and the skew-selfadjointness of $\A$ encode the interface conditions and the boundary condition of a perfect conductor, if the boundaries of $\Omega_1$, $\Omega_2$ are sufficiently regular:
Assume that $\Omega$, $\Omega_1$, $\Omega_2$ have Lipschitz boundaries and denote their outward normal fields by $n$.
Let $\mathscr{D}\subseteq\dom(\A)$ be a subset consisting of functions that are smooth in $\overline{\Omega}_1$ and $\overline{\Omega}_2$ and fix $(u_E, u_H) \in \mathscr{D}$.
Using the divergence theorem on $\overline{\Omega}_1$ and $\overline{\Omega}_2$ separately, we have  for all $v_E \in C_c^\infty(\Omega)$, $v_H \in C^\infty(\overline{\Omega})$
\begin{align*}
  \int_{\Omega} \bigl( \Curl_0 u_E \cdot {v_H} - u_E \cdot {\Curl v_H} \bigr)
  &= \int_{\Omega_1} \Div (u_E \times {v_H}) + \int_{\Omega_2} \Div (u_E \times {v_H}) \\
  &= \int_{\partial\Omega_1} (u_E \times {v_H}) \cdot n + \int_{\partial\Omega_2} (u_E \times {v_H}) \cdot n \\
  &= \int_\Gamma \bigl[(u_E \times {v_H}) \cdot n\bigr]_\Gamma + \int_{\partial\Omega} (u_E \times v_H) \cdot n \\
  &= \int_\Gamma \bigl[n \times u_E\bigr]_\Gamma \cdot {v_H} + \int_{\partial\Omega} (n \times u_E) \cdot v_H,
\end{align*}
and similarly,
\begin{align*}
  \int_{\Omega} \bigl( \Curl u_H \cdot {v_E} - u_H \cdot {\Curl_0 v_E} \bigr)
  &= \int_\Gamma \bigl[n \times u_H\bigr]_\Gamma \cdot {v_E} + \int_{\partial\Omega} (n \times u_H) \cdot v_E. \\
  &= \int_\Gamma \bigl[n \times u_H\bigr]_\Gamma \cdot {v_E}.
\end{align*}
By skew-selfadjointness of $\A$, the left-hand sides must vanish for arbitrary $v_E,v_H$.
Therefore,
\begin{equation}
  \bigl[n \times u_E \bigr]_\Gamma = \bigl[n \times u_H \bigr]_\Gamma = 0
  \quad \text{and}\quad
  (n \times u_E)\rvert_{\partial\Omega} = 0.
  \label{eq:Hcurl_jump}
\end{equation}
The latter identity is the boundary condition of a perfect conductor.
Using the traces in $H(\Curl,\Omega)$ and $H_0(\Curl,\Omega)$, equations \eqref{eq:Hcurl_jump} can be shown to hold for $u_E \in H_0(\Curl,\Omega)$, $u_H\in H(\Curl,\Omega)$ in the sense of traces, see \cite{Leis, BDPW22}.
In absence of this regularity of $\partial\Omega$ and $\Gamma$, conditions \eqref{eq:Hcurl_jump} are interpreted in a generalized sense.
\bigskip

We consider throughout a Cauchy problem for the Maxwell system \eqref{eq:Maxwell_diff}, formulated as an evolutionary problem in (positive) time:
\begin{equation}
  \left\{
    \begin{aligned}
      \partial_t
      \begin{pmatrix}
        D(E) \\ B(H)
      \end{pmatrix} +
      \begin{pmatrix}
        0 & -\Curl \\ \Curl_0 & 0
      \end{pmatrix}
      \begin{pmatrix}
        E \\ H
      \end{pmatrix} &=
      \begin{pmatrix}
        -J \\ 0
      \end{pmatrix}, & t > 0\\
      (E(t), H(t)) &= (E_0(t), H_0(t)), & t \le 0
    \end{aligned}
  \right\}
  \label{eq:Maxwell_pos_init}
\end{equation}
in $L^2(\Omega)^6$, where
\begin{equation*}
  (E_0,H_0) \colon (-\infty,0] \to L^2(\Omega)^6,\quad t\mapsto (E_0(t), H_0(t))
\end{equation*}
is a given history (this is necessary since the material function $(D,B)$ is in general dependent on past values of its argument).
The solution should meet the condition $(E(t),H(t)) \in \dom(\A)$ for all $t$, in order for the jump conditions of $E,H$ to be fulfilled.

For the divergence equations for $D,B$ one finds that they are largely redundant, in the sense that they follow from \eqref{eq:Maxwell_pos_init} and suitable initial values.
Indeed (cf.~\cite{SchnaubeltSpitz2022}), applying $\Div$ to the first line in \eqref{eq:Maxwell_pos_init} and integrating, it follows that $\Div D(t) = \varrho(t)$ holds for $t \ge 0$ if and only if $\varrho$ and $J$ are related by
\begin{equation*}
  \varrho(t) = \varrho(0) - \int_0^t \Div J(s)\dd s.
\end{equation*}
Similarly, it follows from second line in \eqref{eq:Maxwell_pos_init} that $\Div B$ is constant for all $t > 0$, so if $\Div B(0) = 0$, then $\Div B(t) = 0$ holds for all $t > 0$.

Regarding the jump conditions
\begin{equation}
  [n\cdot D]_\Gamma = [n\cdot B]_\Gamma = 0,
  \label{eq:div_jump_BD}
\end{equation}
it suffices that they are fulfilled at time $t = 0$; then \eqref{eq:div_jump_BD} follows for all $t > 0$ by taking derivatives in time and using the structure of the remaining equations (see again~\cite{SchnaubeltSpitz2022}).

More generally, the jump conditions \eqref{eq:div_jump_BD} are a property of the domain of both the operators $\Div\colon H(\Div,\Omega) \subseteq L^2(\Omega)^3\to L^2(\Omega)$ and $\Div_0\colon H_0(\Div,\Omega) \subseteq L^2(\Omega)^3\to L^2(\Omega)$, and can thus be interpreted as a regularity condition, see Sections \ref{sec:initial_values}, \ref{sec:higher_spatial_regularity}.
Here $\Div = -\Grad_0^*$, $\Div_0 = -\Grad^*$ are defined similarly to $\Curl, \Curl_0$ in terms of the usual weak gradient $\Grad\colon H^1(\Omega) \subseteq L^2(\Omega)\to L^2(\Omega)^3$ and the weak gradient with zero boundary condition, $\Grad_0\colon H_0^1(\Omega) \subseteq L^2(\Omega)\to L^2(\Omega)^3$.

We have
\begin{equation*}
  H(\Div,\Omega) = \{ u \in L^2(\Omega)^3 : \Div u \in L^2(\Omega) \}
\end{equation*}
and $H_0(\Div,\Omega)$ is the closure of $C_c^\infty(\Omega)^3$ with respect to the norm $u\mapsto (\lVert u \rVert_{L^2}^2 + \lVert \nabla\cdot u \rVert_{L^2}^2)^{1/2}$.

\begin{rem}
  The interface setting will play no role in the solution theory established in the next sections, as the transmission conditions are naturally embedded into the domain $\dom(\A)$, and is thus independent of any inhomogeneities in the material.
  Only in the context of higher regularity in space (Section \ref{sec:higher_spatial_regularity}, and Section \ref{sec:ExpStab_H2} for exponential stability in those spaces) will we need to take the regularity of the interface into account.
\end{rem}

\subsection{Linear evolutionary equations}

In the following, we provide a short overview of the theory of evolutionary equations. For details, see \cite{STW22}.
We will first consider a purely linear material function (for example taking $k = 1$ in \eqref{eq:Polarization} and $\snl(u) = u$, see Example~\ref{ex:LinearPermittivity} below).
In this case, the linear Maxwell system fits into the category of abstract evolutionary equations of the form
\begin{equation}
  \bigl( \partial_t M(\partial_t) + A \bigr) u = g
  \label{eq:evo_lin_1order}
\end{equation}
with given data $g$, understood as an operator equation in the weighted Hilbert space
\begin{equation*}
  \Hw_\varrho(\R,\HX) \coleq \{ u \in L^2_{\mathrm{loc}}(\R, \HX) : \exp(-\varrho\,\mult)u \in L^2(\R,\HX) \}
\end{equation*}
with the inner product
\begin{equation*}
  {\langle u, v \rangle}_{\varrho,0} \coleq \int_\R {\langle u(t), v(t) \rangle}_{\HX}\, e^{-2\varrho t} \dd t
\end{equation*}
and norm $\lVert \cdot \rVert_{\varrho,0} = \sqrt{\langle \cdot,\,\cdot \rangle_{\varrho,0}}$,
for some $\varrho \in \R$.
Here $\HX$ is a Hilbert space and $A\colon \dom(A) \subseteq \HX \to \HX$ a densely defined and closed operator, extended to a subset of $\Hw_\varrho(\R,\HX)$ via $(A u)(t) \coleq A(u(t))$.
The time derivative $\partial_t\colon \dom(\partial_t) \subseteq \Hw(\R,\HX)\to \Hw(\R,\HX)$ is understood in the weak sense and is a densely defined and closed operator, where
\begin{equation*}
  \dom(\partial_t) = H^1_\varrho(\R,\HX) \coleq \{ u \in L^2_\varrho(\R,\HX) : u' \in L^2(\R,\HX) \}.
\end{equation*}
In most cases $M(\partial_t)$ will denote a convolution operator, but more generally it is a linear material law:
\begin{defi}
  \label{def:MaterialLaw}
  A \emph{linear material law} is an analytic mapping $M\colon \mathrm{dom}(M) \subseteq \C \to \mathcal{B}(\HX)$ into the space $\BLO(\HX)$ of bounded linear operators on $\HX$ (with norm denoted by $\lVert \mult \rVert$), which is uniformly bounded on a right half-plane, \ie,
  \begin{align*}
    \exists \varrho_0 \in \R\colon\quad
    \sup_{\Re z > \varrho_0}\lVert M(z) \rVert < \infty.
  \end{align*}
  The operator $M(\partial_t)$ is defined by the composition
  \begin{equation*}
    M(\partial_t) = \mathcal{L}_{\varrho}^* M(i\,\mult + \varrho)\mathcal{L}_\varrho,
  \end{equation*}
  with $M(i\,\mult + \varrho)$ defined by $(M(i\,\mult + \varrho) \varphi)(t) = M(it + \varrho)\varphi(t)$, and where $\mathcal{L}_{\varrho}\colon \Hw_\varrho(\R,\HX) \to L^2(\R,\HX)$ is the unitary extension of the Fourier--Laplace transform
  \begin{equation*}
    (\mathcal{L}_\varrho \varphi)(t) \coleq \frac{1}{\sqrt{2\pi}} \int_\R \varphi(s)\, e^{-(it + \varrho)s}\dd s.
  \end{equation*}
\end{defi}

The solution theory for \eqref{eq:evo_lin_1order} is established by Theorem \ref{thm:linear_solution_theory} and is closely tied to the concept of causality.
\begin{defi}
For $a \in \R$ we denote by $\theta_a^+ \colon \Hw_\varrho(\R,\HX) \to \Hw_\varrho(\R,\HX)$ the multiplication operator defined by
\begin{equation*}
  \theta_a^+ u(t) = 1_{(a, \infty)}(t) u(t) = 
  \begin{cases}
    u(t), & t > a \\ 0, & t \le a.
  \end{cases}
\end{equation*}
A mapping $f$ on $\Hw_\varrho(\R,\HX)$ is called \emph{(forward) causal}, if for all $a \in \R$ the implication
\begin{equation*}
  (1 - \theta_a^+)(u - v) = 0 \implies (1 - \theta_a^+)(f(u) - f(v)) = 0
\end{equation*}
holds, \ie, if $u, v$ agree on $(-\infty,a]$, then so do $f(u), f(v)$.
\end{defi}

By a consequence of the Paley--Wiener theorem (\cite[Theorem 8.1.2]{STW22}), if $M$ is a linear material law, then $M(\partial_t)$ is a causal operator on $\Hw_\varrho(\R,\HX)$ for $\varrho > \varrho_0$.
It is a key observation that $\partial_t$ is boundedly invertible for $\varrho \ne 0$, and causally invertible for $\varrho > 0$. In the latter case $\partial_t^{-1}\colon \Hw_\varrho(\R,\HX) \to \Hw_\varrho(\R,\HX)$ is given by
\begin{equation}
  (\partial_t^{-1}u)(t) = \int_{-\infty}^t u(s)\dd s
  \label{eq:Antiderivative}
\end{equation}
and satisfies $\lVert \partial_t^{-1} \rVert \le 1/\varrho$ (see \cite[Sec.~3.2]{STW22}). We use the symbol $\partial_t^{-1}$ exclusively to denote the causal map given by \eqref{eq:Antiderivative}.

\begin{thm}[Picard's Theorem, see e.g. {\cite[Theorem 6.2.1]{STW22}}]
  \label{thm:linear_solution_theory}
  Let $A\colon \dom(A) \subseteq \HX \to \HX$ be skew-selfadjoint and $M$ a linear material law, for which $zM(z)$ is strictly accretive on a half-plane $\C_{\Re > \varrho_0}$ with $\varrho_0 \in \R$, in the sense that
  \begin{equation}
    \exists c > 0\ \forall z \in \C_{\Re > \varrho_0}:
    \quad \Re zM(z) \ge c
    \label{eq:MaterialLawAccretivity}
  \end{equation}
  (\ie, $\Re \langle zM(z)x, x \rangle \ge c\lVert x \rVert_{\HX}^2$ for all $x \in \HX$).
  Then for all $\varrho > \varrho_0$ the operator ${\partial_tM(\partial_t) + A}$ is closable and
  \begin{equation*}
    S_\varrho \coleq {(\overline{\partial_tM(\partial_t) + A})}^{-1}\colon \Hw_\varrho(\R,\HX)\to \Hw_\varrho(\R,\HX)
  \end{equation*}
  is well-defined and bounded, with $\lVert S_\varrho \rVert_{\Hw_\varrho\to \Hw_\varrho} \le 1/c$.
  Moreover, $S_\varrho$ is causal and for all $g \in \Hw_\varrho(\R,\HX)$ the following implications hold:
  \begin{enumerate}
    \item[(i)] If $g \in H^1_\varrho(\R,\HX)$, then $S_\varrho g \in H^1_\varrho(\R,\HX) \cap \dom(A)$.
    \item[(ii)] If $g \in \Hw_\varrho(\R,\HX)\cap \Hw_{\varrho'}(\R,\HX)$, then $S_\varrho g = S_{\varrho'}g \in \Hw_\varrho(\R,\HX)\cap \Hw_{\varrho'}(\R,\HX)$ for all $\varrho,\varrho' > \varrho_0$.
  \end{enumerate}
\end{thm}

\begin{rem}
  If $M, A$ satisfy the assumptions of Theorem~\ref{thm:linear_solution_theory}, the solution $u$ of \eqref{eq:evo_lin_1order} is explicitly given using the spectral representation of the time-derivative:
\begin{equation*}
  u = {(\overline{\partial_tM(\partial_t) + A})}^{-1}g = \mathcal{L}_\varrho^{-1}\bigl((i\,\mult + \varrho) M(i\,\mult + \varrho) + A\bigr)^{-1}\mathcal{L}_\varrho g
\end{equation*}
for $g \in \Hw_\varrho(\R,\HX)$, $\varrho > \varrho_0$.
As such, Theorem \ref{thm:linear_solution_theory} provides sufficient conditions for the operator
\begin{equation*}
  {(\mult\, M(\mult) + A)}^{-1}\colon \C_{\Re > \varrho_0} \cap \dom(M) \to \BLO(\HX), \quad z \mapsto (zM(z) + A)^{-1}
\end{equation*}
to have a bounded and analytic extension on $\C_{\Re > \varrho_0}$.
In this case, we say that the problem \eqref{eq:evo_lin_1order} is \emph{well-posed in the range of spaces $L^2_\varrho(\R,\HX)$, $\varrho > \varrho_0$}, or simply \emph{well-posed}, implicitly presuming the existence of such $\varrho_0 \in \R$.
\end{rem}

\subsection{Well-posedness of the nonlinear Maxwell system}
\label{sec:well_posedness}

Nonlinear, (uniformly) Lipschitz-continuous perturbations of linear equations can be treated by a Banach fixed-point argument; see \cite{SuessWaurick} for a similar argument.
Precisely, inspired by \cite[Section 4.2]{STW22}, we consider the following class of nonlinearities.

\begin{defi}
  A function $f \colon \dom(f)\subseteq \bigcap_{\varrho \ge \varrho_0}\Hw_\varrho(\R,\HX) \to \bigcap_{\varrho \ge \varrho_0}\Hw_\varrho(\R,\HX)$ is called \emph{uniformly Lipschitz-continuous}, if $\dom(f)$ is dense in $\Hw_\varrho(\R,\HX)$ and $f$ extends to a Lipschitz-continuous map $f_\varrho \colon \Hw_\varrho(\R,\HX) \to \Hw_\varrho(\R,\HX)$ such that $f_\varrho = f_{\varrho'}$ on $\Hw_\varrho(\R,\HX)\cap \Hw_{\varrho'}(\R,\HX)$ for all $\varrho,\varrho' > \varrho_0$.
  In this case we simply write $f\colon \Hw_{\varrho}(\R,\HX) \to \Hw_{\varrho}(\R,\HX)$ for $\varrho > \varrho_0$.
\end{defi}

\begin{prop}
  \label{prop:Lipschitz_perturb}
  Let $\varrho_0, d \in \R_{> 0}$, and let $M$ be a linear material law satisfying
  \begin{equation}
    \forall z \in \C_{\Re > \varrho_0}\colon\quad \Re zM(z) \ge \frac{\Re z}{d}.
    \label{eq:Lipschitz_perturb}
  \end{equation}
  Furthermore, let $f\colon \Hw_{\varrho}(\R,\HX)\to \Hw_{\varrho}(\R,\HX)$ be causal and uniformly Lipschitz-continuous for all $\varrho > \varrho_0$ with
  \begin{equation*}
    \limsup_{\varrho \to +\infty} \frac{d}{\varrho}\lVert f \rVert_{\Lip(\Hw_\varrho \to \Hw_\varrho)}  < 1.
  \end{equation*}
  Then there exists $\varrho_1 \ge \varrho_0$ such that for all $\varrho > \varrho_1$ the problem
  $(\partial_tM(\partial_t) + A)u = f(u)$
  possesses a unique solution $u \in \Hw_{\varrho}(\R,\HX)$, which is independent of $\varrho$.
\end{prop}

\begin{proof}
  Denote again $S_\varrho = {(\overline{\partial_tM(\partial_t) + A})}^{-1}$.
  Due to \eqref{eq:Lipschitz_perturb} the material law $M$ satisfies a stronger variant of \eqref{eq:MaterialLawAccretivity}, namely $\Re zM(z) \ge \varrho/d$ for $\Re z > \varrho > \varrho_0$, hence $S_\varrho$ is bounded in $L^2_\varrho(\R,\HX)$ for all $\varrho > \varrho_0$, with $\lVert S_\varrho\rVert \le d/\varrho$.
  By Theorem \ref{thm:linear_solution_theory} and Lipschitz continuity of $f$ we can estimate
  \begin{equation*}
  \begin{aligned}
    \lVert S_\varrho f(u) - S_\varrho f(v) \rVert_{\varrho,0}
    &= \lVert {(\overline{\partial_tM(\partial_t) + A})}^{-1}(f(u) - f(v)) \rVert_{\varrho,0} \\
    &\le \frac{d}{\varrho} \lVert f \rVert_{\Lip} \lVert u - v \rVert_{\varrho,0}.
  \end{aligned}
  \end{equation*}
  Hence, $S_{\varrho}f(\cdot)\colon \Hw_{\varrho}(\R,\HX) \to \Hw_{\varrho}(\R,\HX)$ becomes a contraction for large $\varrho$.
  By Theorem~\ref{thm:linear_solution_theory} (ii) and the assumption on $f$, the unique fixed point is independent of $\varrho$.
\end{proof}

From now on we apply this result to the general nonlinear Maxwell system \eqref{eq:Maxwell_pos_init} setting $\HX = \mathcal{H}\times \mathcal{H}$.
To this end, we isolate the linear part of the polarization and take
\begin{equation*}
  D(E) = \epsilon(\partial_t)E + P_{\mathrm{nl}}(E),
\end{equation*}
where $\epsilon(\cdot)$ is a linear material law (the linear permittivity) and the resulting nonlinear system takes the form
\begin{equation}
  \left( \partial_t \begin{pmatrix} \epsilon(\partial_t) & 0 \\ 0 & \mu \end{pmatrix} + \begin{pmatrix} 0 & -\Curl \\ \Curl_0 & 0 \end{pmatrix} \right)
  \begin{pmatrix} E \\ H \end{pmatrix}
  = \begin{pmatrix} -\partial_tP_{\mathrm{nl}}(E) \\ 0 \end{pmatrix}
  + \begin{pmatrix} \varPhi \\ \varPsi \end{pmatrix}.
  \label{eq:Maxwell_nonlinear}
\end{equation}
The assumptions of Proposition~\ref{prop:Lipschitz_perturb} are satisfied if:
\begin{itemize}
  \item $\epsilon$ satisfies $\Re z\epsilon(z) \ge c_\epsilon\Re z$ for $\Re z > \varrho_0$ and some $c_\epsilon > 0$.
  \item $\partial_tP_\mathrm{nl}\colon \Hw_\varrho(\R,\Hs)\to \Hw_\varrho(\R,\Hs)$ is uniformly Lipschitz continuous  for $\varrho > \varrho_0$.
  \item The permeability $\mu \in \BLO(\Hs)$ is a bounded, selfadjoint operator which is uniformly positive definite, \ie, $\mu(x) \ge c_\mu > 0$ for all $x \in \Omega$.
  \item $\varPhi, \varPsi \in L^2_\varrho(\R,\Hs)$ ($\varrho > \varrho_0$) are arbitrary inhomogeneities (used below to encode the history of the system, see \eqref{eq:Maxwell_evo_inhomogeneity}).
\end{itemize}
Thus, the \emph{linearized Maxwell system}
\begin{equation}
  \left( \partial_t \begin{pmatrix} \epsilon(\partial_t) & 0 \\ 0 & \mu \end{pmatrix} + \begin{pmatrix} 0 & -\Curl \\ \Curl_0 & 0 \end{pmatrix} \right)
  \begin{pmatrix} E \\ H \end{pmatrix} =
  \begin{pmatrix} \varPhi \\ \varPsi \end{pmatrix}
  \label{eq:Mw1ord}
\end{equation}
is well-posed by Theorem~\ref{thm:linear_solution_theory}, since
  $\Re z \Bigl(\!
  \begin{smallmatrix}
    \epsilon(z) & 0 \\ 0 & \mu
  \end{smallmatrix}
  \Bigr)\ge \min\{c_\epsilon, c_\mu\}$
  for all $\varrho > \max\{1,\varrho_0\}$.

\begin{ex}
  \label{ex:LinearPermittivity}
  In most cases we will consider $\epsilon(\partial_t)E = \epsilon_0E + \chi\ast E$ and $\mu H = \mu_0H$, where $\epsilon_0,\mu_0 > 0$ are the vacuum permittivity and vacuum permeability, and $\chi\ast$ denotes the time convolution with the linear electric susceptibility tensor $\chi\colon \R \to \BLO(\Hs)$, with $\supp\chi \subseteq[0,\infty)$ due to causality, such that
  \begin{equation*}
    \chi\ast E = \int_0^\infty \chi(\mult - s) E(s) \dd s.
  \end{equation*}
  The simplest case in an interface setting is given by
  \begin{equation*}
    \chi(t) = \chi_1(t)1_{\Omega_1} + \chi_2(t)1_{\Omega_2}
  \end{equation*}
  with scalar-valued $\chi_1,\chi_2$ supported in $[0,\infty)$ and $1_{\Omega_i}$ being the characteristic function on $\Omega_i$.
  Let
  \begin{equation}
    \tilde{\chi}(\omega) = \frac{1}{\sqrt{2\pi}} \int_\R \chi(s)e^{i\omega s}\dd s
    \label{eq:FreqSusceptibility}
  \end{equation}
  denote the frequency-dependent susceptibility.
  Since
  \begin{equation*}
    (\mathcal{L}_\varrho \chi)(\xi)
    = \frac{1}{\sqrt{2\pi}} \int_{\R} \chi(s) e^{-(\varrho + i\xi)s} \dd s
    = \tilde{\chi}(i(\varrho + i\xi)),
  \end{equation*}
  we have $\epsilon(z) = \epsilon_0 + \tilde{\chi}(iz)$.
  If now $\chi_1,\chi_2 \in L^2_{\varrho_\chi}(\R)$ for some $\varrho_\chi \in \R$, then $\lvert \tilde{\chi}(iz) \rvert$ is bounded on $\C_{\Re > \varrho_0}$ for all $\varrho_0 > \varrho_\chi$, which follows from
  \begin{equation*}
    \lvert \tilde{\chi}(i(\varrho + i\xi)) \rvert^2
    = \Bigl| \int_\R \chi(s) e^{-(\varrho + i\xi)s}\dd s \Bigr|^2
    \le \lVert \chi \rVert_{L^2_{\varrho_\chi}}^2 \int_0^\infty e^{-2(\varrho - \varrho_\chi)s}\dd s
    \le \lVert \chi \rVert_{L^2_{\varrho_\chi}}^2 \int_0^\infty e^{-2(\varrho_0 - \varrho_\chi)s}\dd s
  \end{equation*}
  for $\varrho \ge \varrho_0$.
  Hence $\epsilon$ is a linear material law.
  The condition $\Re(z\epsilon(z)) \ge c_\epsilon \Re z \ge c_\epsilon \varrho_0$ above is satisfied for some $\varrho_0 > 0$ if $z \mapsto \lvert z\tilde{\chi}(iz)\rvert$ is bounded on $\C_{\Re > \varrho_1}$ for some $\varrho_1 \in \R$ (for example, this is the case for the Drude--Lorentz model in Appendix~\ref{sec:app1}).
  The extension to the case in which $\chi_1(t),\chi_2(t)\in L^\infty(\Omega)^{3\times 3}$ are matrix-valued is straightforward; here we impose the condition $\lVert \chi_1(\cdot) \rVert_{L^\infty}, \lVert \chi_2(\cdot) \rVert_{L^\infty} \in L^2_{\varrho_\chi}(\R)$.
\end{ex}

For the nonlinear part we similarly assume that
\begin{equation}
  P_\mathrm{nl}(E) = \int_{\R} \kappa(\mult-s)\, \snl(E(s)) \dd s,
  \label{eq:PolSimple}
\end{equation}
with $\kappa\colon \R \to \mathcal{B}(\Hs)$, $\supp \kappa \subseteq [0,\infty)$, and $\snl\colon \Hs\to\Hs$ Lipschitz-continuous.
Since $\partial_tP_\mathrm{nl}(E)$ is the only term on the right-hand side of \eqref{eq:Maxwell_nonlinear} depending on $E$ and
\begin{equation*}
  \partial_tP_\mathrm{nl}(E)(t) = \kappa(0^+)q(E(t)) + \int_0^\infty \kappa'(t-s)q(E(s))\dd s,
\end{equation*}
the needed Lipschitz-continuity in $E$ is implied by conditions on $\kappa$ and its derivative $\kappa'$. These conditions are provided next.

\begin{lem}
  \label{lem:Maxwell_nonlinear_AbstractLip}
  Let $\kappa\colon [0,\infty) \to \BLO(\Hs)$ be continuous and differentiable in $(0,\infty)$ and denote by $\kappa'\colon \R \to \BLO(\Hs)$ the zero extension of its derivative.
  Suppose $\lVert \kappa'(\mult) \rVert \coleq \lVert \kappa'(\mult) \rVert_{\BLO(\Hs)}$ is measurable, and let $\varrho_\kappa \in \R$ be such that $\lVert \kappa'(\mult) \rVert \in L^1_{\varrho_\kappa}(\R)$, \ie,
  \begin{equation}
    L_\kappa \coleq \int_\R \left\| \kappa'(s) \right\| e^{-\varrho_\kappa s}\dd s < \infty.
    \label{eq:kernel_decay}
  \end{equation}
  Furthermore, assume that $\snl$ is a Lipschitz-continuous map $\snl\colon \Hs \to \Hs$.
  Then
  \begin{equation*}
    Z\colon \Hw_\varrho(\R,\Hs) \to \Hw_\varrho(\R,\Hs), \quad Z(u)(t) = \int_\R \kappa'(t-s)\, \snl(u(s)) \dd s
  \end{equation*}
  is Lipschitz continuous, uniformly in $\varrho \ge \varrho_\kappa$, with $\lVert Z \rVert_{\Lip} \le L_\kappa\lVert q \rVert_\Lip$.
\end{lem}

\begin{proof}
  (cf.~\cite{MilaniPicard}) We compute for $u,v \in \Hw_{\varrho}(\R,\Hs)$,
  \begin{equation}
  \begin{aligned}
    \left\lVert Z(u) - Z(v) \right\rVert_{\varrho,0}^2
    &\le \left\lVert \int_\R \lVert \kappa'(t-s) \rVert \left| \snl(u(s)) - \snl(v(s)) \right|  \dd s \right\rVert_{\varrho,0}^2 e^{-2\varrho t}\dd t\\
    &\le \lVert \snl \rVert_{\Lip}^2 \int_\R \left( \int_\R \lVert \kappa'(t-s) \rVert {\lVert u(s) - v(s) \rVert}_{\Hs} \dd s \right)^2 e^{-2\varrho t}\dd t \\
    &\stackrel{(\star)}{\le} \lVert \snl \rVert_{\Lip}^2 L_\kappa \int_\R \biggl( \int_\R \lVert \kappa'(t-s) \rVert {\lVert u(s) - v(s) \rVert}_{\Hs}^2\, e^{\varrho_\kappa (t-s)} \dd s \biggr)\, e^{-2\varrho t}\dd t \\
    &= \lVert \snl \rVert_{\Lip}^2 L_\kappa \int_\R \int_\R \lVert \kappa'(t-s) \rVert e^{-\varrho_\kappa (t-s)} \cdot \\
    &\qquad \cdot \underbrace{e^{-2(\varrho - \varrho_\kappa) (t-s)}}_{\le 1\, \text{for } t-s \ge 0}\, \dd t\,
    {\lVert u(s) - v(s) \rVert}_{\Hs}^2\, e^{-2\varrho s}\, \dd s \\
    &\le \lVert \snl \rVert_{\Lip}^2 L_\kappa \int_\R \lVert \kappa'(r) \rVert e^{-\varrho_\kappa r} \dd r \int_\R {\lVert u(s) - v(s) \rVert}_{\Hs}^2\, e^{-2\varrho s}\, \dd s \\
    &= \lVert \snl \rVert_{\Lip}^2 L_\kappa^2\, \lVert u - v \rVert_{\varrho,0}^2,
  \end{aligned}
  \label{eq:convolution_kappa_Lipschitz}
  \end{equation}
  where $(\star)$ follows after writing $\lVert \kappa'(t-s) \rVert = {\lVert \kappa'(t-s) \rVert}^{\frac{1}{2} + \frac{1}{2}} e^{-\varrho_\kappa(t-s)(\frac{1}{2} - \frac{1}{2})}$ and applying the Cauchy--Schwarz inequality.
\end{proof}

Lemma~\ref{lem:Maxwell_nonlinear_AbstractLip} yields the uniform Lipschitz-continuity of
\begin{equation*}
  E\mapsto \partial_tP_\mathrm{nl}(E) = \kappa(0^+)\snl(E(\cdot)) + Z(E),
\end{equation*}
implying for $\varrho > \max\{ \varrho_\kappa, 0\}$ that the solution operator
\begin{equation*}
  S_\varrho^\mathrm{nl}\colon
  \begin{pmatrix}
    E \\ H
  \end{pmatrix} \mapsto
  \left(
      \overline{
        \partial_t
        \begin{pmatrix}
          \epsilon(\partial_t) & 0 \\ 0 & \mu
        \end{pmatrix} +
        \begin{pmatrix}
          0 & -\Curl \\ \Curl_0 & 0
        \end{pmatrix}
      }
      \right)^{-1}
      \left(
      \begin{pmatrix}
        -\partial_tP_\mathrm{nl}(E) \\ 0
      \end{pmatrix} +
      \begin{pmatrix}
        \varPhi \\ \varPsi
      \end{pmatrix}
      \right)
\end{equation*}
is causal and Lipschitz-continuous for any $\varPhi,\varPsi \in \Hw_{\varrho}(\R,\Hs)$, with Lipschitz constant at most
\begin{equation}
  \lVert S_{\varrho}^\mathrm{nl} \rVert_{\Lip}
  \le \frac{c}{\varrho} {\lVert \snl \rVert}_{\Lip} (|\kappa(0^+)| + L_\kappa),
  \label{eq:sol_F_Lip}
\end{equation}
where $c$ is given by the condition \eqref{eq:Lipschitz_perturb} imposed on $\epsilon(\cdot)$ and $\mu$.
Choosing $\varrho \ge \varrho_\kappa$ large enough, $S_{\varrho}^\mathrm{nl}$ becomes a contraction on $\Hw_\varrho(\R,\Hs\times\Hs)$.

\begin{ex}[Saturable nonlinearity]
  \label{ex:SaturableNonlinearity1}
  Let $k\in\N_{\ge 2}$ and $\tau > 0$ and consider $\snl\colon \Hs \to \Hs$ given by
  \begin{equation*}
    \snl(u)(x) = \frac{|u(x)|^{k-1}}{1 + \tau|u(x)|^{k-1}}u(x) \eqcol V(|u(x)|)u(x).
  \end{equation*}
  (For $k = 3$ this is a saturable version of the Kerr-type nonlinearity $E \mapsto |E|^2E$.)
  Since $\R^3\ni \xi\mapsto V(|\xi|)\xi$ is smooth and asymptotically linear, it is Lipschitz-continuous, hence $q\colon \Hs \to \Hs$ is Lipschitz-continuous.
  Thus, $P_\mathrm{nl}$ defined as in \eqref{eq:PolSimple} with $\kappa$ as in Lemma \ref{lem:Maxwell_nonlinear_AbstractLip}, fulfills the necessary assumptions of the lemma.
\end{ex}

\subsubsection{Local well-posedness}
Let $X$ be a general Hilbert space.
The uniform Lipschitz-continuity in the range of spaces $\Hw_\varrho(\R,\HX)$ imposed on the nonlinearity may seem restrictive (in particular, nonlinearities growing at a superlinear rate are excluded as candidates for $q$). In fact, this condition can be replaced by Lipschitz continuity on closed subsets in $\Hw_\varrho(\R,\HX)$, which eventually (for large $\varrho$) grow large enough to include given data.
To illustrate this, we formulate the following refinement of Proposition \ref{prop:Lipschitz_perturb}.

\begin{prop}
  \label{prop:FixPointRefinement}
  Let $A: \dom(A) \subset \HX \to \HX$ be skew-selfadjoint and $M$ a linear material law with $\Re zM(z) \ge \Re z/d$ for $\Re z > \varrho_0$.
  Let $f\colon L^2_{\varrho}(\R, \HX) \to L^2_{\varrho}(\R, \HX)$ be a causal nonlinear map satisfying $f(0) = 0$ and let $c, \alpha > 0$ be such that for $\varrho > \varrho_0$
  \begin{equation}
    \lVert f(u) - f(v) \rVert_{\varrho,0} \le c \bigl( \lVert u \rVert_{\varrho,0} + \lVert v \rVert_{\varrho,0} \bigr)^\alpha \lVert u - v \rVert_{\varrho,0}
    \label{eq:LocLip_Refinement}
  \end{equation}
  for all $u,v \in L^2_\varrho(\R,\HX)$.
  Suppose $g \in L^2_{\varrho}(\R,\HX)$ is such that $\lVert g \rVert_{\varrho,0} = o(\varrho^{1+\frac{1}{\alpha}})$ as $\varrho \to \infty$.
  Then the equation $(\partial_tM(\partial_t) + A)u = f(u) + g$ admits a unique solution $u \in L^2_\varrho(\R,\HX)$ for large $\varrho > \varrho_0$.
\end{prop}

\begin{proof}
  Denote by
  $S_\varrho \coleq \bigl(\overline{\partial_tM(\partial_t) + A}\bigr)^{-1}\colon \Hw_\varrho(\R,\HX)\to \Hw_\varrho(\R,\HX)$
  the linear solution operator.
  Then, using $\lVert S_\varrho \rVert \le d/\varrho$, the Lipschitz constant of $S_\varrho(f(\mult) + g)$ on a closed ball $B_r := \{ u \in L^2_{\varrho}(\R,\Hs) : \lVert u \rVert_{\varrho,0} \le r \}$ for $r > 0$ can be estimated by
  \begin{equation*}
    L_{\varrho,r} := \sup_{u,v \in B_r, u\ne v}\frac{\lVert S_\varrho f(u) - S_\varrho f(v) \rVert_{\varrho,0}}{\lVert u - v \rVert_{\varrho,0}}
    \le \frac{cd}{\varrho} \left( 2r \right)^{\alpha}
  \end{equation*}
  for $\varrho > \varrho_0$, thus $L_{\varrho,r} < 1$ if $r < \frac{1}{2} \bigl(\frac{\varrho}{cd}\bigr)^{1/\alpha}$.
  In order to have $S_\varrho(f(u) + g) \in B_r$ for all $u \in B_r$, we demand that
  \begin{align*}
    \lVert S_\varrho(f(u) + g) \rVert_{\varrho,0}
    \le \frac{d}{\varrho} \bigl(c\lVert u \rVert_{\varrho,0}^{\alpha + 1} + \lVert g \rVert_{\varrho,0}\bigr)
    \le \frac{d}{\varrho}(cr^{\alpha+1} + \lVert g \rVert_{\varrho,0}) \stackrel{!}{\le} r.
  \end{align*}
  Replacing $r$ with $\frac{1}{2}\bigl(\frac{\varrho}{cd}\bigr)^{1/\alpha}$ in the last inequality leads to the condition
  \begin{equation*}
    \frac{\lVert g \rVert_{\varrho,0}}{c}
    \stackrel{!}{<} \frac{1}{2}\frac{\varrho}{cd} \Bigl( \frac{\varrho}{cd}\Bigr)^{\frac{1}{\alpha}} - \Bigl( \frac{1}{2}\Bigl( \frac{\varrho}{cd}\Bigr)^{\frac{1}{\alpha}} \Bigr)^{\alpha + 1}
    = \frac{1}{2}\Bigl( \frac{\varrho}{cd}\Bigr)^{1 + \frac{1}{\alpha}} \bigl( 1 - 2^{-\alpha} \bigr),
  \end{equation*}
  which is fulfilled by assumption on $g$ for large $\varrho > 0$. This establishes $S_\varrho(f(\mult) + g)$ as a contraction on $B_r$ for some $r < \frac{1}{2}\bigl(\frac{\varrho}{cd}\bigr)^{1/\alpha}$.
\end{proof}

\begin{rem}
  The condition $\lVert g \rVert_{\varrho,0} = o(\varrho^{1 + \frac{1}{\alpha}})$ as $\varrho \to \infty$ of Proposition~\ref{prop:FixPointRefinement} is satisfied if $g \in \Hw_{\varrho}(\R,\Hs)$ for some $\varrho \in \R$ with $\supp g \in [0,\infty)$. This latter assumption, in fact $\lVert g \rVert_{\varrho,0} = O(1)$, is justified in Section~\ref{sec:initial_values}.
\end{rem}

Consider now a quadratic nonlinearity of the form
\begin{equation}
  f(u)(t) = \int_\R \int_\R K(t-\tau_1, t-\tau_2)
  q(u(\tau_1), u(\tau_2)) \dd{\tau_1}\dd{\tau_2},
  \label{eq:QuadraticNonlinearity}
\end{equation}
where $K\colon \R^2 \to \BLO(\HX)$ is an operator-valued kernel with $\supp K \subseteq [0,\infty)^2$ (to ensure causality), and where $q\colon \HX\times\HX \to \HX$ is a bounded bilinear map, \ie,
  $\lVert q(u,v) \rVert_{\HX} \le C_q \lVert u \rVert_{\HX} \lVert v \rVert_{\HX}$ for some $C_q > 0$.
In analogy to (2.12) we impose the following integrability conditions,
\begin{equation}
  \left.
  \begin{aligned}
    L_K &\coleq \iint \lVert K (\tau_1,\tau_2)\rVert\, e^{-\varrho_K(\tau_1 + \tau_2)} \dd{\tau_1}\dd{\tau_2} < \infty\\
    \ell_K &\coleq \sup_{\tau_1,\tau_2\in \R}
    \int \left\|K (t-\tau_1,t-\tau_2)\right\|
    e^{-\varrho_K(2t-\tau_1-\tau_2)} \dd{\tau_1}\dd{\tau_2} < \infty
  \end{aligned}
  \right\}
  \label{eq:KernelIntegrableQuadratic}
\end{equation}
for some $\varrho_K \in \R$.
Using the same strategy as in (2.13) we can show for all $\varrho \ge \varrho_K$ that $f$ maps $\Hw_{\varrho}(\R,\HX)$ into $\Hw_{2\varrho}(\R,\HX)$;
indeed,
\begin{multline*}
  \int_\R \left\lVert \iint K(t-\tau_1, t-\tau_2)
  q(u(\tau_1), v(\tau_2)) \dd{\tau_1}\dd{\tau_2} \right\rVert_\HX^2 e^{-4\varrho t} \dd t\\
  \begin{aligned}
    &\le L_KC_q^2 \int_\R \left( \iint \lVert K(t-\tau_1, t-\tau_2) \rVert e^{\varrho_K(2t-\tau_1-\tau_2)} \lVert u(\tau_1) \rVert_\HX^2 \lVert v(\tau_2) \rVert_\HX^2 \dd{\tau_1} \dd{\tau_2}\right) e^{-4\varrho t} \dd t\\
    &\le L_K C_q^2\iint \biggl( \int_\R \lVert K(t-\tau_1, t-\tau_2) \rVert e^{-\varrho_K(2t-\tau_1-\tau_2)} \underbrace{e^{-2(\varrho-\varrho_K)(2t-\tau_1-\tau_2)}}_{\le 1 \text{ for $\tau_1,\tau_2\le t$}} \dd t \biggr) \cdot\\
    &\qquad \cdot \lVert u(\tau_1) \rVert_\HX^2 e^{-2\varrho\tau_1} \dd{\tau_1} \lVert v(\tau_2) \rVert_\HX^2 e^{-2\varrho\tau_2} \dd{\tau_2}\\
    &\le L_K\ell_K C_q^2\lVert u \rVert_{\varrho,0}^2 \lVert v \rVert_{\varrho,0}^2.
  \end{aligned}
\end{multline*}
This computation makes it clear, however, that the mapping property
$f\colon \Hw_\varrho(\R,\HX) \to \Hw_\varrho(\R,\HX)$ cannot be obtained in general for $\varrho > 0$. (Yet if $\varrho_K < 0$, then $f$ leaves a subspace of $L^2_{\varrho}(\R,\HX)$, $\varrho_K \le \varrho < 0$, invariant; this fact is relevant in combination with the notion of exponential stability, see Section~\ref{sec:exp_stab}).
\bigskip

\begin{rem}[Nonlinearity with a cutoff in time]
  We now show how the fixed-point argument underlying Proposition~\ref{prop:FixPointRefinement} can still be applied to a modified version of the nonlinearity above to obtain a local well-posedness result for $\varrho \ge \varrho_K > 0$.
  In detail, we modify the kernel $K$ by applying a cutoff in the $t$-variable. 
  Suppose again $K\colon \R^2 \to \BLO(\HX)$ is causal map satisfying \eqref{eq:KernelIntegrableQuadratic}
  and that, in addition,
  \begin{equation}
    d_{K} \coleq \mathop{\mathrm{ess\,sup}}_{\tau_1,\tau_2 \ge 0}\left\lVert K(\tau_1,\tau_2)\right\rVert e^{-\varrho_K(\tau_1 + \tau_2)} < \infty.
    \label{eq:KernelBoundedQuadratic}
  \end{equation}
  For $T > 0$ we define $K_T\colon \R^3 \to \BLO(\HX)$ and $f_T\colon \Hw_\varrho(\R,X) \to \Hw_\varrho(\R,X)$ by
  \begin{align*}
    K_T (t,\tau_1,\tau_2) &\coleq 1_{(-\infty,T]}(t)K(\tau_1,\tau_2)\\
    f_T(u)(t) &\coleq
    \int_\R \int_\R K_T(t,t-\tau_1, t-\tau_2)
      q(u(\tau_1), u(\tau_2)) \dd{\tau_1}\dd{\tau_2},
  \end{align*}
  and observe that
  \begin{align*}
    \int_\R\int_\R \left\lVert K_T(t,\tau_1,\tau_2) \right\rVert e^{-\varrho_K(\tau_1 + \tau_2)}\dd{\tau_1}\dd{\tau_2}
    &\le L_K\\
    \int_\R \left\lVert K_T(t,\tau_1,\tau_2)\right\rVert e^{-\varrho_K(\tau_1 + \tau_2)}e^{2\varrho t} \dd t
    \le \int_0^T d_{K} e^{2\varrho t} \dd t
    &\le Te^{2\varrho T}\,d_{K},\quad
    \text{if } \tau_1,\tau_2 \ge 0
  \end{align*}
  for all $\varrho, t, T > 0$.
  Now modifying the estimate above we obtain for $\varrho \ge \varrho_K$ and $u,v \in \Hw_\varrho(\R,\HX)$
  \begin{multline*}
    \int_\R \left\lVert \iint K_T(t,t-\tau_1, t-\tau_2)\,
    q(u(\tau_1), v(\tau_2)) \dd{\tau_1}\dd{\tau_2} \right\rVert_\HX^2 e^{-2\varrho t} \dd t\\
    \begin{aligned}
      &\le L_KC_q^2 \int_0^T \left( \iint \left\lVert K(t-\tau_1, t-\tau_2) \right\rVert e^{\varrho_K(2t-\tau_1-\tau_2)} \lVert u(\tau_1) \rVert_\HX^2 \lVert v(\tau_2) \rVert_\HX^2 \dd{\tau_1} \dd{\tau_2}\right) e^{-2\varrho t} \dd t\\
      &\le L_KC_q^2 \iint \biggl( \int_0^T \lVert K(t-\tau_1, t-\tau_2) \rVert e^{-\varrho_K(2t-\tau_1-\tau_2)} 
      \underbrace{e^{2\varrho_K(2t-\tau_1-\tau_2) + 2\varrho(\tau_1 + \tau_2 - t)}}_{\le e^{2\varrho t} \text{ for $\tau_1,\tau_2\le t$}} \dd t \biggr) \cdot\\
      &\qquad \cdot \lVert u(\tau_1) \rVert_\HX^2 e^{-2\varrho\tau_1} \dd{\tau_1} \lVert v(\tau_2) \rVert_\HX^2 e^{-2\varrho\tau_2} \dd{\tau_2}\\
      &\le Te^{2\varrho T}d_{K} L_KC_q^2 \lVert u \rVert_{\varrho,0}^2 \lVert v \rVert_{\varrho,0}^2.
    \end{aligned}
  \end{multline*}
  Since by bilinearity of $q$ we have
  \begin{multline}
    \lVert q(u(\tau_1), u(\tau_2)) - q(v(\tau_1), v(\tau_2)) \rVert_{\HX}\\
    \begin{aligned}
      &\le \lVert q(u(\tau_1), u(\tau_2) -  v(\tau_2)) \rVert_{\HX}
      + \lVert q(u(\tau_1) - v(\tau_1), v(\tau_2)) \rVert_{\HX}\\
      &\le C_q\left(\lVert u(\tau_1) \rVert_{\HX} \lVert u(\tau_2) - v(\tau_2) \rVert_{\HX}
      + \lVert u(\tau_1) - v(\tau_1) \rVert_{\HX} \lVert v(\tau_2) \rVert_{\HX}\right),
    \end{aligned}
    \label{eq:spatialLocLip}
  \end{multline}
  the estimate above produces
  \begin{equation}
    \lVert f_T(u) - f_T(v) \rVert_{\varrho,0}
    \le \sqrt{T}e^{\varrho T}C_q \sqrt{d_KL_K}
    \left( \lVert u \rVert_{\varrho,0} + \lVert v \rVert_{\varrho,0} \right)
    \lVert u - v \rVert_{\varrho,0}.
    \label{eq:NonlinearityLocLip_T}
  \end{equation}
  Hence, given $\varrho \ge \varrho_K$, the parameters $T, r > 0$ can be chosen small enough so that $u\mapsto S_\varrho (f_T(u) + g)$ becomes a contraction on a closed ball with radius $r$ in $\Hw_\varrho(\R,X)$, provided that the data $g \in \Hw_\varrho(\R,\HX)$ is small enough.
\end{rem}
  
\begin{ex}
  \label{ex:NonlocalNonlinearity}
  As an application to the Maxwell system, let $P_\mathrm{nl} = P^{(2)}$, where $P^{(2)}$ is a fully nonlocal quadratic polarization given by
  \begin{equation*}
    P^{(2)}(E)(t) \coleq 
    \int_\R \int_\R \kappa(t-\tau_1,t-\tau_2)\, q(E(\tau_1), E(\tau_2))\dd{\tau_1}\dd{\tau_2}.
  \end{equation*}
  Here we assume that the spatial map $q$ is defined via a tensor $\Lambda = (\Lambda_{ijk})_{i,j,k\in\{1,2,3\}}$ with $\Lambda_{ijk}\in L^2(\Omega^3)$,
  \begin{align*}
    q(u,v)(x) &= \int_{\Omega}\int_{\Omega} \Lambda(x,y,y')u(y)v(y')\dd{y}\dd{y'}\\
      &\coleq \biggl(\int_\Omega\int_\Omega \sum_{j,k \in \{1,2,3\}}\Lambda_{ijk}(x,y,y') u_j(y) v_k(y')\dd{y}\dd{y'}\biggr)_{i=1,2,3},
  \end{align*}
  and that $\kappa \in C^1(\R^2,\R^{3\times 3})$ is a matrix-valued map with compact support in $[0,\infty)^2$.
  By the Cauchy--Schwarz inequality we have the pointwise estimate
  \begin{align*}
    \Bigl\lvert \int_\Omega\int_\Omega \Lambda_{ijk}(x,y,y')u_j(y)v_k(y') \dd{y}\dd{y'} \Bigr\rvert
    \le \lVert \Lambda_{ijk}(x,\mult,\mult) \rVert_{L^2(\Omega^2)} \lVert u_j \rVert_{L^2(\Omega)} \lVert v_k \rVert_{L^2(\Omega)},
  \end{align*}
  from which we obtain $\lVert q(u,v) \rVert_{\Hs} \le C_q \lVert u \rVert_\Hs \lVert v \rVert_\Hs$ with an appropriate constant $C_q$. Hence, $q$ is a bilinear, bounded map and satisfies \eqref{eq:spatialLocLip} with $X = L^2(\Omega)^3$.
  Furthermore, formally computing the derivative of $P^{(2)}(E)$ gives
  \begin{equation*}
    \begin{aligned}
      \partial_tP^{(2)}(E)(t)
      &= \int_{-\infty}^t \kappa(t-\tau_1,0)\,\snl(E(\tau_1),E(t))\dd{\tau_1}
      + \int_{-\infty}^t \kappa(0,t-\tau_2)\,\snl(E(t),E(\tau_2))\dd{\tau_2}\\
      &\quad+ \int_{-\infty}^t\int_{-\infty}^t (\partial_1 + \partial_2)\kappa(t-\tau_1,t-\tau_2)\,\snl(E(\tau_1),E(\tau_2))\dd{\tau_1}\dd{\tau_2}\\
      &= \int_{-\infty}^t\int_{-\infty}^t (\partial_1 + \partial_2)\kappa(t-\tau_1,t-\tau_2)\,\snl(E(\tau_1),E(\tau_2))\dd{\tau_1}\dd{\tau_2}
    \end{aligned}
  \end{equation*}
  since $\kappa(0,\mult) = \kappa(\mult,0) = 0$, \ie, $\partial_tP^{(2)}$ is again a fully nonlocal map of the form \eqref{eq:QuadraticNonlinearity}, where $K \coleq (\partial_1 + \partial_2)\kappa$ is continuous with compact support in $[0,\infty)$.
  As such, $K$ satisfies \eqref{eq:KernelIntegrableQuadratic} and \eqref{eq:KernelBoundedQuadratic} for arbitrary $\varrho_K \in \R$.
  Thus, the cutoff version of $\partial_tP^{(2)}$ defined by
  \begin{equation*}
    (\partial_tP^{(2)})_T \coleq 1_{(-\infty,T]}(\partial_tP^{(2)})
  \end{equation*}
  satisfies \eqref{eq:NonlinearityLocLip_T}.
  \medskip

  The same principle applies to multilinear maps in general:
  Let $n \in \N_{\ge 2}$ and $q\colon (\Hs)^{n} \to \Hs$ be a bounded $n$-linear map. Let $\kappa\colon \R^n \to \R^{n\times n}$ be supported in $[0,\infty)^n$ with $\kappa(s_1,\ldots,s_n) = 0$ whenever $s_j = 0$ for some $j\in \{1,\ldots,n\}$.
  Defining $P^{(n)}$ as in \eqref{eq:Polarization} by
  \begin{equation*}
    P^{(n)}(E)(t) = \int\cdots\int \kappa(t-\tau_1,\ldots,t-\tau_n) q(E(\tau_1),\ldots,E(\tau_n))\dd{\tau_1}\cdots\dd{\tau_n},
  \end{equation*}
  we find that, if $\kappa$ satisfies an integrability condition similar to \eqref{eq:KernelIntegrableQuadratic}, then $(\partial_tP^{(n)})_T = 1_{(-\infty,T]}\partial_tP^{(n)}$, $T > 0$, satisfies
  \begin{equation}
    \lVert (\partial_tP^{(n)})_T(u) - (\partial_tP^{(n)})_T(v) \rVert_{\varrho,0}
    \le \sqrt{T} e^{(n-1)\varrho T}C
    \bigl( \lVert u \rVert_{\varrho,0} + \lVert v \rVert_{\varrho,0} \bigr)^{n-1}
    \lVert u - v \rVert_{\varrho,0}.
    \label{eq:MultilinearLocLip_T}
  \end{equation}
  Reasoning as above in (ii) we obtain the following result.

  \begin{prop}[Well-posedness of the Maxwell system with fully nonlocal multilinear polarization and a cutoff in time]
    \label{prop:LocalWP_multilinear}
    Suppose the linear system \eqref{eq:Mw1ord} is well-posed in $\Hw_\varrho(\R,\Hs)$, $\varrho > \varrho_0$. Then, for each $\varrho > \max\{\varrho_0, \varrho_K\}$, the nonlinear system
    \begin{equation*}
      \left(
      \partial_t
      \begin{pmatrix}
        \epsilon(\partial_t) & 0 \\ 0 & \mu
      \end{pmatrix} + 
      \begin{pmatrix}
        0 & \Curl \\ \Curl_0 & 0
      \end{pmatrix} \right)
      \begin{pmatrix} E \\ H \end{pmatrix} =
      \begin{pmatrix} -(\partial_tP^{(n)})_T(E) \\ 0 \end{pmatrix} +
      \begin{pmatrix} \varPhi \\ \varPsi \end{pmatrix}
    \end{equation*}
    admits a unique solution $(E,H)\in \Hw_{\varrho}(\R,\Hs)^2$ for small $T > 0$ and small data $\varPhi,\varPsi \in \Hw_\varrho(\R,\Hs)$.
  \end{prop}
\end{ex}

\subsection{Initial values}
\label{sec:initial_values}

In order to apply the well-posedness theory to system \eqref{eq:Maxwell_pos_init}, it remains to discuss how the history of the electromagnetic field can be incorporated into the framework.
First we mention the following result concerning regularity in time (see \cite[Section 3.1]{PicardMcGhee} or \cite[Section 6.3]{STW22}).
(We subsequently use the notation $f \lesssim g$ or equivalently $g \gtrsim f$ to denote $f \le Cg$ for some $C > 0$ independent of $f,g$.)

\begin{prop}
  \label{prop:time_regularity}
  Let $\left(\partial_tM(\partial_t) + A \right)u = g$ be well-posed in the range of spaces $\Hw_{\varrho}(\R,\HX)$ for $\varrho > \varrho_0$ with $\varrho_0 \in \R$.
  If $g \in H^1_{\varrho}(\R,\HX)$, then $u = \bigl(\overline{\partial_tM(\partial_t) + A}\bigr)^{-1} g \in H^1_{\varrho}(\R,\HX)$, with continuous dependence on the data:
  \begin{equation*}
    \lVert u \rVert_{\varrho,0} \lesssim \lVert g \rVert_{\varrho,0},\quad
    \lVert \partial_tu \rVert_{\varrho,0} \lesssim \lVert \partial_tg \rVert_{\varrho,0}.
  \end{equation*}
  In fact, $\partial_tu = {(\overline{\partial_tM(\partial_t) + A})}^{-1}\partial_tg$.
  Moreover,
  \begin{equation*}
    u \in C_\varrho(\R,\HX) \coleq \{ f \in C(\R,\HX) : \sup_{t \in \R} \left\| f(t)\right\|_Xe^{-\varrho t} < \infty \}
  \end{equation*}
  by the Sobolev embedding theorem.
\end{prop}

\begin{figure}[htbp]
  \centering
    \begin{tikzpicture}
      \draw[->] (-5,0) -- (5.1,0) node[below] {$t$};
      \draw[help lines] (0,0) -- (0,3);
      \draw[green!40!blue, very thick,rounded corners] (-5,.8) -- ++(.5,-.1) -- ++(.5,-.2) -- ++(1,.2) -- ++(1,.5) node[above left] {$\phi$} -- ++(1,.7) -- ++(.7,.3) -- ++(.3,0) node[behind path] (a) {};
      \draw[green!40!blue, very thick,rounded corners] (0,0) -- (5,0);
      \draw[magenta, very thick,rounded corners] (a.center) -- ++(.3,0) node[above right] {$u$} -- ++(.7,-.1) -- ++(1,-.5)  node (b) {} -- ++(1,-.3) -- ++(1,-.7) -- ++(.7,-.3) -- ++(.3,0);
      \draw[magenta, very thick,rounded corners] (-5,0) -- (0,0);
      \draw[black!80!white, ultra thick, dash dot, rounded corners] (-5,0) -- (.1,0) -- ++(.2,.2) -- ++($(a.center) + (.4,-1)$) node[below right] {$\tilde{u}$} -- ++(.2,.3) -- ++(.2,.2) -- ++(.3,.0) -- (b.center) -- ++(1,-.3) -- ++(1,-.7) -- ++(.7,-.3) -- ++(.3,0);
      \draw [very thick,gray] (0,-.15) -- (0,-.3) -- node[below] {$\supp \theta^+\eta$} (1.5,-.3) --++(0,.15);
    \end{tikzpicture}
  \caption{Schematic for the conversion of the Cauchy problem to an evolutionary equation.}
  \label{fig:Cauchy2Evo}
\end{figure}
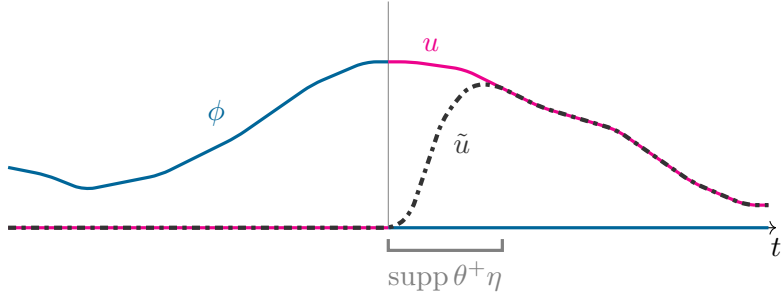

Consider a general Cauchy problem for the Maxwell equations,
\begin{equation}
  \left\{
    \begin{aligned}
      \partial_t\mathcal{M}(U)(t) + \A U(t) &= 0, & t > 0 \\
      U(t) &= \phi(t), & t \le 0
    \end{aligned}
  \right\}
  \label{eq:Maxwell_matrix_pos}
\end{equation}
for a given history $\phi\colon (-\infty,0] \to L^2(\Omega)^6$.
For simplicity of this model problem we have set $J = 0$ and assume
\begin{equation*}
  \mathcal{M}(U) = M_0U + \mathcal{G}(U),\quad \text{with}\quad \mathcal{G}(U) = \chi\ast\snl(U),
\end{equation*}
where $M_0$ is selfadjoint and uniformly positive definite, $\chi$ is rapidly decaying, smooth, and $\supp \chi \subseteq (0,\infty)$, and $\snl\colon L^2(\Omega)^6\to L^2(\Omega)^6$ is Lipschitz-continuous with $q(0) = 0$.
We want to convert \eqref{eq:Maxwell_matrix_pos} into a nonlinear evolutionary equation 
in $L^2_\varrho(\R,L^2(\Omega)^6)$ (we note however that the derivation below is not strictly tied to the Maxwell system).
To this end, suppose $U\in C(\R,L^2(\Omega)^6)$ is a continuous solution of \eqref{eq:Maxwell_matrix_pos}.
Let $\theta^+ \coleq \theta^+_0$ denote multiplication with the Heaviside step function, then the projection
\begin{equation*}
  u \coleq \theta^+U
\end{equation*}
separates the “unknown” solution $u$ with $\supp u \subseteq [0,\infty)$ from the given history $\phi$, which we extend trivially to the whole line, thus $ \phi = (1 - \theta^+)\phi$.
With $U = u + \phi$ we also have $\snl(U(t)) = \snl(u(t)) + \snl(\phi(t))$ for all $t \in \R$, and therefore in fact $\mathcal{M}(U)(t) = \mathcal{M}(u)(t) + \mathcal{M}(\phi)(t)$.
Interpreting now $\partial_t$ in the distributional sense, we use the formula
\begin{equation*}
  \partial_t(\theta^+\varphi) = \theta^+\partial_t\varphi + \varphi(0^+)\delta_0
\end{equation*}
to extract from \eqref{eq:Maxwell_matrix_pos} an equation for $u$ on the whole real line:
\begin{align}
  0 = \theta^+\bigl[ \partial_t\mathcal{M}(U) + \A U \bigr]
  &= \partial_t(\theta^+\mathcal{M}(U)) - \mathcal{M}(U)(0^+)\delta_0 + \A \theta^+ U \nonumber\\
  &= \partial_t(\theta^+\mathcal{M}(u)) + \partial_t(\theta^+\mathcal{M}(\phi)) - \mathcal{M}(\phi)(0^-)\delta_0 + \A u \nonumber\\
  &= \partial_t\mathcal{M}(u)  + \A u + \partial_t(\theta^+\mathcal{G}(\phi)) - M_0\phi(0^-)\delta_0 - \mathcal{G}(\phi)(0^-)\delta_0 \nonumber\\
  &= \partial_t\mathcal{M}(u)  + \A u + \theta^+\partial_t\mathcal{G}(\phi)  + \mathcal{G}(\phi)(0^+)\delta_0 - M_0\phi(0^-)\delta_0 - \mathcal{G}(\phi)(0^-)\delta_0 \nonumber\\
  &= \partial_t\mathcal{M}(u)  + \A u + \theta^+\partial_t\mathcal{G}(\phi) - M_0\phi(0^-)\delta_0, \label{eq:evo_line_with_jump}
\end{align}
where we used $\mathcal{G}(\phi)(0^-) = \mathcal{G}(\phi)(0^+)$.
The $\delta_0$-term can be removed by smoothing the jump of $u$ at $t = 0$:
Choose $\eta \in C_c^\infty(\R)$ with $\eta(0) = 1$, and set
\begin{equation*}
  \phi^+ := \phi(0^-)\theta^+\eta,\quad \tilde{u} \coleq u - \phi^+,
\end{equation*}
see Figure~\ref{fig:Cauchy2Evo}. Then,
\begin{align*}
  \partial_t\mathcal{M}(u)
  = \partial_t\mathcal{M}(\tilde{u} + \phi^+)
  &= \partial_t\bigl( M_0\tilde{u} + M_0\phi^+ + \mathcal{G}(\tilde{u} + \phi^+) \bigr)\\
  &= \partial_t \bigl(M_0\tilde{u} + \mathcal{G}(\tilde{u} + \phi^+)\bigr) + \theta^+\partial_tM_0\phi^+ + M_0\phi^+(0^+)\delta_0.
\end{align*}
Thus, using that $\phi^+(0^+) = \phi(0^-)$, \eqref{eq:evo_line_with_jump} becomes
\begin{align*}
  0 &= \partial_t\mathcal{M}(\tilde{u} + \phi^+) + \A \tilde{u} + \A \phi^+ + \theta^+\partial_t\mathcal{G}(\phi) - M_0\phi(0^-)\delta_0\\
  &= \partial_t \bigl(M_0\tilde{u} + \mathcal{G}(\tilde{u} + \phi^+)\bigr) + \A \tilde{u} + \theta^+\partial_tM_0\phi^+ + \theta^+\partial_t\mathcal{G}(\phi) + \A \phi^+.
\end{align*}
Finally, the last identity can be written as
\begin{equation}
  \bigl(\partial_t M_0 + \A\bigr)\tilde{u} = -\partial_t\mathcal{G}(\tilde{u} + \phi^+) + g_\phi
  \label{eq:evo_line_proper_reformulation}
\end{equation}
where
\begin{equation*}
  g_\phi \coleq -\theta^+\bigl[ \partial_t\bigl( M_0\phi^+ + \mathcal{G}(\phi)\bigr) + \A \phi^+ \bigr].
\end{equation*}
Now \eqref{eq:evo_line_proper_reformulation} is a proper reformulation of \eqref{eq:Maxwell_matrix_pos} as an operator equation in $\Hw_\varrho(\R,L^2(\Omega)^6)$.
The well-posedness follows by Proposition \ref{prop:Lipschitz_perturb} from the Lipschitz continuity of $\tilde{u}\mapsto\partial_t\mathcal{G}(\tilde{u} + \phi^+)$.
Since $\phi^+ = 0$ on $(-\infty,0]$, the causality of the solution operator and the fixed-point iteration implies $\tilde{u} = 0$ on $(-\infty,0]$.

\begin{rem}
  The initial value theory in \cite{Trostorff_habil} for linear systems uses $\eta \equiv 1$, however, the present choice $\eta \in C_c^\infty(\R)$ is more convenient since it also works in the context of exponential stability, \ie, if the system is well-posed for $\varrho = -\nu < 0$, we have $\tilde{u} \in \Hw_{-\nu}(\R,\Hs)$ if and only if $u \in \Hw_{-\nu}(\R,\Hs)$; see also \cite[Section 4]{TrostorffLinExpStab}.
\end{rem}

\begin{rem}[A posteriori justification]
  If $g_\phi \in H^1_{\varrho}(\R,L^2(\Omega)^6)$, then solutions of \eqref{eq:evo_line_proper_reformulation} generate continuous solutions of \eqref{eq:Maxwell_matrix_pos}.
  Indeed, in this case Proposition \ref{prop:time_regularity} justifies $\tilde{u} \in H^1_{\varrho}(\R,L^2(\Omega)^6)$, and since $\phi - \phi^+$ is continuous, $U = \tilde{u} + (\phi - \phi^+)  \in C(\R,L^2(\Omega)^6)$.
Assuming the history $\phi$ is sufficiently regular, $\phi \in H^1_\varrho((-\infty,0], L^2(\Omega)^6)$ with $\phi(0^-) \in \dom(\A)$, then $g_\phi \in H^1_\varrho((0,\infty),L^2(\Omega)^6)$.
In this case, since $g_\phi = 0$ on $(-\infty,0]$, a necessary and sufficient condition for $g_\phi \in H^1_\varrho(\R,L^2(\Omega)^6)$ is the continuity of $g_\phi$ in $t = 0$, \ie,
\begin{equation}
  g_\phi(0^+) = \bigl[\partial_t(M_0\phi^+ + \mathcal{G}(\phi)) + \A\phi\bigr]_{t = 0} = 0.
  \label{eq:Maxwell_initial_history}
\end{equation}
This,  after a slight modification of $\phi^+$, can be interpreted as: $\phi$ must be a solution of the Maxwell system in $t = 0$.
Indeed, under the assumption that $\phi$ is differentiable in $t = 0$, let
\begin{equation*}
  \phi^+ = \phi(0^-)\theta^+\eta + (\partial_t\phi)(0^-)\theta^+\gamma,
\end{equation*}
where $\eta \in C_c^\infty(\R)$, $\eta(0) = 1, \eta'(0) = 0$ is as before, and $\gamma \in C_c^\infty(\R)$ satisfies $\gamma(0) = 0, \gamma'(0) = 1$.
Then, for $g_\phi$ defined as above, \eqref{eq:Maxwell_initial_history} becomes
\begin{equation}
  g_\phi(0^+) = \bigl[M_0(\partial_t\phi)(0^-) + \partial_t\mathcal{G}(\phi)(0^-) + \A\phi(0^-)\bigr]
  =  \partial_t\mathcal{M}(\phi)(0) + \A\phi(0)
  = 0.
  \label{eq:Maxwell_initial_history_improved}
\end{equation}
\end{rem}

\begin{ex}
  Let us formulate the above transformation $U \mapsto \tilde{u}$ in the original Maxwell variables $E,H$. Hence, consider
  \begin{equation*}
    \begin{aligned}
      \partial_tD(E) - \Curl H &= 0\\
      \partial_t\mu H - \Curl_0 E &= 0
    \end{aligned}
  \end{equation*}
  for $t > 0$, where
  $D(E) = \epsilon(\partial_t)E + P_\mathrm{nl}(E) = \epsilon_0E + \chi\ast E + \kappa\ast \snl(E)$.
  To simplify the notation, we denote the zero extension of the history of the fields by $E_0, H_0$. 
  Setting $E_0^+ \coleq E_0(0^-)\theta^+\eta,\ H_0^+ \coleq H_0(0^-)\theta^+\eta$ ($\eta \in C_c^\infty(\R)$ as before), we see that the resulting evolutionary system takes the form of \eqref{eq:Maxwell_nonlinear},
  \begin{equation}
    \left( \partial_t \begin{pmatrix} \epsilon(\partial_t) & 0 \\ 0 & \mu \end{pmatrix} + \begin{pmatrix} 0 & -\Curl \\ \Curl_0 & 0 \end{pmatrix} \right)
    \begin{pmatrix} \tilde{E} \\ \tilde{H} \end{pmatrix}
    = \begin{pmatrix} -\partial_t\tilde{P}_{\mathrm{nl}}(\tilde{E}) \\ 0 \end{pmatrix}
    + \begin{pmatrix} \varPhi \\ \varPsi \end{pmatrix},
    \label{eq:Maxwell_evo_proper_reformulation}
  \end{equation}
  with $\tilde{E} \coleq E-E_0^+$, $\tilde{H} \coleq H-H_0^+$, $\tilde{P}_\mathrm{nl}(\tilde{E}) \coleq P_{\mathrm{nl}}(\tilde{E} + E_0^+) = P_\mathrm{nl}(E)$, and the role of $g_\phi$ is played by
  \begin{equation}
    \begin{aligned}
      \varPhi &= -\theta^+\bigl[ \partial_t\bigl(\epsilon_0E_0^+ + \chi\ast E_0 + \kappa\ast\snl(E_0(\mult))\bigr) - \Curl H_0^+\bigr]\\
      \varPsi &= -\theta^+\bigl[ \partial_t\mu H_0^+ + \Curl_0 E_0^+\bigr].
    \end{aligned}
    \label{eq:Maxwell_evo_inhomogeneity}
  \end{equation}
  Thus, the history of $H$ only enters the equation via the initial value $H_0(0^-)$.
  To ensure $\Div B(H) = \Div\mu H = 0$ for $t \in (0,\infty)$, it suffices that $\Div\mu H_0(0^-) = 0$.
  In general, we will assume at least that $\mu H_0(0^-) \in H_0(\Div,\Omega)$.

  In order to apply Proposition \ref{prop:Lipschitz_perturb} or Proposition \ref{prop:FixPointRefinement} to the Maxwell problem, we need $\partial_t\tilde{P}_\mathrm{nl}(0) = 0$. This can, however, always be achieved
  by substituting $-\partial_t\tilde{P}_\mathrm{nl}(\tilde{E}) + \varPhi = -\partial_t(\tilde{P}_\mathrm{nl}(\tilde{E}) - \tilde{P}_\mathrm{nl}(0)) + (\varPhi - \partial_t\tilde{P}_\mathrm{nl}(0))$ on the right-hand side of \eqref{eq:Maxwell_evo_proper_reformulation}.
\end{ex}

Henceforth we shall drop the tilde and write $E, H$ instead of $\tilde{E}, \tilde{H}$, as well as $P_\mathrm{nl}$ instead of $\tilde{P}_\mathrm{nl}$, and always assume that the system is given in the evolutionary form \eqref{eq:Maxwell_evo_proper_reformulation}, where $\varPhi, \varPsi \in H^1_{\varrho}(\R,\Hs)$ are supported in $[0,\infty)$.
There is also no loss in assuming that $J = 0$, since a nonzero $J$ (supported in $[0,\infty)$ and $J \in H^1_\varrho(\R,\Hs)$) can be incorporated into the inhomogeneity $\varPhi$.

\subsection{Higher spatial regularity}
\label{sec:higher_spatial_regularity}
The interface setting and the choice of the domain $\Omega$ can be safely ignored in the results established so far (as long as the material laws are bounded linear operators on $L^2(\Omega)^3$).
In particular, $\Omega$ can be bounded or unbounded.
However, the heterogeneity of the material plays a more important role if tools relying on spatial regularity are used.
While spatial regularity is interesting in its own right, working in higher order Sobolev spaces also allows to control other types of nonlinearities for which Lipschitz-continuity fails in $L^2$ (compare in particular the local Lipschitz estimates in Theorem \ref{thm:NonlinFixL2} and \ref{thm:NonlinFixH2}).

From now on we assume that $\Omega = \Omega_1\sqcup\Gamma\sqcup\Omega_2$ is a \emph{bounded} domain with interface $\Gamma$ (the boundedness of $\Omega$ is a necessary requirement of Proposition~\ref{prop:spatial_regularity}).
We want to establish conditions that allow the solution $E,H$ for $k \in \N$ to lie (pointwise almost everywhere in time) in the space
\begin{equation*}
  \spatH{k} \coleq H^k(\Omega_1)^3\oplus H^k(\Omega_2)^3,
\end{equation*}
the latter being just the product $H^k(\Omega_1)^3\times H^k(\Omega_2)^3$ endowed with the sum-norm
\begin{equation*}
  \begin{aligned}
    \lVert (u_1,u_2) \rVert_{H^k(\Omega_1)\oplus H^k(\Omega_2)} &= \lVert u_1 \rVert_{H^k(\Omega_1)} + \lVert u_2 \rVert_{H^k(\Omega_2)}.
  \end{aligned}
\end{equation*}
(Note that for functions $u_1 \in H^k(\Omega_1)$, $u_2 \in H^k(\Omega_2)$ we identify the pair $(u_1,u_2)$ with the sum $u_1 + u_2$ of their zero extensions on $\Omega$; in particular, $\spatH{k}$ is a subspace of $\spatH{0} = \Hs = L^2(\Omega)^3$.)

We are in particular interested in $k \ge 2$ since in this case $H^k(\R^d)$ is a multiplication algebra in dimension $d \le 3$.
By extension, this carries over to $H^{k}(\Omega_i)$, \ie, for each bilinear map $b\colon \R^3\times \R^3 \to \R^3$ there exists a constant $C_b > 0$ such that
\begin{equation}
  \lVert b(u, v) \rVert_{H^{k}} \le C_b \lVert u \rVert_{H^{k}} \lVert v \rVert_{H^{k}}
  \label{eq:H2mult_algebra}
\end{equation}
for all $u,v \in H^{k}(\Omega_i)^3$.
The subsequent analysis relies on the following spatial regularity result adapted from \cite{weber81}.

\begin{prop}
  \label{prop:spatial_regularity}
  For some $k\in \N$ and $i \in \{1, 2\}$, let $\Omega_i\subset\R^3$ be bounded Lipschitz domains separated by the interface $\Gamma = \partial\Omega_1\cap\partial\Omega_2$ such that the complement of $\Omega \coleq \Omega_1\cup\Omega_2$ is simply connected, and $\partial\Omega_i$ are of class $C^{k+1}$.
  Let $\epsilon$ be a bounded matrix-valued function on $\Omega$, which, up to a complex factor, is Hermitian and uniformly positive definite, with $\epsilon, \epsilon^{-1} \in C^{k+1}(\overline{\Omega_i})^{3\times 3}$.
  Let $u \in L^2(\Omega)^3$ satisfy either of the conditions
  \begin{enumerate}
    \item[1.] $u \in H_0(\Curl,\Omega)$ and $\epsilon u \in H(\Div,\Omega)$, or
    \item[2.] $u \in H(\Curl,\Omega)$ and $\epsilon u \in H_0(\Div,\Omega)$.
  \end{enumerate}
  If $\ell \in \N$, $\ell \le k$, and $\Curl u \in H^{\ell-1}(\Omega_i)^3$, $\Div \epsilon u \in H^{\ell-1}(\Omega_i)$ for $i = 1,2$, then $u \in \spatH{\ell}$ and
  \begin{equation}
    \sum_{i = 1,2} \lVert u \rVert_{H^{\ell}(\Omega_i)}
    \le C_{\epsilon,\Omega_1,\Omega_2} \sum_{i=1,2} \bigl(\lVert u \rVert_{L^2(\Omega_i)} + \lVert \Curl u \rVert_{H^{\ell-1}(\Omega_i)} + \lVert \Div \epsilon u \rVert_{H^{\ell-1}(\Omega_i)}\bigr).
    \label{eq:weber_estimate_pw}
  \end{equation}
\end{prop}

\begin{rem}
  Note that the regularity conditions on $\partial\Omega_i$ of Proposition~\ref{prop:spatial_regularity} imply that the set $\partial\Omega\cup\partial\Omega_1\cup\partial\Omega_2$ splits into disjoint surfaces with positive distance from each other.
  Thus, a completely flat interface is prohibited under this setup;
  indeed, a straight interface inside a bounded domain $\Omega$ would generate corners in $\partial\Omega_1$ and $\partial\Omega_2$.
  This technical regularity assumption is used in \cite{weber81} to formulate the conditions on $u$ equivalently as a boundary value problem and a transmission problem.
  This setting can be possibly generalized. One prospective approach would be to use the embeddings
  \begin{equation*}
    \begin{aligned}
      H_0(\Curl,\Omega)\cap H(\Div,\Omega) &\subseteq H^1(\Omega)^3\\
      H(\Curl,\Omega)\cap H_0(\Div,\Omega) &\subseteq H^1(\Omega)^3,
    \end{aligned}
  \end{equation*}
  which hold for a $C^{1,1}$ regular boundary $\partial \Omega$
  or if $\Omega$ is a convex polyhedron, see \cite[§3.4, §3,5]{GiraultRaviart}.
\end{rem}

\begin{rem}
  \label{rem:eps_decomp}
  The conditions on $\epsilon$ imply that, upon multiplication with a complex phase, $\langle \epsilon \cdot, \cdot \rangle$ and $\langle \epsilon^{-1} \cdot, \cdot \rangle$ define equivalent inner products on $L^2(\Omega)^3$.
  With $\ker(\Div) \coleq \{ u \in H(\Div,\Omega) : \Div u = 0\}$ we have
  \begin{equation*}
    L^2(\Omega)^3
    = \Grad H_0^1(\Omega) \oplus \ker(\Div)
    = \Grad H_0^1(\Omega) \oplus_\epsilon \epsilon^{-1}\ker(\Div),
  \end{equation*}
  where $\oplus_\epsilon$ denotes orthogonal decomposition with respect to $\langle \epsilon \cdot, \cdot \rangle$.
  Moreover, the connectedness of $\R^3\setminus\Omega$ implies that for every $F \in \ker(\Div)$ there exists $A \in H(\Curl,\Omega) \cap \ker(\Div)$ with $F = \Curl A$.
  Thus, every $u \in L^2(\Omega)^3$ admits a Helmholtz decomposition
  \begin{equation*}
    u = \Grad f + \epsilon^{-1} \Curl A
  \end{equation*}
  with $f \in H^1_0(\Omega)$, $A \in H(\Curl,\Omega)$ and $\lVert \Grad f \rVert_{L^2} + \lVert \Curl A \rVert_{L^2} \le C_{\epsilon,\Omega} \lVert u \rVert_{L^2}$ (see \cite[Lemma 3.8]{weber80}, \cite[Lemma 3.7]{weber81}).
\end{rem}

The next result provides a conclusion similar to Proposition \ref{prop:spatial_regularity} if $\epsilon u$ is allowed to jump in normal direction across the interface.

\begin{prop}
  \label{prop:spatial_regularity_jump}
  Let $\Omega_i$ and $\epsilon$ satisfy the assumptions of Proposition \ref{prop:spatial_regularity} for some $k\in \N$ and assume for $\ell\in\N$, $\ell \le k$, $i = 1,2$ that $\Curl u \in H^{\ell-1}(\Omega_i)^3$, $\Div \epsilon u \in H^{\ell-1}(\Omega_i)$, and 
  \begin{equation*}
    u \in H_0(\Curl,\Omega),\quad \epsilon u + v \in H(\Div,\Omega)
  \end{equation*}
  for some $v \in \spatH{\ell}$.
  Then $u \in \spatH{\ell}$ and
  \begin{equation*}
    \sum_{i=1,2} \lVert u \rVert_{H^\ell(\Omega_i)} \le C_{\epsilon,\Omega_1,\Omega_2}
    \sum_{i=1,2} \bigl(\lVert u \rVert_{L^2(\Omega_i)} + \lVert \Curl u \rVert_{H^{\ell-1}(\Omega_i)} + \lVert \Div \epsilon u \rVert_{H^{\ell-1}(\Omega_i)} + \lVert v \rVert_{H^\ell(\Omega_i)}\bigr).
  \end{equation*}
\end{prop}

\begin{proof}
  In view of Remark \ref{rem:eps_decomp} we can write
  \begin{equation*}
    \epsilon^{-1}v = \Grad f + \epsilon^{-1}\Curl A
  \end{equation*}
  with $f \in H^1_0(\Omega)$, $A \in H(\Curl,\Omega)$ and $\lVert \Grad f \rVert_{L^2} \le C_{\epsilon,\Omega} \lVert \epsilon^{-1} \rVert \lVert v \rVert_{L^2}$, where $\lVert \epsilon^{-1} \rVert = \lVert \epsilon^{-1} \rVert_{L^\infty(\Omega)^{3\times 3}}$.
  In fact (see Lemma 5 in \cite[Chapter IX, §1.4]{DautrayLions3}), $\Grad f \in H^\ell(\Omega_i)^3$ ($i = 1,2$) since $\epsilon^{-1}v \in H^{\ell}(\Omega_i)^3$.
  Moreover, with $\Grad H_0^1(\Omega) \subseteq \ker(\Curl_0)$ and $\Curl H(\Curl,\Omega) \subseteq \ker(\Div)$ we have $\Curl_0 \Grad f = 0$ and $\Div \Curl A = 0$. Hence,
  \begin{equation*}
    \Curl_0 u = \Curl_0{} (u + \Grad f)\quad \text{and}\quad
    \Div{}(\epsilon u + v) = \Div \epsilon (u + \Grad f),
  \end{equation*}
  leading to $u + \Grad f \in H_0(\Curl,\Omega)$ and $\epsilon (u + \Grad f) \in H(\Div,\Omega)$.
  Applying Proposition~\ref{prop:spatial_regularity} now yields $u + \Grad f \in H^\ell(\Omega_i)^3$, together with the estimate
  \begin{multline*}
    \sum_{i=1,2}\lVert u + \Grad f \rVert_{H^\ell(\Omega_i)}\\
    \begin{aligned}
      &\le C_{\epsilon,\Omega_1,\Omega_2} \sum_{i=1,2} \bigl(\lVert u + \Grad f \rVert_{L^2(\Omega_i)} + \lVert \Curl (u + \Grad f) \rVert_{H^{\ell-1}(\Omega_i)} + \lVert \Div \epsilon(u + \Grad f) \rVert_{H^{\ell-1}(\Omega_i)}\bigr)\\
      &\le C_{\epsilon,\Omega_1,\Omega_2} \sum_{i=1,2} \bigl(\lVert u \rVert_{L^2(\Omega_i)} + \lVert \Curl u \rVert_{H^{\ell-1}(\Omega_i)} + \lVert \Div \epsilon u \rVert_{H^{\ell-1}(\Omega_i)} + \lVert \Grad f \rVert_{H^{\ell}(\Omega_i)}\bigr).
    \end{aligned}
  \end{multline*}
  Since $\Grad f \in H^\ell(\Omega_i)^3$, also $u \in H^\ell(\Omega_i)^3$ and the claim follows from
  \begin{equation*}
    \sum_{i=1,2}\lVert u + \Grad f \rVert_{H^\ell(\Omega_i)}
    \ge \sum_{i=1,2}\bigl(\lVert u \rVert_{H^\ell(\Omega_i)} - \lVert \Grad f \rVert_{H^\ell(\Omega_i)}\bigr)
    \ge C_\epsilon \sum_{i=1,2}\bigl(\lVert u \rVert_{H^\ell(\Omega_i)} - \lVert v \rVert_{H^\ell(\Omega_i)}\bigr).\qedhere
  \end{equation*}
\end{proof}

\begin{rem}
  \label{rem:uniform_bound_eps}
  To apply Proposition \ref{prop:spatial_regularity} to functions depending on time (or rather frequency in Fourier space), it is important to have uniform estimates, \ie, the dependence of the constants on $\epsilon$ should be removed.
  For simplicity, suppose that $\epsilon$ is constant on each $\Omega_i$ and
  let $u$ be given as in Proposition~\ref{prop:spatial_regularity}.
  In this case, a closer look at the proof in \cite{weber81} reveals for $\alpha\in\N^3$, $|\alpha| = \ell$, that
  \begin{equation*}
    \lVert \partial^\alpha(\epsilon u) \rVert_{L^2(\Omega_i)}
    \le C \bigl( \lVert \Curl u \rVert_{H^{\ell-1}(\Omega_i)} + \lVert \Div \epsilon u \rVert_{H^{\ell-1}(\Omega_i)} \bigr)
  \end{equation*}
  holds, with $C$ independent of $\epsilon$.
  One then easily obtains the estimate
  \begin{equation}
    \lVert u \rVert_{H^\ell(\Omega_i)}
    \le C(1 + \lVert \epsilon^{-1} \rVert) \bigl( \lVert u \rVert_{L^2(\Omega_i)} + \lVert \Curl u \rVert_{H^{\ell-1}(\Omega_i)} + \lVert \Div \epsilon u \rVert_{H^{\ell-1}(\Omega_i)} \bigr)
    \label{eq:estimate_eps_independent}
  \end{equation}
  with $C$ independent of $\epsilon$.

  Analogously, if, instead, $u$ and $v$ are chosen as in Proposition \ref{prop:spatial_regularity_jump}, then we have
  \begin{equation}
    \lVert u \rVert_{H^\ell(\Omega_i)}
    \le C(1 + \lVert \epsilon^{-1} \rVert) \bigl( \lVert u \rVert_{L^2(\Omega_i)} + \lVert \Curl u \rVert_{H^{\ell-1}(\Omega_i)} + \lVert \Div \epsilon u \rVert_{H^{\ell-1}(\Omega_i)} + \lVert v \rVert_{H^{\ell}(\Omega_i)} \bigr),
    \label{eq:estimate_eps_independent_jump}
  \end{equation}
  again with $C$ independent of $\epsilon$.

  For a generalization, suppose now that instead $\epsilon$ is not constant, but smooth, $\epsilon \in C^{k+1}(\Omega, \C^{3\times 3})$, and that all derivatives of $\epsilon, \epsilon^{-1}$ up to order $k+1$ are bounded.
  Then we obtain similar estimates of the form
  \begin{equation}
    \lVert u \rVert_{H^\ell(\Omega_i)}
    \le C_\epsilon \bigl( \lVert u \rVert_{L^2(\Omega_i)} + \lVert \Curl u \rVert_{H^{\ell-1}(\Omega_i)} + \lVert \Div \epsilon u \rVert_{H^{\ell-1}(\Omega_i)} + \lVert v \rVert_{H^{\ell}(\Omega_i)} \bigr),
    \label{eq:estimate_eps_independent_jump_higher}
  \end{equation}
  where $C_\epsilon$ depends on the norms of $\epsilon^{-1}$ and its derivatives.
\end{rem}

For simplicity, we assume for the rest of this section homogeneous materials in each $\Omega_i$, $i\in \{1,2\}$, where $\epsilon(\partial_t)$ and $\mu$ satisfy the following condition.
\begin{description}
  \item[(M1)\label{itm:MatCond_split}]
    The mapping $z\mapsto \begin{pmatrix} \epsilon(z) & 0 \\[-.2em] 0 & \mu \end{pmatrix}$ is a linear material law.
    There exist $\epsilon_1,\epsilon_2\colon \dom(\epsilon) \to \C^{3\times 3}$, such that for all $z \in \dom(\epsilon)$, up to a complex factor, $\epsilon_1(z), \epsilon_2(z)$ are Hermitian and positive definite, and
    \begin{equation*}
      \epsilon(z) = \epsilon_1(z) 1_{\Omega_1} + \epsilon_2(z) 1_{\Omega_2}.
    \end{equation*}
    Furthermore, $\mu = \mu_1{1}_{\Omega_1} + \mu_2{1}_{\Omega_2}$ with $\mu_1,\mu_2 \in \C^{3\times 3}$ Hermitian and uniformly positive definite.
\end{description}
Note in this case that $\epsilon(z)$ and $\mu$ satisfy the smoothness assumptions of Proposition \ref{prop:spatial_regularity} with $k = \infty$.

Under the condition \nameref{itm:MatCond_split} on $\epsilon(\partial_t)$ and $\mu$, we consider next the linearized system (in evolutionary form)
\begin{equation}
  \left(
  \partial_t
  \begin{pmatrix}
    \epsilon(\partial_t) & 0 \\ 0 & \mu
  \end{pmatrix} + 
  \begin{pmatrix}
    0 & -\Curl \\ \Curl_0 & 0
  \end{pmatrix}
  \right)
  \begin{pmatrix}
    E \\ H
  \end{pmatrix} =
  \begin{pmatrix}
    \varPhi - \ph \\ \varPsi
  \end{pmatrix}.
  \label{eq:MwLinear_jump}
\end{equation}

\begin{rem}
  The term $\ph$ acts purely as a placeholder for the nonlinearity $\partial_tP_\mathrm{nl}(E)$, which is added in Corollary \ref{cor:H2_regularity_nonlinear}.
  The jump of the material law at the interface leads to a jump of the electric field $E$ in normal direction, hence,
  if for $D(E) = \epsilon(\partial_t)E + P_\mathrm{nl}(E)$ we impose $\Div D = 0$ and hence $[n\cdot D]_\Gamma = 0$, then, in general, the normal components of $P_\mathrm{nl}(E)$ as well as $\ph = \partial_tP_\mathrm{nl}(E)$ are also discontinuous at the interface.
\end{rem}

For $\ell\in\N_0$ denote by $\lVert \cdot \rVert_{\varrho,\ell}$ the norm in $L^2_{\varrho}(\R,\spatH{\ell})$.

\begin{thm}[$H^2$-regularity for the linear Maxwell system]
  \label{thm:Maxwell_spatial_regularity}
  Impose the conditions on $\Omega$ of Proposition \ref{prop:spatial_regularity} with $k = 2$, and condition \nameref{itm:MatCond_split} on $\epsilon$ and $\mu$.
  Let $\varrho_0, c \in (0,\infty)$ be such that 
  $\Re z\epsilon(z) \ge c$ for $\Re z > \varrho_0$,
  making the system \eqref{eq:MwLinear_jump}
  well-posed in $L^2_{\varrho}(\R, \Hs)$, $\varrho > \varrho_0$.
  If additionally $\epsilon(z),\epsilon(z)^{-1}$ are uniformly bounded for $\Re z > \varrho_0$ and if the data $\varPhi,\varPsi,\ph$ fulfill
  \begin{align*}
    \varPhi &\in H^2_\varrho(\R, \Hs) \cap L^2_\varrho(\R, \spatH{1} \cap H(\Div,\Omega)),&
    \partial_t^{-1}\varPhi &\in L^2_\varrho(\R, H(\Div,\Omega)),&
    \Div\partial_t^{-1}\varPhi &\in L^2_\varrho(\R, \spatH{1})\\
    \varPsi &\in H^2_\varrho(\R, \Hs) \cap L^2_\varrho(\R, \spatH{1}\cap H_0(\Div,\Omega)),&
    \partial_t^{-1}\varPsi &\in L^2_\varrho(\R, H_0(\Div,\Omega)),&
    \Div\partial_t^{-1}\varPsi &\in L^2_\varrho(\R, \spatH{1}),
  \end{align*}
  and $\ph \in H^2_\varrho(\R, \spatH{2})$, then $E \in L^2_\varrho(\R, \spatH{2}\cap H_0(\Curl))$, $H \in L^2_\varrho(\R, \spatH{2}\cap H(\Curl))$ for all $\varrho > \varrho_0$ and
  \begin{equation}
    \begin{aligned}
      \lVert E \rVert_{\varrho,2} &\lesssim V(\varPhi,\varPsi) + \sum_{j=0}^2 \lVert \partial_t^j\ph \rVert_{\varrho,0} + \lVert \ph \rVert_{\varrho,2}\\
      \lVert H \rVert_{\varrho,2} &\lesssim V(\varPsi,\varPhi) + \sum_{j=0}^2 \lVert \partial_t^j\ph \rVert_{\varrho,0} + \lVert \ph \rVert_{\varrho,1},
      \label{eq:Maxwell_spatial_regularity}
    \end{aligned}
  \end{equation}
  where $V(f,g) = \sum_{j=0}^2 \lVert \partial_t^jf \rVert_{\varrho,0} + \lVert \partial_t^jg \rVert_{\varrho,0} + \lVert f \rVert_{\varrho,1} + \lVert \Div \partial_t^{-1} g \rVert_{\varrho,1}$ and
  with constants independent of $E,H,\varPhi,\varPsi,\ph, \varrho$.
\end{thm}

\begin{proof}
  As a result of the solution theory (Theorem \ref{thm:linear_solution_theory}) and the time-regularity (Proposition \ref{prop:time_regularity}), since $\varPhi,\varPsi,\ph \in H^2_\varrho(\R, \Hs)$, the solution fulfills $E,H \in H^2_\varrho(\R, \Hs)$ with continuous dependence on the data,
  \begin{equation}
    \lVert \partial_t^j E \rVert_{\varrho,0}, \lVert \partial_t^j H \rVert_{\varrho,0} 
    \le C \left(\lVert \partial_t^j(\varPhi - \ph) \rVert_{\varrho,0} + \lVert \partial_t^j\varPsi \rVert_{\varrho,0} \right),\quad j \in \{0,1,2\}.
    \label{eq:continuous_dependence_partial_t}
  \end{equation}
  Moreover, $(E,H) \in L^2_{\varrho}(\R, \dom(\A))$ and
  \begin{equation}
    \left.
    \begin{aligned}
      \partial_t \epsilon(\partial_t)E - \Curl H &= \varPhi - \ph \\
      \partial_t\mu H + \Curl_0 E &= \varPsi
    \end{aligned}
    \right\}
    \label{eq:Maxwell_H0}
  \end{equation}
  hold in $\Hw_{\varrho}(\R,L^2(\Omega)^3)$.
  Taking the Fourier--Laplace transform, we obtain from \eqref{eq:Maxwell_H0} pointwise for almost all $\xi \in \R$,
  \begin{equation}
    \left.
      \begin{aligned}
        \Curl \mathcal{L}_\varrho H(\xi)
        &= \epsilon(\varrho + i\xi)\mathcal{L}_\varrho (\partial_tE)(\xi) - \mathcal{L}_\varrho\varPhi(\xi) + \mathcal{L}_\varrho \ph(\xi) \\
        \Div \mu \mathcal{L}_\varrho H(\xi) &= \Div \mathcal{L}_\varrho(\partial_t^{-1}\varPsi)(\xi)\\
        \Curl_0 \mathcal{L}_\varrho E(\xi) &= -\mu \mathcal{L}_\varrho (\partial_tH)(\xi) + \mathcal{L}_\varrho\varPsi(\xi)\\
        \Div{}(\epsilon(\varrho + i\xi) \mathcal{L}_\varrho E(\xi) + \mathcal{L}_\varrho (\partial_t^{-1}\ph)(\xi)) &= \Div \mathcal{L}_\varrho(\partial_t^{-1}\varPhi)(\xi)
      \end{aligned}
    \right\}
    \label{eq:CurlDiv_Fourier}
  \end{equation}
  as identities in $L^2(\Omega)^3$.
  Now Proposition \ref{prop:spatial_regularity} (with $u = (\mathcal{L}_\varrho H)(\xi)$, using $u\in H_0(\Div) \cap H(\Curl)$) and Proposition \ref{prop:spatial_regularity_jump} (with $u = (\mathcal{L}_\varrho E)(\xi), v = (\mathcal{L}_\varrho \partial_t^{-1}\ph)(\xi)$) yield $\mathcal{L}_\varrho E(\xi),\mathcal{L}_\varrho H(\xi) \in H^1(\Omega_i)^3$.
  The same conclusion can be drawn for $\mathcal{L}_\varrho(\partial_tE), \mathcal{L}_\varrho(\partial_tH)$; indeed, apply $\partial_t$ to \eqref{eq:Maxwell_H0} and take the Fourier--Laplace transform to obtain the following identities in $L^2(\Omega)^3$, again for almost all $\xi \in \R$:
  \begin{equation}
    \left.
      \begin{aligned}
        \Curl \mathcal{L}_\varrho (\partial_tH)(\xi)
        &= \epsilon(\varrho + i\xi)\mathcal{L}_\varrho (\partial_t^2E)(\xi) - \mathcal{L}_\varrho(\partial_t\varPhi)(\xi) + \mathcal{L}_\varrho (\partial_t\ph)(\xi) \\
        \Div \mu \mathcal{L}_\varrho (\partial_tH)(\xi) &= \Div \mathcal{L}_\varrho\varPsi(\xi)\\
        \Curl_0 \mathcal{L}_\varrho (\partial_tE)(\xi) &= -\mu \mathcal{L}_\varrho (\partial_t^2H)(\xi) + \mathcal{L}_\varrho(\partial_t\varPsi)(\xi)\\
        \Div{}(\epsilon(\varrho + i\xi) \mathcal{L}_\varrho (\partial_tE)(\xi) + \mathcal{L}_\varrho \ph(\xi)) &= \Div \mathcal{L}_\varrho\varPhi(\xi),
      \end{aligned}
    \right\}
    \label{eq:CurlDiv_Fourier_partial_t}
  \end{equation}
  hence $\mathcal{L}_\varrho (\partial_tE)(\xi), \mathcal{L}_\varrho (\partial_tH)(\xi) \in H^1(\Omega_i)^3$ follows by the same argument.
  This, together with the assumptions on the data $\varPhi, \ph, \varPsi$ implies that \eqref{eq:CurlDiv_Fourier} are in fact identities in $H^1(\Omega_1)^3\oplus H^1(\Omega_2)^3$, which yields $\mathcal{L}_\varrho E(\xi),\mathcal{L}_\varrho H(\xi) \in H^2(\Omega_i)^3$ for almost all $\xi \in \R$, once more by Proposition \ref{prop:spatial_regularity} and Proposition \ref{prop:spatial_regularity_jump},
  together with the estimates (cf.\ Remark \ref{rem:uniform_bound_eps})
  \begin{align*}
    \lVert \mathcal{L}_\varrho E(\xi) \rVert_{H^2}
    &\le C(1 + \lVert \epsilon(\varrho + i\xi)^{-1} \rVert)
    \bigl( \lVert \mathcal{L}_\varrho E(\xi) \rVert_{L^2} + \lVert \Curl_0 \mathcal{L}_\varrho E(\xi) \rVert_{H^1}\\
    &\hspace{4.5cm} + \lVert \Div \epsilon(\varrho + i\xi)\mathcal{L}_\varrho E(\xi) \rVert_{H^1} + \lVert \mathcal{L}_\varrho \ph(\xi) \rVert_{H^2} \bigr)\\
    \lVert \mathcal{L}_\varrho H(\xi) \rVert_{H^2}
    &\le C(1 + \lVert \mu^{-1} \rVert)
    \left( \lVert \mathcal{L}_\varrho H(\xi) \rVert_{L^2} + \lVert \Curl \mathcal{L}_\varrho H(\xi) \rVert_{H^1} + \lVert \Div \mu\mathcal{L}_\varrho H(\xi) \rVert_{H^1} \right),
  \end{align*}
  where the norms are taken in $\Omega_i$.
  After integration, using the boundedness of $\mu^{-1}, \epsilon(\mult)^{-1}$ and the Plancherel theorem, we have
  \begin{equation}
    \left.
    \begin{aligned}
      \lVert E \rVert_{\varrho,2}
      &\lesssim \bigl(1 + \sup_{\Re z \ge \varrho_0}\lVert \epsilon(z)^{-1} \rVert\bigr) \left(\lVert E \rVert_{\varrho,0}
        + \lVert \Curl_0 E \rVert_{\varrho,1}
        + \lVert \Div \epsilon(\partial_t)E \rVert_{\varrho,1}
        + \lVert \ph \rVert_{\varrho,2}\right)\\
      \lVert H \rVert_{\varrho,2}
      &\lesssim (1 + \lVert \mu^{-1} \rVert) \left(\lVert H \rVert_{\varrho,0}
        + \lVert \Curl H \rVert_{\varrho,1}
        + \lVert \Div \mu H \rVert_{\varrho,1}\right)
    \end{aligned}
    \right\}
    \label{eq:DivCurl_H2_estimate}
  \end{equation}
  with constants independent of $E,H,\ph, \varrho$.
  Now we can employ \eqref{eq:Maxwell_H0} and estimates of type \eqref{eq:estimate_eps_independent}, \eqref{eq:estimate_eps_independent_jump} recursively to replace the $\Curl$- and $\Div$-terms on the right-hand side:
  \begin{align*}
    \lVert \Curl_0 E \rVert_{\varrho,1}
    &= \lVert -\mu\partial_tH + \varPsi  \rVert_{\varrho,1}\\
    &\le \lVert \partial_tH \rVert_{\varrho,1}
      + \lVert \varPsi  \rVert_{\varrho,1}\\
    &\lesssim \lVert \partial_tH \rVert_{\varrho,0}
      + \lVert \Curl \partial_tH \rVert_{\varrho,0}
      + \lVert \Div \mu\partial_tH \rVert_{\varrho,0}
      + \lVert \varPsi  \rVert_{\varrho,1}\\
    &\lesssim \lVert \partial_tH \rVert_{\varrho,0}
      + \lVert \epsilon(\partial_t)\partial_t^2E - \partial_t\varPhi + \partial_t\ph \rVert_{\varrho,0}
      + \lVert \Div \varPsi \rVert_{\varrho,0}
      + \lVert \varPsi  \rVert_{\varrho,1}\\
    &\lesssim \lVert \partial_tH \rVert_{\varrho,0}
      + \lVert \partial_t^2E \rVert_{\varrho,0}
      + \lVert \partial_t\varPhi \rVert_{\varrho,0}
      + \lVert \partial_t\ph \rVert_{\varrho,0}
      + \lVert \varPsi \rVert_{\varrho,1}\\
    \lVert \Div \epsilon(\partial_t)E \rVert_{\varrho,1}
    &\lesssim \lVert \Div \partial_t^{-1}\varPhi \rVert_{\varrho,1} + \lVert \Div \partial_t^{-1}\ph \rVert_{\varrho,1}\\
    &\le \lVert \Div \partial_t^{-1}\varPhi \rVert_{\varrho,1} + \lVert \partial_t^{-1}\ph \rVert_{\varrho,2}\\
    &\le \lVert \Div \partial_t^{-1}\varPhi \rVert_{\varrho,1} + \varrho^{-1}\lVert \ph \rVert_{\varrho,2}
    \intertext{and similarly,}
    \lVert \Curl H \rVert_{\varrho,1}
    &\lesssim \lVert \epsilon(\partial_t)\partial_tE - \varPhi + \ph \rVert_{\varrho,1}\\
    &\lesssim \lVert \partial_tE \rVert_{\varrho,1} + \lVert \varPhi \rVert_{\varrho,1} + \lVert \ph \rVert_{\varrho,1}\\
    &\lesssim \lVert \partial_tE \rVert_{\varrho,0} + \lVert \Curl\partial_tE \rVert_{\varrho,0} + \lVert \Div \epsilon(\partial_t)\partial_tE \rVert_{\varrho,0} + \lVert \varPhi \rVert_{\varrho,1} + \lVert \ph \rVert_{\varrho,1}\\
    &\lesssim \lVert \partial_tE \rVert_{\varrho,0} + \lVert -\mu\partial_t^2H + \partial_t\varPsi \rVert_{\varrho,0} + \lVert \Div (\varPhi -  \ph) \rVert_{\varrho,0} + \lVert \varPhi \rVert_{\varrho,1} + \lVert \ph \rVert_{\varrho,1}\\
    &\lesssim \lVert \partial_tE \rVert_{\varrho,0}
      + \lVert \partial_t^2H \rVert_{\varrho,0}
      + \lVert \partial_t\varPsi \rVert_{\varrho,0}
      + \lVert \varPhi \rVert_{\varrho,1}
      + \lVert \ph \rVert_{\varrho,1}\\
    \lVert \Div \mu H \rVert_{\varrho,1}
    &= \lVert \Div \partial_t^{-1}\varPsi \rVert_{\varrho,1}.
  \end{align*}
  Finally estimating $\lVert \partial_t^jE \rVert_{\varrho,0}, \lVert \partial_t^jH \rVert_{\varrho,0}$ with the help of \eqref{eq:continuous_dependence_partial_t} we obtain \eqref{eq:Maxwell_spatial_regularity}.
\end{proof}

\begin{rem}
  \label{rem:HigherRegularity_generalization_k}
  Theorem~\ref{thm:Maxwell_spatial_regularity} shows that one can trade regularity in time of order $k$ for regularity in space of the same order.
  Here we have fixed $k = 2$, but spatial $H^k$-regularity of the solution can be achieved for any $k \ge 1$.
  To this end, we must assume sufficient regularity of the data, for instance
  \begin{equation*}
    \begin{aligned}
      \varPhi &\in \bigcap_{j=0}^k H^j_\varrho(\R,\Hs^{k-j-1} \cap H(\Div,\Omega))\\
      \varPsi &\in \bigcap_{j=0}^k H^j_\varrho(\R,\Hs^{k-j-1} \cap H_0(\Div,\Omega))\\
      N &\in \bigcap_{j=0}^k H^j_\varrho(\R,\Hs^{k-j}).
    \end{aligned}
  \end{equation*}
  By proceeding inductively as in the proof of Theorem~\ref{thm:Maxwell_spatial_regularity}, one obtains estimates similar to \eqref{eq:Maxwell_spatial_regularity} of $\left\|E\right\|_{\varrho,k}$, $\left\|H\right\|_{\varrho,k}$ in terms of the higher Sobolev norms of $\varPhi,\varPsi,N$.
\end{rem}

We now extend the regularity to the nonlinear case, replacing $\ph$ by a map $\partial_tP_{\mathrm{nl}}(\cdot)$.

\begin{cor}[$H^2$-regularity for the nonlinear Maxwell system I]
  \label{cor:H2_regularity_nonlinear}
  Suppose that the domain $\Omega$ and the linear material law $\epsilon$ satisfy the conditions of Theorem~\ref{thm:Maxwell_spatial_regularity} and impose the regularity assumptions on the data $\varPhi, \varPsi \in L^2_\varrho(\R,L^2(\Omega)^3)$.
  Let $P_{\mathrm{nl}}$ be a nonlinearity for which
  \begin{align*}
    \partial_t^jP_\mathrm{nl}&\colon L^2_\varrho(\R,\spatH{2}) \to L^2_\varrho(\R,\spatH{2})\quad (j \in \{1,2,3\})
  \end{align*}
  are uniformly Lipschitz-continuous and satisfy $\partial_t^jP_\mathrm{nl}(0) = 0$ for $j \in \{1,2,3\}$, and
  \begin{equation*}
    L_{\varrho} \coleq
    \sum_{j=1}^3\lVert\partial_t^jP_\mathrm{nl}\rVert_{\Lip(L^2_\varrho(\R,\spatH{2}) \to L^2_\varrho(\R,\spatH{2}))}
    = o(1)\quad \text{as } \varrho \to \infty.
  \end{equation*}
  Then, for $\varrho_1 \ge \varrho_0$ large enough, the system
  \begin{equation}
    \left( \partial_t
      \begin{pmatrix}
        \epsilon(\partial_t) & 0 \\ 0 & \mu
      \end{pmatrix} +
      \begin{pmatrix}
        0 & -\Curl \\ \Curl_0 & 0
      \end{pmatrix}
    \right)
    \begin{pmatrix}
      E \\ H
    \end{pmatrix} =
    \begin{pmatrix}
      \varPhi - \partial_tP_{\mathrm{nl}}(E) \\ \varPsi
    \end{pmatrix}
    \label{eq:Maxwell_nl_H2}
  \end{equation}
  admits a unique solution $(E,H)$ with $E \in L^2_\varrho(\R,\spatH{2}\cap H_0(\Curl,\Omega))$, $H \in L^2_\varrho(\R,\spatH{2}\cap H(\Curl,\Omega))$ for $\varrho > \varrho_1$.
\end{cor}

\begin{proof}
  By assumption, $\partial_tP_\mathrm{nl}$ maps $\Hw_\varrho(\R,\spatH{2})$ into $H^2_\varrho(\R,\spatH{2})$.
  Define the map $T_\varrho$ ($\varrho \gg 0$) by
  \begin{equation}
    T_\varrho (u) \coleq 
    \Biggl(
      \overline{
      \partial_t
      \begin{pmatrix}
        \epsilon(\partial_t) & 0 \\ 0 & \mu
      \end{pmatrix} +
      \begin{pmatrix}
        0 & -\Curl \\ \Curl_0 & 0
      \end{pmatrix}
      }
    \Biggr)^{-1}
    \begin{pmatrix}
      \varPhi - \partial_tP_{\mathrm{nl}}(u) \\ \varPsi
    \end{pmatrix}.
    \label{eq:NonlinSolutionOperator}
  \end{equation}
  Then by Theorem \ref{thm:Maxwell_spatial_regularity}, $T_\varrho$ maps $L^2_\varrho(\R,\spatH{2})$ into $L^2_\varrho(\R,\spatH{2})^2$ and from \eqref{eq:Maxwell_spatial_regularity} we have the estimate
  \begin{equation*}
    \lVert T_\varrho(u) - T_\varrho(v) \rVert_{\varrho,2}
    \le C L_\varrho \lVert u - v \rVert_{\varrho,2}
  \end{equation*}
  with $C$ independent of $\varrho$.
  By assumption on $L_\varrho$, the product $CL_\varrho$ is smaller than $1$ for large $\varrho$, thus, denoting by $\pi_1,\pi_2:\C^3\times\C^3 \to \C^3$ the projections
  $\pi_1(u,v) = u$ and $\pi_2(u,v) = v$, the map
  $\pi_1\circ T_\varrho$ is a contraction on $L^2_\varrho(\R,\spatH{2})$ for large $\varrho > \varrho_0$.
  The solution to \eqref{eq:Maxwell_nl_H2} is given by the unique fixed point $E = \pi_1T_\varrho(E)$, $H = \pi_2T_\varrho(E)$.
\end{proof}

The following is a variant of Corollary~\ref{cor:H2_regularity_nonlinear} in the spirit of Proposition~\ref{prop:FixPointRefinement}.
\begin{cor}[$H^2$-regularity for the nonlinear Maxwell system II]
  \label{cor:H2_regularity_nonlinear_LocLip}
  Suppose that the domain $\Omega$ and the linear material law $\epsilon$ satisfy the conditions of Theorem~\ref{thm:Maxwell_spatial_regularity}.
  Impose the regularity assumptions on the data $\varPhi, \varPsi \in L^2_\varrho(\R,L^2(\Omega)^3)$ and suppose further that $\varPhi, \varPsi$ are supported in $[0,\infty)$.
  Let $P_{\mathrm{nl}}$ be a causal nonlinearity for which
  \begin{align*}
    \partial_t^jP_\mathrm{nl}&\colon L^2_\varrho(\R,\spatH{2}) \to L^2_\varrho(\R,\spatH{2})\quad (j \in \{1,2,3\})
  \end{align*}
  satisfy $\partial_t^jP_\mathrm{nl}(0) = 0$ as well as the local Lipschitz estimate
  \begin{equation*}
    \lVert \partial_t^jP_\mathrm{nl}(u) - \partial_t^jP_\mathrm{nl}(v) \rVert_{\varrho,2}
    \le C \left(\lVert u \rVert_{\varrho,2} + \lVert v \rVert_{\varrho,2}\right)^\alpha \lVert u - v \rVert_{\varrho,2}
  \end{equation*}
  with $C, \alpha > 0$.
  Then, for $\varrho_1 \ge \varrho_0$ large enough, the system \eqref{eq:Maxwell_nl_H2}
  admits a unique solution $(E,H)$ with $E \in L^2_\varrho(\R,\spatH{2}\cap H_0(\Curl,\Omega))$, $H \in L^2_\varrho(\R,\spatH{2}\cap H(\Curl,\Omega))$ for $\varrho > \varrho_1$.
\end{cor}

\begin{proof}
  The proof follows a similar idea to that of Proposition~\ref{prop:FixPointRefinement}. 
  Note that since $\supp\varPhi, \supp\varPsi \subseteq [0,\infty)$, $V_\varrho(\varPhi,\varPsi) = V(\varPhi,\varPsi)$ defined in Theorem~\ref{thm:Maxwell_spatial_regularity} fulfills  $V_\varrho(\varPhi,\varPsi) = o(1)$ as $\varrho \to \infty$.
  Consider the map $E\mapsto F_\varrho(E) \coleq \pi_1T_\varrho(E)$, where $T_\varrho$ is defined in \eqref{eq:NonlinSolutionOperator} and $\pi_1$ denotes the projection $\pi_1(E,H) = E$. Then by \eqref{eq:Maxwell_spatial_regularity} we have
  \begin{align*}
    \lVert F_\varrho(U) \rVert_{\varrho,2}
    &\lesssim V_\varrho(\varPhi, \varPsi) + \sum_{j=1}^3 \lVert \partial_t^jP_\mathrm{nl}(U) \rVert_{\varrho,2}
    \le c \bigl(V_\varrho(\varPhi, \varPsi) + \lVert U \rVert_{\varrho,2}^{\alpha + 1}\bigr),\\
    \lVert F_\varrho(U) - F_\varrho(V) \rVert_{\varrho,2}
    &\lesssim \sum_{j=1}^3 \lVert \partial_t^jP_\mathrm{nl}(U) - \partial_t^jP_\mathrm{nl}(V) \rVert_{\varrho,2}
    \le d(\lVert U \rVert_{\varrho,2} + \lVert V \rVert_{\varrho,2})^\alpha \lVert U - V \rVert_{\varrho,2}
  \end{align*}
  with $c,d > 0$.
  We conclude that $F_\varrho$ is a contraction on a closed ball $B_r \subset \Hw_\varrho(\R,\spatH{2})$ of sufficiently small radius $r > 0$, if $\varrho > 0$ is sufficiently large.
\end{proof}

\begin{ex}
  \label{ex:quad_nonlocal_H2}
  We adapt the quadratic model
  \begin{equation*}
    P^{(2)}(E)(t) = \iint \kappa(t-s_1, t-s_2) q(E(s_1), E(s_2)) \dd{s_1} \dd{s_2},
  \end{equation*}
  in Example~\ref{ex:NonlocalNonlinearity} to the $\spatH{2}$-setting.
  For simplicity we assume that $\kappa \in C^3((0,\infty)^2)^{3\times 3}$ is compactly supported, hence each time derivative of $P^{(2)}(E)$ is again of the form above:
  \begin{equation*}
    \partial_t^\ell P^{(2)}(E)(t) = \iint (\partial_1 + \partial_2)^\ell\kappa(t-s_1,t-s_2)\,\snl(E(s_1),E(s_2))\dd{s_1}\dd{s_2},\quad \ell \in \N,
  \end{equation*}
  with $(\partial_1 + \partial_2)^\ell\kappa \in C((0,\infty)^2)^{3\times 3}$ bounded.

  As for the spatial nonlinearity, by the algebra property of $\spatH{2}$, each bilinear map $q\colon \R^3\times\R^3 \to \R^3$ extends to a bilinear and bounded map $q\colon \spatH{2}\times\spatH{2} \to \spatH{2}$, and we have
  \begin{equation*}
    \lVert q(u(\tau_1), u(\tau_2)) - q(v(\tau_1), v(\tau_2)) \rVert_{\spatH{2}}
    \lesssim \lVert u(\tau_1) - v(\tau_1) \rVert_{\spatH{2}} \lVert u(\tau_2) \rVert_{\spatH{2}}
    + \lVert u(\tau_2) - v(\tau_2) \rVert_{\spatH{2}} \lVert u(\tau_1) \rVert_{\spatH{2}}.
  \end{equation*}
  This can be generalized to $x$-dependent bilinear maps like
  \begin{equation*}
    \tilde{q}(u,v)(x) = M(x)q(u(x),v(x))
  \end{equation*}
  with $q$ as above and $M$ of class $H^{2}$ on each $\Omega_i$ ($i \in \{1,2\}$), as well as to $q$ like in Example~\ref{ex:NonlocalNonlinearity}, \ie,
  \begin{equation*}
    {q}(u,v)(x) = \int_{\Omega}\int_{\Omega} \Lambda(x,y_1,y_2)u(y_1)v(y_2)\dd{y_1}\dd{y_2}
  \end{equation*}
  with the same regularity assumption on $\Lambda = (\Lambda_{ijk})_{i,j,k\in\{1,2,3\}}$.

  Now using a \emph{smooth} cutoff in time, \ie,
  \begin{equation*}
    P^{(2)}_T \coleq \eta\cdot P^{(2)},
    \quad \eta \in C_c^\infty(\R)
  \end{equation*}
  and denoting $T \coleq \left|\supp \eta\right|$, we obtain for $j \in \{1,2,3\}$ and $\varrho \in \R$ similar to \eqref{eq:MultilinearLocLip_T}
  \begin{equation*}
    \lVert \partial_t^jP^{(2)}_T(U) - \partial_t^jP^{(2)}_T(V) \rVert_{\varrho,2}
    \le C\sqrt{T}e^{\varrho T}\left( \lVert U \rVert_{\varrho,2} + \lVert V \rVert_{\varrho,2} \right)\lVert U - V \rVert_{\varrho,2}
  \end{equation*}
  with $C$ independent of $U,V, \varrho$.
\end{ex}

\section{Exponential stability of the Maxwell system on a bounded domain}
\label{sec:exp_stab}

In \cite{Trostorff_habil}, exponential stability for linear equations
\begin{equation}
  \bigl( \partial_t M(\partial_t) + A \bigr) u = f
  \label{eq:evo_eq_generic}
\end{equation}
was investigated.
Assuming \eqref{eq:evo_eq_generic} is well-posed in the range of spaces $L^2_\varrho(\R,\HX)$, $\varrho > \varrho_0$,
the equation is said to be \emph{exponentially stable with decay rate $\nu_0 > 0$}, if the implication
\begin{equation}
  f \in \Hw_\varrho(\R,\HX) \cap \Hw_{-\nu}(\R,\HX) \implies u \in \Hw_{-\nu}(\R,\HX)\qquad (0 \le \nu < \nu_0,\ \varrho > \varrho_0)
  \label{eq:ExpStab_definition}
\end{equation}
holds.
If $s_0(M,A)$ denotes the infimum over $\varrho \in \R$ such that the equation is well-posed in $\Hw_{\varrho}(\R,\HX)$, then exponential stability with rate $\nu_0$ is essentially equivalent (under some natural assumptions on the domain of $M$, \cite[Theorem 2.1.3]{Trostorff_habil}) to $s_0(M,A) \le -\nu_0$.
This implies in particular the continuous dependence on the data in the $L^2_{-\nu}$-norm, \ie,
$\lVert u \rVert_{-\nu,0} \le K \lVert f \rVert_{-\nu,0}$ with $K$ independent of $f\in \Hw_{\varrho}(\R,\HX)\cap\Hw_{-\nu}(\R,\HX)$.

There are two abstract criteria to ensure exponential stability of the linear system \eqref{eq:evo_lin_1order}. The first requires strict and uniform accretivity of $z M(z)$.
(See also \cite[Chapter 11]{STW22}.)
Recall that a linear operator $A\colon \dom(A) \subset \HX \to \HX$ on a Hilbert space $\HX$ is called $m$-accretive, if $\Re \langle Ax,x \rangle \ge 0$ for all $x\in \dom(A)$ and $A + \lambda$ is onto for all $\lambda \in \C_{\Re > 0}$. In particular, every skew-selfadjoint operator is $m$-accretive.

\begin{prop}[{\cite[Theorem 2.1.5]{Trostorff_habil}}]
  \label{prop:exp_stab_criterion_1base}
    Let $A$ be m-accretive, let $\nu_0 > 0$ be such that $\C_{\Re > -\nu_0} \setminus \dom(M)$ is discrete, and assume that
    \begin{equation*}
      \exists c > 0\, \forall z \in \C_{\Re > -\nu_0} \cap \dom(M):\quad
      \Re z M(z) \ge c.
    \end{equation*}
    Then the linear problem \eqref{eq:evo_eq_generic} is well-posed and exponentially stable with decay rate $\nu_0$.
\end{prop}

The latter condition can be relaxed if one assumes that $A$ is invertible.
Let $B[0,\delta] = \{ z \in \C : |z| \le \delta \}$ denote the closed disk with radius $\delta > 0$.

\begin{prop}[{\cite[Theorem 2.1.6]{Trostorff_habil}}]
  \label{prop:exp_stab_criterion_1inv}
  Let $A$ be m-accretive and boundedly invertible, let $\nu_0 > 0$ be such that $\C_{\Re > -\nu_0} \setminus \dom(M)$ is discrete, and assume that for some $\delta > 0$
  \begin{itemize}
      \item[(i)]
          $\displaystyle \sup_{z \in \dom(M) \cap B[0,\delta]} \lVert z M(z) \rVert < \frac{1}{\lVert A^{-1} \rVert}$
      \item[(ii)]
        $\displaystyle \exists c > 0\, \forall z \in \C_{\Re > -\nu_0} \cap \dom(M) \setminus B[0,\delta]:\ \Re z M(z) \ge c.$
  \end{itemize}
  Then the linear problem \eqref{eq:evo_eq_generic} is well-posed and exponentially stable with decay rate $\nu_0$.
\end{prop}

These criteria can be applied to second-order equations of the form
\begin{equation}
    \bigl( \partial_t^2 M(\partial_t) + C^*C \bigr)u = f
    \label{eq:evo_2order}
\end{equation}
where $C\colon \dom(C) \subseteq \HX_0 \to \HX_1$ is assumed boundedly invertible and $M(z) = M_0(z) + z^{-1}M_1(z)$, with $M_0, M_1 \colon \dom(M) \subseteq \C \to \mathcal{B}(\HX_0)$ analytic and bounded.
The strategy relies on the substitution
\begin{equation*}
    v = \partial_t u + d u, \qquad q = -Cu \qquad (d > 0)
\end{equation*}
which converts \eqref{eq:evo_2order} into the first-order system
\begin{equation}
  \left( \partial_t \mathcal{M}_d(\partial_t) + \begin{pmatrix} 0 & -C^* \\ C & 0 \end{pmatrix} \right) \begin{pmatrix} v \\ q \end{pmatrix} = \begin{pmatrix} f \\ 0 \end{pmatrix},
  \label{eq:M_d_evo_sys}
\end{equation}
where
\begin{equation*}
  \mathcal{M}_d(z) =
    \begin{pmatrix} M(z) & 0 \\ 0 & 1 \end{pmatrix}
    + z^{-1} d \begin{pmatrix} -M_0(z) & \bigl( M_1(z) - d M_0(z) \bigr) C^{-1} \\ 0 & 1 \end{pmatrix}.
\end{equation*}
The accretivity properties of $z M(z)$ are inherited by $z\mathcal{M}_d(z)$, provided $d > 0$ is sufficiently small, and Proposition \ref{prop:exp_stab_criterion_1base} resp.\ Proposition \ref{prop:exp_stab_criterion_1inv}  are applicable, see \cite[§2.2]{Trostorff_habil}, \cite[Section 11.4]{STW22}.
We mention explicitly the following result.

\begin{thm}
  \label{thm:exp_stab_M_d}
  Let $M$ be given by $M(z) \coleq M_0(z) + z^{-1}M_1(z)$, where $M_0, M_1\colon \dom(M) \subseteq \C \to \BLO(\HX_0)$ are analytic and bounded, $\C_{\Re>-\tilde{\nu}}\setminus\dom(M)$ is discrete for some $\tilde{\nu} > 0$, and $\lim_{z\to 0}M_1(z) = 0$.
  If the condition
  \begin{equation*}
    \forall \delta > 0\ \exists \nu,c > 0\ \forall z\in\dom(M)\cap\C_{\Re>-\nu}\setminus B[0,\delta]:\quad
    \Re zM(z) \ge c
  \end{equation*}
  is met, then there exist $d_0, \nu_0 > 0$ such that system \eqref{eq:M_d_evo_sys} with $d = d_0$ is exponentially stable with decay rate $\nu_0$.
\end{thm}

\begin{rem}
  \label{rem:u_H1_2ord}
  We note that $v,q \in \Hw_{-\nu}(\R,\HX_0)$ already implies $u \in H^1_{-\nu}(\R,\HX_0)$ and $Cu \in \Hw_{-\nu}(\R,H_1)$, with
  \begin{equation*}
    \lVert u \rVert_{-\nu,0}, \lVert \partial_tu \rVert_{-\nu,0}, \lVert Cu \rVert_{-\nu,0}
    \lesssim \lVert f \rVert_{-\nu,0}.
  \end{equation*}
\end{rem}

\begin{rem}
  We will subsequently use Theorem 3.3 to formulate stability results for Maxwell's equations.
  As before, the presence of an interface will play no role at first, and only be of importance in the later Section \ref{sec:ExpStab_H2} when higher spacial regularity is involved.
\end{rem}

\subsection{Exponential stability of the \texorpdfstring{$E$}{E}-field via the second-order formulation}
\label{sec:ExpStab_SecondOrder}

A first observation is that 
the material law $M(\partial_t) = \Bigl(\begin{smallmatrix} \epsilon(\partial_t) & 0 \\ 0 & \mu \end{smallmatrix}\Bigr)$
associated with the linearized system
\begin{equation}
  \left( \partial_t \begin{pmatrix} \epsilon(\partial_t) & 0 \\ 0 & \mu \end{pmatrix} + \begin{pmatrix} 0 & -\Curl \\ \Curl_0 & 0 \end{pmatrix} \right)
  \begin{pmatrix} E \\ H \end{pmatrix} =
  \begin{pmatrix} \varPhi \\ \varPsi \end{pmatrix}
  \label{eq:Mw1ord_2}
\end{equation}
does not fulfill any of the strict accretivity conditions above, since
\begin{equation*}
  \Re \left< z M(z) \binom{0}{u}, \binom{0}{u} \right> =  (\Re z)\langle \mu u, u \rangle = 0,
\end{equation*}
whenever $\Re z = 0$, independently of $\epsilon$.
However, as we show next, in the second-order formulation under suitable assumptions on $\epsilon$ the accretivity assumption of Theorem~\ref{thm:exp_stab_M_d} is satisfied.

Recall that well-posedness of \eqref{eq:Mw1ord_2} means that
\begin{equation*}
  \left( \overline{\partial_t \begin{pmatrix} \epsilon(\partial_t) & 0 \\ 0 & \mu \end{pmatrix} + \begin{pmatrix} 0 & -\Curl \\ \Curl_0 & 0 \end{pmatrix}} \right)
  \begin{pmatrix} E \\ H \end{pmatrix} =
  \begin{pmatrix} \varPhi \\ \varPsi \end{pmatrix}
\end{equation*}
holds as an equation in $L^2_\varrho(\R,\Hs)^2$ for $\varrho$ large enough.
If the data satisfy the regularity assumption $\varPhi,\varPsi \in H^1_\varrho(\R,\Hs)$, then by Theorem \ref{thm:linear_solution_theory} (i) we can drop the closure bar and \eqref{eq:Mw1ord_2} itself holds in $L^2_\varrho(\R,\Hs)^2$.
Applying $\partial_t$ to the first line, we can insert the second line via $\partial_t\Curl H = \Curl\partial_t H = \Curl\mu^{-1}\varPsi - \Curl\mu^{-1}\Curl_0 E$ to obtain the second-order system
\begin{equation}
  (\partial_t^2 \epsilon(\partial_t) + \Curl \mu^{-1} \Curl_0) E = \partial_t\varPhi + \Curl \mu^{-1}\varPsi.
  \label{eq:Mw2ordE}
\end{equation}
We will impose the following conditions on the permittivity.

\begin{description}
  \item[(M2)\label{itm:MatCond_accr1}]
    For all $\delta > 0$ there exist $\nu > 0$ and $c > 0$ such that
    \begin{equation*}
      \forall z \in \C_{\Re > -\nu} \setminus B[0,\delta]: \quad \Re z\epsilon(z) \ge c.
    \end{equation*}
  \item[(M3)\label{itm:MatCond_accr2}] 
    $\epsilon(\partial_t) = \epsilon_0 + \chi(\partial_t)$ with $\epsilon_0 \in \BLO(\Hs)$ selfadjoint and uniformly positive definite,
    $\lim_{z\to 0} z\chi(z) = 0$, and
    there exists $\nu_1 > 0$ such that $z\mapsto \chi(z)$ and $z\mapsto z\chi(z)$ are bounded in $\C_{\Re > -\nu_1}$ and
    \begin{equation*}
      \forall z \in \C_{\Re > -\nu_1}:\quad  \epsilon_0 + \Re \chi(z) \ge c_1
    \end{equation*}
    with $c_1 > 0$.
\end{description}
Under these assumptions, $\epsilon(z) = \epsilon_0 + z^{-1}(z\chi(z))$ satisfies the conditions of Theorem~\ref{thm:exp_stab_M_d} with $M_0 = \epsilon_0$ and $M_1(z) = z\chi(z)$.

\begin{rem}
  An example of a physically relevant $\epsilon(\partial_t)$ compatible with \nameref{itm:MatCond_accr1} and \nameref{itm:MatCond_accr2} is the one given by the Drude--Lorentz model, see Appendix \ref{sec:app1}.
\end{rem}

Conditions \nameref{itm:MatCond_accr1} and \nameref{itm:MatCond_accr2} are sufficient to obtain a notion of exponential stability for the second-order system \eqref{eq:Mw2ordE} on a bounded domain $\Omega$, which is similar to that in \eqref{eq:ExpStab_definition}.
Note that \eqref{eq:Mw2ordE} is not yet of the form \eqref{eq:evo_2order}, because $\Curl, \Curl_0$ are not invertible.
We need two preparatory results, in which we follow a strategy akin to that in \cite{TroWau14}.

\begin{lem}
  \label{lem:Projections_invertible}
  Let $\mathscr{K}_0,\mathscr{K}_1$ be Hilbert spaces and $C \colon \dom(C) \subseteq \mathscr{K}_0 \to \mathscr{K}_1$ a densely defined and closed operator with closed range. Let $\mu\in \BLO(\mathscr{K}_1)$ be selfadjoint and uniformly positive definite.
  Denote by $\iota_k \colon \ker(C)^\bot \hookrightarrow \mathscr{K}_0$ the canonical embedding.
  Then
  \begin{equation*}
    \iota_k^*C^*\mu\, C\iota_k \colon \dom(\iota_k^*C^*\mu\, C\iota_k)\subseteq \ker(C)^\bot \to \ker(C)^\bot
  \end{equation*}
  is selfadjoint, continuously invertible, and nonnegative.
\end{lem}

\begin{proof}
  Let $\iota_r:\ran(C) \hookrightarrow \mathscr{K}_1$ be the canonical embedding. Then $\iota_r^*C\iota_k$ is injective, surjective, and closed, thus, by the closed graph theorem, continuously invertible. Note also that since $\ran(C)$ is closed, $\ran(C^*)$ is closed (see \cite[Theorem 2.19]{brezis}).
  Now it is not difficult to see that
  \begin{equation*}
    \iota_k^*C^*\iota_r\iota_r^*\,\mu\,\iota_r\iota_r^*C\iota_k = \iota_k^*C^*\mu\, C\iota_k
  \end{equation*}
  (for $C^*\iota_r\iota_r^* = C^*$ use that $\ker(C^*)^\bot = \ran(C)$).
  It follows that $\iota_k^*C^*\mu\,C\iota_k$ is the composition of the three continuously invertible mappings
  \begin{equation*}
    \iota_k^*C^*\iota_r,\quad \iota_r^*\mu\,\iota_r,\quad \iota_r^*C\iota_k,
  \end{equation*}
  and that $(\iota_k^*C^*\mu\,C\iota_k)^* = \iota_k^*C^*\mu\,C\iota_k$. Thus we obtain continuous invertibility and selfadjointness. The nonnegativity follows from the nonnegativity of $\mu$.
\end{proof}

\begin{rem}
  \label{rem:Projections_invertible_selfadjoint_op}
  In the situation of the previous lemma we have
  \begin{equation*}
    \iota_k^*C^*\mu\,C\iota_k = B^*B
  \end{equation*}
  for some continuously invertible operator $B \colon \dom(B) \subseteq \ker(C)^\bot \to \ker(C)^\bot$.
  Indeed, this is a direct consequence of the lemma in conjunction with the spectral theorem for unbounded selfadjoint operators.
\end{rem}

\begin{lem}
  \label{lem:proj_positive}
  Let $\mathscr{H}$ be a Hilbert space and $\mathscr{H}_0 \subset \mathscr{H}$ a closed subspace.
  Denote by $\iota_0\colon \mathscr{H}_0 \hookrightarrow \mathscr{H}$ and $\iota_1\colon \mathscr{H}_0^\bot \hookrightarrow \mathscr{H}$ the canonical embeddings.
  Let $T \in \BLO(\mathscr{H})$ be a bounded linear operator and define
  \begin{equation*}
    T_{jk} \coleq \iota_j^* T \iota_k^{\phantom{*}}\quad \text{for } j,k \in \{0,1\}.
  \end{equation*}
  If $\Re T \ge d$ for some $d > 0$, then also
  \begin{equation*}
    \Re T_{11}^{\phantom{1}} \ge d, \qquad \Re \bigl( T_{00}^{\phantom{1}} - T_{01}^{\phantom{1}}T_{11}^{-1}T_{10}^{\phantom{1}} \bigr) \ge d.
  \end{equation*}
\end{lem}

\begin{proof}
  For $\phi \in \mathscr{H}_0^\bot$ we compute
  \begin{align*}
    \Re \langle T_{11}\phi, \phi \rangle
    = \Re \langle T \iota_1 \phi, \iota_1 \phi \rangle
    \ge d \langle \iota_1\phi, \iota_1\phi \rangle = d \lVert \phi \rVert^2,
  \end{align*}
  confirming $\Re T_{11} \ge d$. In particular, $T_{11}$ is invertible.
  As an operator on $\mathscr{H}_0\oplus\mathscr{H}_0^\bot$ we can identify
  \begin{equation*}
    T =
    \begin{pmatrix}
      T_{00} & T_{01} \\ T_{10} & T_{11}
    \end{pmatrix},
  \end{equation*}
  and setting
  \begin{equation*}
    Q =
    \begin{pmatrix}
      1 & 0 \\[.4em] -\bigl(T_{01}^{\phantom{1}}T_{11}^{-1}\bigr)^* & 1
    \end{pmatrix}, \quad
    Q^* =
    \begin{pmatrix}
      1 & -T_{01}^{\phantom{1}}T_{11}^{-1} \\[.4em] 0 & 1
    \end{pmatrix}, \quad
    R =
    \begin{pmatrix}
      T_{00}^{\phantom{1}} - T_{01}^{\phantom{1}}T_{11}^{-1}T_{10}^{\phantom{1}} & 0 \\[.4em] T_{10}^{\phantom{1}} - T_{11}^{\phantom{1}}\bigl(T_{11}^{-1}\bigr)^{*}T_{01}^{*} & T_{11}^{\phantom{1}}
    \end{pmatrix}
  \end{equation*}
    we have the factorization
    \begin{equation*}
      R = Q^*TQ.
    \end{equation*}
    Now we compute for $\phi \in \mathscr{H}_0$,
    \begin{align*}
      \Re\langle (T_{00}^{\phantom{1}} - T^{\phantom{1}}_{01}T_{11}^{-1}T_{10}^{\phantom{1}}) \phi, \phi \rangle
      &= \Re\langle R\begin{pmatrix} \phi \\ 0 \end{pmatrix}, \begin{pmatrix} \phi \\ 0 \end{pmatrix} \rangle \\
      &= \Re\langle Q^*TQ\begin{pmatrix} \phi \\ 0 \end{pmatrix}, \begin{pmatrix} \phi \\ 0 \end{pmatrix} \rangle \\
      &= \Re\langle TQ\begin{pmatrix} \phi \\ 0 \end{pmatrix}, Q\begin{pmatrix} \phi \\ 0 \end{pmatrix} \rangle
      \ge d \langle Q\begin{pmatrix} \phi \\ 0 \end{pmatrix}, Q\begin{pmatrix} \phi \\ 0 \end{pmatrix} \rangle
      \ge d \lVert \phi \rVert^2. \qedhere
    \end{align*}
\end{proof}

In the following, we apply this result to $\mathscr{H} = \Hs = L^2(\Omega)^3$, where $\Omega$ is a bounded weak Lipschitz domain (a bounded domain with a local Lipschitz boundary), and $\mathscr{H}_0 = \ran(\Curl)$.
As a consequence of the Picard--Weber--Weck selection theorem (see Lemma \ref{lem:RanCurl_closed}), $\mathscr{H}_0$ is a closed subspace of $\Hs$. 
Lemma~\ref{lem:proj_positive} can be applied to the second-order formulation \eqref{eq:Mw2ordE} if $\varPhi, \varPsi$ are regular enough, as we show next.
Let $\Pi_0 \colon \Hs \to \ker(\Curl_0)$ denote the canonical projection.

\begin{thm}[Exponential stability for the $E$-field of the linear Maxwell equations]
  \label{thm:Maxwell_ExpStab_SecondOrder}
  Let $\Omega \subset \R^3$ be a bounded, weak Lipschitz domain.
  Assume that the material law $\epsilon(\partial_t) = \epsilon_0 + {\chi}(\partial_t)$ satisfies \nameref{itm:MatCond_accr1} and \nameref{itm:MatCond_accr2}, and that $\mu\in\BLO(\Hs)$ is selfadjoint and uniformly positive definite.
  Let $(\varPhi,\varPsi) \in \Hw_\varrho(\R,\Hs\times\Hs)$ for large $\varrho > 0$, and let $(E,H) \in L^2_{\varrho}(\R,\Hs\times\Hs)$ be the unique solution to \eqref{eq:Mw1ord} provided by Theorem~\ref{thm:linear_solution_theory}.
  Let $g \coleq \partial_t\varPhi + \Curl\mu^{-1}\varPsi$ and $h \coleq\Pi_0\, \partial_t^{-1}\varPhi$.
  Then there exists $\nu_0 > 0$ such that if $\nu < \nu_0$, $g \in L^2_{-\nu}(\R,\Hs)$, and $h \in L^2_{-\nu}(\R, \ker(\Curl_0))$, then $E \in L^2_{-\nu}(\R, \Hs)$ and
  \begin{equation}
    \lVert E \rVert_{-\nu,0}
    \le K \bigl( \lVert g \rVert_{-\nu,0} + \lVert h \rVert_{-\nu,0} \bigr)
    \label{eq:expstab_lin_L2L2_estimate}
  \end{equation}
  with $K > 0$ independent of $E$, $g$, $h$.
\end{thm}

\begin{proof}
  Let $\varPhi,\, \varPsi \in C_c^\infty(\R,\Hs)$, and let $E,H \in \Hw_{\varrho}(\R,\Hs)$ for large $\varrho > 0$ fulfill
  \begin{equation*}
    \Biggl(
      \overline{
        \partial_t
        \begin{pmatrix}
          \epsilon(\partial_t) & 0 \\ 0 & \mu
        \end{pmatrix} +
        \begin{pmatrix}
          0 & -\Curl \\ \Curl_0 & 0
        \end{pmatrix}
      }
    \Biggr)
    \begin{pmatrix}
      E \\ H
    \end{pmatrix} =
    \begin{pmatrix}
      \varPhi \\ \varPsi
    \end{pmatrix}.
  \end{equation*}
  Due to the regularity of the right-hand side we can drop the closure bar and obtain the second-order system \eqref{eq:Mw2ordE} for $E$, thus
  \begin{equation}
    (\partial_t^2\epsilon(\partial_t) + \Curl\mu^{-1}\Curl_0)E = g.
    \label{eq:Mw_expstab2ord}
  \end{equation}
  The aim is to show that $E \in \Hw_\varrho(\R,\Hs)$ additionally satisfies $E \in \Hw_{-\nu}(\R,\Hs)$, if $g \in \Hw_{-\nu}(\R,\Hs)$ and $h \in \Hw_{-\nu}(\R,\ker(\Curl_0))$ for small $\nu > 0$.
  To this end, we set
  \begin{equation*}
      \mathscr{H}_0 \coleq \ran(\Curl), \quad
      \mathscr{H}_0^\bot = \ker(\Curl_0).
  \end{equation*}
    Denote by $\iota_0\colon \mathscr{H}_0 \hookrightarrow \Hs$, $\iota_1\colon \mathscr{H}_0^\bot \hookrightarrow \Hs$ the canonical embeddings and define
  \begin{equation*}
    E_j \coleq \iota_j^*E,\quad g_j \coleq \iota_j^*g, \quad
    \epsilon_{jk}(\partial_t) \coleq \iota_j^*\epsilon(\partial_t)\iota_k\quad \text{for } j,k \in \{0,1\}
  \end{equation*}
  (note that $\iota_1^\ast = \Pi_0$, thus $h = \partial_t^{-2}g_1$).
  Then \eqref{eq:Mw_expstab2ord} can be written equivalently in the form
  \begin{equation}
    \left[
      \partial_t^2
      \begin{pmatrix}
        \epsilon_{00}(\partial_t) & \epsilon_{01}(\partial_t) \\ \epsilon_{10}(\partial_t) & \epsilon_{11}(\partial_t)
      \end{pmatrix} +
      \begin{pmatrix}
        \iota_0^*\Curl\mu^{-1}\Curl_0\iota_0 & 0 \\ 0 & 0
      \end{pmatrix}
    \right]
    \begin{pmatrix}
      E_0 \\ E_1
    \end{pmatrix} =
    \begin{pmatrix}
      g_0 \\ g_1
    \end{pmatrix}.
    \label{eq:big_proj_system}
  \end{equation}
  Next we observe that, by assumption \nameref{itm:MatCond_accr2},
  \begin{equation*}
    \partial_t\epsilon_{01}(\partial_t) \left( \partial_t\epsilon_{11}(\partial_t) \right)^{-1}\colon
    \Hw_{-\nu}(\R,\mathscr{H}_0^\bot) \to \Hw_{-\nu}(\R,\mathscr{H}_0)
  \end{equation*}
  is a bounded linear operator for $\nu < \nu_1$.
  Indeed, since $\epsilon$ is bounded on $\C_{\Re>-\nu_1}$, so is $\epsilon_{01}$, and with Lemma \ref{lem:proj_positive} also $\epsilon_{11}(\cdot)^{-1}$. Hence,
  \begin{equation*}
    z \mapsto z\epsilon_{01}(z) \left( z\epsilon_{11}(z) \right)^{-1} = \epsilon_{01}(z) \left( \epsilon_{11}(z) \right)^{-1}
  \end{equation*}
  is bounded on $\C_{\Re > -\nu_1}$ as desired.
  Now apply this operator to the second equation in \eqref{eq:big_proj_system} and subtract the result from the first equation to obtain
  \begin{multline}
    \Biggl[
      \partial_t
      \begin{pmatrix}
        \partial_t\bigl( \epsilon_{00}(\partial_t) - \epsilon_{01}(\partial_t)\left( \epsilon_{11}(\partial_t) \right)^{-1}\epsilon_{10}(\partial_t) \bigr) & 0 \\ \partial_t\epsilon_{10}(\partial_t) & \partial_t\epsilon_{11}(\partial_t)
      \end{pmatrix} \\
      +
      \begin{pmatrix}
        \iota_0^*\Curl\mu^{-1}\Curl_0\iota_0 & 0 \\ 0 & 0
      \end{pmatrix}
    \Biggr]
    \begin{pmatrix}
      E_0 \\ E_1
    \end{pmatrix} =
    \begin{pmatrix}
      g_0 - \partial_t\epsilon_{01}(\partial_t)\left( \partial_t\epsilon_{11}(\partial_t) \right)^{-1}g_1 \\ g_1
    \end{pmatrix}.
    \label{eq:big_proj_system_shifted}
  \end{multline}
  Regarding the equation for $E_0$, we have that $\iota_0^*\Curl\mu^{-1}\Curl_0\iota_0$ is selfadjoint, continuously invertible, and nonnegative on $\mathscr{H}_0$.
  In view of Remark \ref{rem:Projections_invertible_selfadjoint_op} there exists a densely defined and boundedly invertible operator $C$ such that
  \begin{equation*}
    C^*C = \iota_0^*\Curl\mu^{-1}\Curl_0\iota_0.
  \end{equation*}
  Hence, the equation for $E_0$ is of the form \eqref{eq:evo_2order},
  \begin{equation*}
    \bigl(\partial_t^2\tilde{\epsilon}(\partial_t) + C^*C \bigr) E_0 = \tilde{g},
  \end{equation*}
  with $\tilde{g} = g_0 - \partial_t\epsilon_{01}(\partial_t)\left( \partial_t\epsilon_{11}(\partial_t) \right)^{-1}g_1$ and
  $\tilde{\epsilon}(z) = \epsilon_{00}(z) - \epsilon_{01}(z)\left( \epsilon_{11}(z) \right)^{-1}\epsilon_{10}(z)$.
  We verify that $\tilde{\epsilon}$ satisfies the accretivity conditions of Theorem~\ref{thm:exp_stab_M_d}.
  Indeed, with Lemma~\ref{lem:proj_positive} we have $\Re z\tilde{\epsilon}(z) \ge c$ whenever $\Re z\epsilon(z) \ge c$, thus $\tilde{\epsilon}$ fulfills \nameref{itm:MatCond_accr1}.
  Furthermore, we find that $\tilde{\epsilon}(z) = M_0(z) + z^{-1}M_1(z)$ with
  \begin{align*}
    M_0(z) &= \epsilon_{0,00} - \epsilon_{0,01}\epsilon_{11}(z)^{-1}\epsilon_{0,10}\\
    M_1(z) &= z\bigl( \chi_{00}(z)
      - \epsilon_{0,01} \epsilon_{11}(z)^{-1} \chi_{10}(z)
      - \chi_{01}(z)\epsilon_{11}(z)^{-1}\epsilon_{0,10}
      - \chi_{01}(z)\epsilon_{11}(z)^{-1}\chi_{10}(z) \bigr),
  \end{align*}
  where, due to \nameref{itm:MatCond_accr2}, $M_0(\mult)$, $M_1(\mult)$ are uniformly bounded on $\C_{\Re > -\nu_1}$ and $\lim_{z\to 0} M_1(z) = 0$.
  Consequently, there exists ${\nu}_0 \in (0, \nu_1]$ such that
  \begin{equation*}
    \left( \partial_t^2\left( \epsilon_{00}(\partial_t) - \epsilon_{01}(\partial_t)\epsilon_{11}(\partial_t)^{-1}\epsilon_{10}(\partial_t) \right) + C^*C \right)^{-1}
  \end{equation*}
  maps $\Hw_{-{\nu}}(\R,\mathscr{H}_0)$, $\nu < \nu_0$, causally into itself.
  Therefore, if we take $\nu < {\nu}_0$, then
  \begin{equation*}
    \lVert E_0 \rVert_{-\nu,0} \le \bigl\lVert g_0 - \partial_t\epsilon_{01}(\partial_t)\left( \partial_t\epsilon_{11}(\partial_t) \right)^{-1}g_1 \bigr\rVert_{-\nu,0} \le K\lVert g \rVert_{-\nu,0}
  \end{equation*}
  with some $K > 0$ independent of $g$ and $E$.

  Turning to the equation for $E_1$ in \eqref{eq:big_proj_system_shifted},
  \begin{equation*}
    \partial_t^2\epsilon_{10}(\partial_t)E_0 + \partial_t^2\epsilon_{11}(\partial_t)E_1 = g_1,
  \end{equation*}
  we infer
  \begin{align*}
    E_1 &= \left( \partial_t\epsilon_{11}(\partial_t) \right)^{-1}\partial_t^{-1}g_1 - \left( \partial_t\epsilon_{11}(\partial_t) \right)^{-1}\partial_t\epsilon_{10}(\partial_t)E_0 \\
    &= \partial_t\left( \partial_t\epsilon_{11}(\partial_t) \right)^{-1}\partial_t^{-2}g_1 - \left( \partial_t\epsilon_{11}(\partial_t) \right)^{-1}\partial_t\epsilon_{10}(\partial_t)E_0.
  \end{align*}
  Since by assumption $\partial_t^{-2}g_1 = h \in \Hw_{-\nu}(\R,\mathscr{H}_0^\bot) \subseteq \Hw_{-\nu}(\R,\Hs)$, and since
  \begin{equation*}
    \partial_t\left( \partial_t\epsilon_{11}(\partial_t) \right)^{-1}\quad \text{and}\quad \left( \partial_t\epsilon_{11}(\partial_t) \right)^{-1}\partial_t\epsilon_{10}(\partial_t)
  \end{equation*}
  leave $\Hw_{-\nu}(\R,\Hs)$ invariant, we obtain
  \begin{equation*}
    \lVert E_1 \rVert_{-\nu,0}
    \lesssim \lVert \partial_t^{-2}g_1 \rVert_{-\nu,0} + \lVert E_0 \rVert_{-\nu,0}
    \lesssim \lVert \partial_t^{-2}g_1 \rVert_{-\nu,0} + \lVert g_0 \rVert_{-\nu,0}.
  \end{equation*}
  The assertion now follows due to the density of $C_c^\infty(\R,\Hs)$ in $L^2_{-\nu}(\R,\Hs)$.
\end{proof}

Since the proof relies on an application of Theorem~\ref{thm:exp_stab_M_d}, in view of Remark~\ref{rem:u_H1_2ord} in fact a stronger result is implied by Theorem~\ref{thm:Maxwell_ExpStab_SecondOrder}.

\begin{cor}
  \label{cor:Maxwell_ExpStab_SecondOrder_stronger}
  Under the assumptions of Theorem~\ref{thm:Maxwell_ExpStab_SecondOrder}, and with the notation in the proof, the following holds:
  \begin{equation*}
    \begin{aligned}
      E_0, \partial_tE_0, CE_0 &\in L^2_{-\nu}(\R,\mathscr{H}_0)\\
      \lVert E_0 \rVert_{-\nu,0}, \lVert \partial_tE_0 \rVert_{-\nu,0}, \lVert CE_0 \rVert_{-\nu,0} &\lesssim \lVert g \rVert_{-\nu,0}
    \end{aligned}
  \end{equation*}
  and
  \begin{equation*}
    \begin{aligned}
      E_1, \partial_tE_1 &\in L^2_{-\nu}(\R,\mathscr{H}_0^\bot)\\
      \lVert E_1 \rVert_{-\nu,0} &\lesssim \lVert g \rVert_{-\nu,0} + \lVert h \rVert_{-\nu,0}\\
      \lVert \partial_tE_1 \rVert_{-\nu,0}
      &\lesssim \lVert \partial_tE_0 \rVert_{-\nu,0} + \lVert \partial_t^{-1}g_1 \rVert_{-\nu,0}
      \lesssim \lVert g \rVert_{-\nu,0} + \lVert \partial_t^{-1}g_1 \rVert_{-\nu,0}.
    \end{aligned}
  \end{equation*}
\end{cor}

\begin{rem}
  Obtaining exponential decay for the $H$-field through a second-order system is not completely analogous, the problem being that the material law $\epsilon(\partial_t)$ is non-instantaneous and, due to jumps, in general does not commute with $\Curl$.
  We address this issue in Section \ref{sec:ExpStab_FirstOrder}.
\end{rem}

\subsubsection{Exponential stability of the nonlinear second-order system}
\label{sec:NonlinExpStab}

We want to use the results in the linearized case to obtain exponential stability for the nonlinear system \eqref{eq:Maxwell_nonlinear}.
However, the fixed-point argument we employed previously to obtain well-posedness in $L^2_{\varrho}$, for large $\varrho > 0$, cannot be repeated in $L^2_{-\nu}$ (for $\nu < \nu_0$).
This can be seen in \eqref{eq:sol_F_Lip}, where a large $\varrho > 0$ is needed to ensure the contraction property on $\Hw_\varrho$.
As we show next, this problem can be avoided if we restrict ourselves to small solutions.

For $\varepsilon > 0$ we denote by $B_\varepsilon^{\nu}(\R,\Hs) \coleq \{ u \in L^2_{-\nu}(\R,\Hs) : \lVert u \rVert_{-\nu,0} \le \varepsilon \}$ the closed $\varepsilon$-ball in $L^2_{-\nu}(\R,\Hs)$.

\begin{thm}[Small solutions of the nonlinear system in $\Hs$]
  \label{thm:NonlinFixL2}
  Under the conditions of Theorem~\ref{thm:Maxwell_ExpStab_SecondOrder}, assume that for some $\tilde{\nu} > 0$,
  \begin{equation*}
    g = \partial_t \varPhi + \Curl\mu^{-1}\varPsi \in \Hw_{-\tilde{\nu}}(\R,\Hs)\quad
    \text{and}\quad h = \Pi_0\partial_t^{-1}\varPhi \in \Hw_{-\tilde{\nu}}(\R,\ker(\Curl_0)),
  \end{equation*}
  and suppose there exists $C > 0$ such that each $F \in \{ \partial_t^2P_{\mathrm{nl}}, \Pi_0P_{\mathrm{nl}} \}$ satisfies $F(0) = 0$ and
  \begin{equation}
    \forall u,v \in \Hw_{-\nu}(\R,\Hs) :\quad
    \lVert F(u) - F(v) \rVert_{-\nu,0}
      \le C\bigl( \lVert u \rVert_{-\nu,0} + \lVert v \rVert_{-\nu,0} \bigr) \lVert u - v \rVert_{-\nu,0}
    \label{eq:LocLipL2H0}
  \end{equation}
  for all $\nu \le \tilde{\nu}$.
  Then there exist $\nu_0, \varepsilon_0, c_0 > 0$ such that for all $\varepsilon \in (0,\varepsilon_0)$, $\nu < \nu_0$, the following holds:
  If $g, h \in B_{c_0\varepsilon}^\nu(\R,\Hs)$, then the nonlinear second-order system
  \begin{equation*}
    \bigl( \partial_t^2\epsilon(\partial_t) + \Curl\mu^{-1}\Curl_0 \bigr)E = -\partial_t^2P_{\mathrm{nl}}(E) + g
  \end{equation*}
  admits a unique solution $E \in B^{\nu}_{\varepsilon}(\R,\Hs)$.
\end{thm}

\begin{proof}
  By Theorem~\ref{thm:Maxwell_ExpStab_SecondOrder} there exists $\nu_0 \in (0, \tilde{\nu}]$ such that for all $\nu < \nu_0$ the solution operator
  \begin{equation*}
    T \coleq \bigl(\overline{\partial_t^2\epsilon(\partial_t) + \Curl\mu^{-1}\Curl_0}\bigr)^{-1}
  \end{equation*}
  for the linearized second-order system maps the subspace
  \begin{equation*}
    \{ f \in \Hw_{-\nu}(\R,\Hs) : \Pi_0\partial_t^{-2}f \in \Hw_{-\nu}(\R,\ker(\Curl_0)) \}
  \end{equation*}
  into $\Hw_{-\nu}(\R,\Hs)$.
  We show that if $\varepsilon, c_0 > 0$ are small, then the map $S$ given by
  \begin{equation*}
    S(u) \coleq T(g - \partial_t^2P_{\mathrm{nl}}(u))
  \end{equation*}
  is a contraction on $B^{\nu}_{\varepsilon}(\R,\Hs)$ if $g \in B^{\nu}_{c_0\varepsilon}(\R,\Hs)$ and $h = \Pi_0\partial_t^{-2}g \in B^{\nu}_{c_0\varepsilon}(\R,\Hs)$.
  Indeed, for fixed $u \in B^{\nu}_\varepsilon(\R,\Hs)$ we have with \eqref{eq:LocLipL2H0} and the estimate in Theorem~\ref{thm:Maxwell_ExpStab_SecondOrder} (replacing $\varPhi$ by $\varPhi - \partial_tP_{\mathrm{nl}}(u)$),
  \begin{align*}
    \lVert S(u) \rVert_{-\nu,0}
    &\le K\bigl( \lVert \partial_t^2P_{\mathrm{nl}}(u) \rVert_{-\nu,0} + \lVert \Pi_0P_{\mathrm{nl}}(u) \rVert_{-\nu,0} + \lVert g \rVert_{-\nu,0} + \lVert h \rVert_{-\nu,0} \bigr) \\
    &\le K\bigl( 2C\lVert u \rVert_{-\nu,0}^2 + \lVert g \rVert_{-\nu,0} + \lVert h \rVert_{-\nu,0} \bigr) \\
    &\le 2KC\varepsilon^2 + 2Kc_0\varepsilon \\
    &= \varepsilon (2KC\varepsilon + 2Kc_0),
  \end{align*}
  and similarly,
  \begin{align*}
    \lVert S(u) - S(v) \rVert_{-\nu,0}
    &= \lVert T(\partial_t^2P_{\mathrm{nl}}(u) - \partial_t^2P_{\mathrm{nl}}(v) \bigr) \rVert_{-\nu,0} \\
    &\le K \bigl( \lVert \partial_t^2P_{\mathrm{nl}}(u) - \partial_t^2P_{\mathrm{nl}}(v) \rVert_{-\nu,0} + \lVert \Pi_0P_{\mathrm{nl}}(u) - \Pi_0P_{\mathrm{nl}}(v) \rVert_{-\nu,0} \bigr) \\
    &\le 2KC \bigl( \lVert u \rVert_{-\nu,0} + \lVert v \rVert_{-\nu,0} \bigr) \lVert u - v \rVert_{-\nu,0} \\
    &\le 4KC\varepsilon \lVert u - v \rVert_{-\nu,0}
  \end{align*}
  Setting $\varepsilon < \varepsilon_0 \coleq \frac{1}{4KC}$ and $c_0 \coleq \frac{1}{4K}$, the claim follows.
\end{proof}

\begin{ex}
  \label{ex:NonlocalExpStab}
  Consider the fully nonlocal quadratic polarization $P_\mathrm{nl}$ in Example \ref{ex:NonlocalNonlinearity} given by
  \begin{equation*}
    P_\mathrm{nl}(E)(t) = \iint \kappa(t-s_1,t-s_2)\,\snl(E(s_1),E(s_2))\dd{s_1}\dd{s_2}
  \end{equation*}
  with $\snl(u,v) = \int_\Omega \int_\Omega \Lambda(\,\mult\,,y_1,y_2)u(y_1)v(y_2)\dd{y_1}\dd{y_2}$, where $\Lambda \in L^2(\Omega^3)^{3\times 3\times 3}$.
  For simplicity we assume again that $\kappa$ is compactly supported in $(0,\infty)^2$, say $\kappa \in C_c^\infty((0,\infty)^2)^{3\times 3}$, and have thus
  \begin{equation*}
    \partial_t^\ell P_\mathrm{nl}(E)(t) = \iint (\partial_1 + \partial_2)^\ell\kappa(t-s_1,t-s_2)\,\snl(E(s_1),E(s_2))\dd{s_1}\dd{s_2},\quad \ell \in \N,
  \end{equation*}
  with $K_\ell \coleq (\partial_1 + \partial_2)^\ell\kappa \in C_c^\infty((0,\infty)^2)^{3\times 3} \subseteq L^1_{\varrho}(\R^2)^{3\times 3}$ for all $\varrho \in \R$.
  In view of Section~\ref{sec:initial_values} we can assume that $\supp E, \supp H \subseteq [0,\infty)$.
  Since each $\partial_t^\ell P_\mathrm{nl}$ maps $\Hw_{-\nu}(\R,\Hs)$ causally into $\Hw_{-2\nu}(\R,\Hs)$, we have again $\partial_t^\ell P_\mathrm{nl}(E) \in \Hw_{-\nu}(\R,\Hs)$.
  In particular, estimate \eqref{eq:LocLipL2H0} holds for $F = \Pi_0P_\mathrm{nl}$ and $F = \partial_t^2P_\mathrm{nl}$.
\end{ex}

\subsection{Exponential stability of the first-order system}
\label{sec:ExpStab_FirstOrder}

To obtain exponential decay of the $H$-field, we now consider exponential stability for the first-order system.
In fact, the next result is a refinement of Theorem~\ref{thm:Maxwell_ExpStab_SecondOrder}, and also relies on the second-order formulation \eqref{eq:Mw2ordE}.

Denote by $\Pi_0\colon \Hs \to \ker(\Curl_0)$ and $\Pi_1\colon \Hs \to \ker(\Curl)$ the canonical projections.

\begin{thm}[Exponential stability of the linear Maxwell equations]
  \label{thm:Maxwell_ExpStab_FirstOrder}
  Let $\Omega$ be a bounded weak Lipschitz domain, let $\epsilon(\partial_t) = \epsilon_0 + \chi(\partial_t)$ be a linear material law satisfying \nameref{itm:MatCond_accr1} and \nameref{itm:MatCond_accr2}.
  Furthermore, let $\mu$ be selfadjoint and uniformly positive definite.
  Then there exists $\tilde{\nu} > 0$ such that if $\nu < \tilde{\nu}$ and
  \begin{align*}
    \varPhi, \varPsi \in \Hw_{-\nu}(\R,\Hs)\\
    g \coleq \partial_t\varPhi + \Curl\mu^{-1}\varPsi \in L^2_{-\nu}(\R,\Hs)\\
    h \coleq \Pi_0\partial_t^{-1}\varPhi \in \Hw_{-\nu}(\R,\ker(\Curl_0))\\
    w \coleq \Pi_1\partial_t^{-1}\varPsi \in \Hw_{-\nu}(\R,\ker(\Curl)),
  \end{align*}
  then the solution $E,H$ of
  \begin{equation}
    \left( \partial_t
    \begin{pmatrix}
      \epsilon(\partial_t) & 0 \\ 0 & \mu
    \end{pmatrix} +
    \begin{pmatrix}
      0 & -\Curl \\ \Curl_0 & 0
    \end{pmatrix} \right)
    \begin{pmatrix}
      E \\ H
    \end{pmatrix}
    =
    \begin{pmatrix}
      \varPhi \\ \varPsi
    \end{pmatrix}
    \label{eq:Maxwell_linear_Drude}
  \end{equation}
  satisfies $E,H \in H^1_{-\nu}(\R,\Hs)$ and
  \begin{equation}
    \begin{aligned}
      \lVert E \rVert_{-\nu,0} &\lesssim \lVert g \rVert_{-\nu,0} + \lVert h \rVert_{-\nu,0}\\
      \lVert H \rVert_{-\nu,0} &\lesssim \lVert g \rVert_{-\nu,0} + \lVert \varPhi \rVert_{-\nu,0} + \lVert w \rVert_{-\nu,0}\\
      \lVert \partial_tE \rVert_{-\nu,0} &\lesssim \lVert g \rVert_{-\nu,0} + \lVert \varPhi \rVert_{-\nu,0}\\
      \lVert \partial_tH \rVert_{-\nu,0} &\lesssim \lVert g \rVert_{-\nu,0} + \lVert \varPsi \rVert_{-\nu,0}.
    \end{aligned}
    \label{eq:EH_H1_estimates}
  \end{equation}

\end{thm}

\begin{proof}
  We use projections similar to those in the proof of Theorem~\ref{thm:Maxwell_ExpStab_SecondOrder}, now for the full system \eqref{eq:Maxwell_linear_Drude}.
  Consider
  \begin{equation*}
    \begin{aligned}
      \mathscr{H}_0 &\coleq \ran(\Curl) = \ker(\Curl_0)^\bot & \tikz{\draw[{Hooks[right]}->] (0,0)--(1,0) node[midway,above]{$\iota_0$};}\quad \Hs\\
      \mathscr{H}_0^\bot &= \ran(\Curl)^\bot = \ker(\Curl_0) & \tikz{\draw[{Hooks[right]}->] (0,0)--(1,0) node[midway,above]{$\iota_1$};}\quad \Hs\\
      \mathscr{H}_1 &\coleq \ran(\Curl_0) = \ker(\Curl)^\bot & \tikz{\draw[{Hooks[right]}->] (0,0)--(1,0) node[midway,above]{$\tau_0$};}\quad \Hs\\
      \mathscr{H}_1^\bot &= \ran(\Curl_0)^\bot = \ker(\Curl) & \tikz{\draw[{Hooks[right]}->] (0,0)--(1,0) node[midway,above]{$\tau_1$};}\quad \Hs
    \end{aligned}
  \end{equation*}
  (note that $\Pi_0 = \iota_1^*$, $\Pi_1 = \tau_1^*$)
  and observe that each of the maps $\tau_1^*\Curl_0$, $\iota_1^*\Curl$, $\Curl\tau_1$, $\Curl_0\iota_1$ is zero on its corresponding domain.
  Now set for $i,j \in \{0,1\}$
  \begin{equation*}
  \begin{aligned}
    \epsilon_{ij}(\partial_t) &\coleq \iota_i^*\epsilon(\partial_t)\iota_j  & E_i &\coleq \iota_i^*E & \varPhi_i &\coleq \iota_i^* \varPhi\\
    \mu_{ij} &\coleq \tau_i^*\mu\tau_j & H_i &\coleq \tau_i^*H & \varPsi_i &\coleq \tau_i^* \varPsi.
  \end{aligned}
  \end{equation*}
  Then \eqref{eq:Maxwell_linear_Drude} can be written in the form
  \begin{equation}
    \left( \partial_t
    \begin{pmatrix}
      \epsilon_{00}(\partial_t) & \epsilon_{01}(\partial_t) & 0 & 0\\ \epsilon_{10}(\partial_t) & \epsilon_{11}(\partial_t) & 0 & 0\\
      0 & 0 & \mu_{00} & \mu_{01} &\\  0 & 0 & \mu_{10} & \mu_{11}
    \end{pmatrix} +
    \begin{pmatrix}
      0 & 0 & -\iota_0^*\Curl\tau_0 & 0\\
      0 & 0 & 0 & 0\\
      \tau_0^*\Curl_0\iota_0 & 0 & 0 & 0\\
      0 & 0 & 0 & 0 
    \end{pmatrix} \right)
    \begin{pmatrix}
      E_0 \\ E_1 \\ H_0 \\ H_1
    \end{pmatrix} =
    \begin{pmatrix}
      \varPhi_0\\ \varPhi_1\\ \varPsi_0\\ \varPsi_1
    \end{pmatrix}.
    \label{eq:MaxwellProj_bigmatrix}
  \end{equation}
  Here we note that the operator $B \coleq \tau_0^*\Curl_0\iota_0\colon \dom(B) \subseteq \mathscr{H}_0 \to \mathscr{H}_1$ is injective, surjective, and closed, and thus continuously invertible by the closed graph theorem; the same is true for its adjoint $B^* = \iota_0^*\Curl\tau_0$.

  Solving the second and fourth equation in \eqref{eq:MaxwellProj_bigmatrix} for
  \begin{equation}
    \begin{split}
      E_1 &= \epsilon_{11}(\partial_t)^{-1}\partial_t^{-1}\varPhi_1 - \epsilon_{11}(\partial_t)^{-1}\epsilon_{10}(\partial_t)E_0\\
      H_1 &= \mu_{11}^{-1}\partial_t^{-1}\varPsi_1 - \mu_{11}^{-1}\mu_{10}H_0,
    \end{split}
    \label{eq:E1H1_aux}
  \end{equation}
  we obtain the system
  \begin{equation}
    \left( \partial_t
    \begin{pmatrix}
      \tilde{\epsilon}(\partial_t) & 0 \\
      0 & \tilde{\mu}
    \end{pmatrix} +
    \begin{pmatrix}
      0 & -B^*\\
      B & 0
    \end{pmatrix} \right)
    \begin{pmatrix}
      E_0 \\ H_0
    \end{pmatrix} =
    \begin{pmatrix}
      \tilde{\varPhi}\\ \tilde{\varPsi}
    \end{pmatrix},
    \label{eq:Maxwell_linear_Drude_proj}
  \end{equation}
  where
  \begin{equation*}
    \begin{aligned}
      \widetilde{\epsilon}(\partial_t) &\coleq \epsilon_{00}(\partial_t) - \epsilon_{01}(\partial_t)\epsilon_{11}(\partial_t)^{-1}\epsilon_{10}(\partial_t)\\
      \tilde{\mu} &\coleq \mu_{00} - \mu_{01}\mu_{11}^{-1}\mu_{10},
    \end{aligned}
    \qquad
    \begin{aligned}
      \tilde{\varPhi} &\coleq \varPhi_0 - \epsilon_{01}(\partial_t)\epsilon_{11}(\partial_t)^{-1}\varPhi_1\\
      \tilde{\varPsi} &\coleq \varPsi_0 - \mu_{01}\mu_{11}^{-1}\varPsi_1.
    \end{aligned}
  \end{equation*}
  Here $\epsilon_{11}(\partial_t)^{-1}, \tilde{\epsilon}(\partial_t)$, $\mu_{11}^{-1}, \tilde{\mu}$ are bounded in $\Hw_{-\nu}$ for $\nu < \nu_1$ by assumptions on $\epsilon, \mu$ and Lemma~\ref{lem:proj_positive}.
  Moreover, $\lVert\tilde{\varPhi}\rVert_{-\nu,0} \lesssim \lVert\varPhi\rVert_{-\nu,0}$ and $\lVert\tilde{\varPsi}\rVert_{-\nu,0} \lesssim \lVert\varPsi\rVert_{-\nu,0}$.

  The assertions for $E$ follow from Theorem \ref{thm:Maxwell_ExpStab_SecondOrder} (the notation $E_0, E_1$ is the same as in the proof), since the data fulfills the necessary conditions.
  By Corollary~\ref{cor:Maxwell_ExpStab_SecondOrder_stronger} we even have
  $E_0,\partial_tE_0 \in \Hw_{-\nu}(\R,\mathscr{H}_0)$ and
  $CE_0 \in \Hw_{-\nu}(\R,\mathscr{H}_1)$ for $C$ given by $B^*\mu^{-1}B = C^*C$, and the latter also implies $BE_0 \in \Hw_{-\nu}(\R,\mathscr{H}_1)$ since $\mu$ is selfadjoint and boundedly invertible.
  Moreover, the estimates
  \begin{equation*}
    \begin{aligned}
      \lVert E_0 \rVert_{-\nu,0}, \lVert \partial_tE_0 \rVert_{-\nu,0}, \lVert BE_0 \rVert_{-\nu,0} &\lesssim \lVert g \rVert_{-\nu,0}
    \end{aligned}
  \end{equation*}
  hold.
  Using \eqref{eq:E1H1_aux} now provides
  \begin{equation*}
    \begin{aligned}
      \lVert E_1 \rVert_{-\nu,0} &= \lVert \epsilon_{11}(\partial_t)^{-1}\partial_t^{-1}\varPhi_1 - \epsilon_{11}(\partial_t)^{-1}\epsilon_{10}(\partial_t)E_0\rVert_{-\nu,0}
      \lesssim \lVert \partial_t^{-1}\varPhi_1 \rVert_{-\nu,0} + \lVert E_0 \rVert_{-\nu,0}\\
      \lVert \partial_tE_1 \rVert_{-\nu,0} &= \lVert \epsilon_{11}(\partial_t)^{-1}\varPhi_1 - \epsilon_{11}(\partial_t)^{-1}\epsilon_{10}(\partial_t)\partial_tE_0\rVert_{-\nu,0}
      \lesssim \lVert \varPhi_1 \rVert_{-\nu,0} + \lVert \partial_tE_0 \rVert_{-\nu,0},
    \end{aligned}
  \end{equation*}
  and overall
  \begin{equation}
    \begin{aligned}
      \lVert E \rVert_{-\nu,0} &\lesssim \lVert g \rVert_{-\nu,0} + \lVert \partial_t^{-1}\varPhi_1 \rVert_{-\nu,0}
        = \lVert g \rVert_{-\nu,0} + \lVert h \rVert_{-\nu,0}\\
      \lVert \partial_tE \rVert_{-\nu,0} &\lesssim \lVert g \rVert_{-\nu,0} + \lVert \varPhi_1 \rVert_{-\nu,0}
        \lesssim \lVert g \rVert_{-\nu,0} + \lVert \varPhi \rVert_{-\nu,0}.
    \end{aligned}
    \label{eq:E_H1_estimates}
  \end{equation}

  To obtain the assertions for $H$, we first take the first line in \eqref{eq:Maxwell_linear_Drude_proj} and have
  \begin{equation*}
    H_0 = (B^*)^{-1} \bigl(\partial_t\tilde{\epsilon}(\partial_t)E_0 - \tilde{\varPhi}\bigr) \in \Hw_{-\nu}(\R,\mathscr{H}_1),
  \end{equation*}
  by boundedness of $(B^*)^{-1}$, with
  \begin{equation*}
    \begin{aligned}
      \lVert H_0 \rVert_{-\nu,0}
      &\lesssim \lVert \partial_tE_0\rVert_{-\nu,0} + \lVert \tilde{\varPhi} \rVert_{-\nu,0}
      \lesssim \lVert g\rVert_{-\nu,0} + \lVert \varPhi \rVert_{-\nu,0},
    \end{aligned}
  \end{equation*}
  and the second line gives
  \begin{equation*}
    \partial_tH_0 = -\tilde{\mu}^{-1}BE_0 + \tilde{\mu}^{-1}\tilde{\varPsi} \in \Hw_{-\nu}(\R,\mathscr{H}_1),
  \end{equation*}
  with
  \begin{equation*}
    \lVert \partial_tH_0 \rVert_{-\nu,0} \lesssim \lVert BE_0 \rVert_{-\nu,0} + \lVert \tilde{\varPsi} \rVert_{-\nu,0}
  \lesssim \lVert g \rVert_{-\nu,0} + \lVert \varPsi \rVert_{-\nu,0}.
  \end{equation*}
  The corresponding statement for $H_1$ follows using \eqref{eq:E1H1_aux}, since
  \begin{align*}
    H_1 &= \mu_{11}^{-1}\partial_t^{-1}\varPsi_1 - \mu_{11}^{-1}\mu_{10}H_0 \in \Hw_{-\nu}(\R,\mathscr{H}_1^\bot)\\
    \partial_tH_1 &= \mu_{11}^{-1}\varPsi_1 - \mu_{11}^{-1}\mu_{10}\partial_tH_0 \in \Hw_{-\nu}(\R,\mathscr{H}_1^\bot),
  \end{align*}
  with
  \begin{equation*}
    \begin{aligned}
      \lVert H_1\rVert_{-\nu,0} &\lesssim \lVert H_0\rVert_{-\nu,0} + \lVert\partial_t^{-1}\varPsi_1 \rVert_{-\nu,0}\\
      &\lesssim  \lVert g\rVert_{-\nu,0} + \lVert \varPhi \rVert_{-\nu,0} + \lVert\partial_t^{-1}\varPsi_1 \rVert_{-\nu,0}\\
      \lVert \partial_tH_1\rVert_{-\nu,0}
      &\lesssim \lVert \partial_tH_0 \rVert_{-\nu,0} + \lVert \varPsi_1 \rVert_{-\nu,0}\\
      &\lesssim \lVert g \rVert_{-\nu,0} + \lVert \varPsi \rVert_{-\nu,0}.
    \end{aligned}
  \end{equation*}
  We obtain
  \begin{equation}
    \begin{aligned}
      \lVert H \rVert_{-\nu,0}
      &\lesssim \lVert g\rVert_{-\nu,0} + \lVert \varPhi\rVert_{-\nu,0} + \lVert \partial_t^{-1}\varPsi_1\rVert_{-\nu,0}\\
      \lVert \partial_tH \rVert_{-\nu,0}
      &\lesssim \lVert g\rVert_{-\nu,0} + \lVert \varPsi\rVert_{-\nu,0}.
    \end{aligned}
    \label{eq:H_H1_estimates}
  \end{equation}
  Now \eqref{eq:E_H1_estimates} and \eqref{eq:H_H1_estimates} imply \eqref{eq:EH_H1_estimates}.
\end{proof}

\subsubsection{Exponential stability in \texorpdfstring{$\boldsymbol{H^2}$}{H²} for materials constant on each \texorpdfstring{$\boldsymbol{\Omega_i}$}{Ωi}}
\label{sec:ExpStab_H2}

We now study exponential stability in the context of higher spatial regularity.
This is motivated by the fact that the spaces $\Hs^k$, $k\ge 2$, enjoy the algebra property which is useful when polynomial nonlinearities are present.
We will combine the assertion of Theorem~\ref{thm:Maxwell_ExpStab_FirstOrder} with the results in Section~\ref{sec:higher_spatial_regularity}. Thus, we will only consider material laws that satisfy \nameref{itm:MatCond_split}.
\noindent Let us introduce a more succinct notation for admissible data $\varPhi, \varPsi$.
\begin{defi}
  \label{defi:AdmissibleData_VW}
  For $\nu > 0$, the spaces $\mathcal{V}_{-\nu}$ and $\mathcal{W}_{-\nu}$ are defined by
  \begin{equation*}
  \begin{split}
    \varPhi \in \mathcal{V}_{-\nu} &\iff
    \left\{
    \begin{aligned}
      \varPhi &\in H^2_{-\nu}(\R,\Hs) \cap L^2_{-\nu}(\R,\spatH{1} \cap H(\Div,\Omega))\\
      \partial_t^{-1}\varPhi &\in L^2_{-\nu}(\R,H(\Div,\Omega))\\
      \Div\partial_t^{-1}\varPhi &\in L^2_{-\nu}(\R,\spatH{1})
    \end{aligned}
    \right.\\
    \varPsi \in \mathcal{W}_{-\nu} &\iff
    \left\{
    \begin{aligned}
      \varPsi &\in H^2_{-\nu}(\R,\Hs) \cap L^2_{-\nu}(\R,\spatH{1} \cap H_0(\Div,\Omega))\\
      \partial_t^{-1}\varPsi &\in L^2_{-\nu}(\R,H_0(\Div,\Omega))\\
      \Div\partial_t^{-1}\varPsi &\in L^2_{-\nu}(\R,\spatH{1})\\
      \Curl\mu^{-1}\varPsi &\in H^1_{-\nu}(\R,\Hs)
    \end{aligned}
  \right.
  \end{split}
  \end{equation*}
  and are equipped with the norms
  \begin{align*}
    \lVert \varPhi \rVert_{\mathcal{V}_{-\nu}}
    &= \lVert \varPhi \rVert_{-\nu,1}
      + \sum_{j=1,2}\lVert \partial_t^j\varPhi \rVert_{-\nu,0}
      + \lVert \Pi_0\partial_t^{-1}\varPhi \rVert_{-\nu,0}
      + \lVert \Div\partial_t^{-1}\varPhi \rVert_{-\nu,1}\\
    \lVert \varPsi \rVert_{\mathcal{W}_{-\nu}}
    &= \lVert \varPsi \rVert_{-\nu,1}
      + \sum_{j=1,2}\lVert \partial_t^j\varPsi \rVert_{-\nu,0}
      + \lVert \Pi_1\partial_t^{-1}\varPsi \rVert_{-\nu,0}
      + \lVert \Div\partial_t^{-1}\varPsi \rVert_{-\nu,1}
      + \lVert \Curl\mu^{-1}\partial_t\varPsi \rVert_{-\nu,0}.
  \end{align*}
\end{defi}

\noindent Note that the conditions $\Pi_0\partial_t^{-1}\varPhi \in \Hw_{-\nu}(\R,\ker(\Curl_0))$ and $\Pi_1\partial_t^{-1}\varPsi \in \Hw_{-\nu}(\R,\ker(\Curl))$ of Theorem~\ref{thm:Maxwell_ExpStab_FirstOrder} are indeed fulfilled if $\varPhi \in \mathcal{V}_{-\nu}$, $\varPsi \in \mathcal{W}_{-\nu}$.

\begin{thm}[Exponential stability of the linear Maxwell equations in $H^2$]
  \label{thm:Maxwell_ExpStab_H2}
  Let $\Omega = \Omega_1\sqcup\Gamma\sqcup\Omega_2$ be a bounded Lipschitz domain satisfying the regularity assumptions of Proposition \ref{prop:spatial_regularity} with $k = 2$.
  Let $\epsilon(\partial_t),\mu$ be of the form in \nameref{itm:MatCond_split} and let $\epsilon(\partial_t)$ satisfy \nameref{itm:MatCond_accr1} and \nameref{itm:MatCond_accr2} with some $\nu_1 > 0$.
  Under the assumptions on $\varPhi, \varPsi$ of Theorem~\ref{thm:Maxwell_spatial_regularity}, suppose in addition for $\nu < \nu_1$ that $\varPhi \in \mathcal{V}_{-\nu}$, $\varPsi \in \mathcal{W}_{-\nu}$ and
  \begin{equation*}
    \ph \in L^2_{-\nu}(\R,\spatH{2}) \quad \text{with} \quad \partial_t^j\ph\in \Hw_{-\nu}(\R,\Hs) \quad \text{for all } j \in \{-1, 0, 1, 2\}.
  \end{equation*}
  Then there exists $\tilde{\nu} \in (0,\nu_1]$ such that the solution $(E,H)$ of the system
  \begin{equation*}
    \left( \partial_t
    \begin{pmatrix}
      \epsilon(\partial_t) & 0 \\ 0 & \mu
    \end{pmatrix} +
    \begin{pmatrix}
      0 & -\Curl \\ \Curl_0 & 0
    \end{pmatrix} \right)
    \begin{pmatrix}
      E \\ H
    \end{pmatrix}
    =
    \begin{pmatrix}
      \varPhi - \ph \\ \varPsi
    \end{pmatrix}
  \end{equation*}
  fulfills $E,H \in L^2_{-\nu}(\R,\spatH{2})$ for $0 < \nu < \tilde{\nu}$ and
  \begin{equation}
    \lVert E \rVert_{-\nu,2}, \lVert H \rVert_{-\nu,2} \lesssim
    \lVert \varPhi \rVert_{\mathcal{V}_{-\nu}}
    + \lVert \varPsi \rVert_{\mathcal{W}_{-\nu}}
    + \lVert \ph \rVert_{-\nu,2} + \sum_{j=-1}^2\lVert \partial_t^j\ph \rVert_{-\nu,0}.
    \label{eq:ExpStabH2_estimate}
  \end{equation}
\end{thm}

\begin{proof}
  The claim follows from similar estimates as in the proof of Theorem~\ref{thm:Maxwell_spatial_regularity}.
  Let $0 < \nu < \tilde{\nu}$.
  Then, we have \eqref{eq:DivCurl_H2_estimate} with $\varrho = -\nu$, \ie,
  \begin{equation}
    \begin{aligned}
      \lVert E \rVert_{-\nu,2} &\lesssim \bigl(1 + \sup_{\Re z > -\nu_1}\lVert \epsilon(z)^{-1} \rVert \bigr) \bigl( \lVert E \rVert_{-\nu,0} + \lVert \Curl_0 E \rVert_{-\nu,1} + \lVert \Div \epsilon(\partial_t) E \rVert_{-\nu,1} + \lVert \ph \rVert_{-\nu,2}\bigr)\\
      \lVert H \rVert_{-\nu,2} &\lesssim \bigl(1 + \lVert \mu^{-1} \rVert \bigr) \bigl( \lVert H \rVert_{-\nu,0} + \lVert \Curl H \rVert_{-\nu,1} + \lVert \Div \mu H \rVert_{-\nu,1} \bigr).
    \end{aligned}
    \label{eq:H2_epsmu_estimate}
  \end{equation}
  We will estimate the norms on the right-hand side recursively.

  By assumption on $\ph$ and since $\varPhi\in \mathcal{V}_{-\nu}$, $\varPsi\in \mathcal{W}_{-\nu}$, Theorem~\ref{thm:Maxwell_ExpStab_FirstOrder} applies (with $\varPhi$ replaced by $\varPhi - \ph$) and yields $E,\partial_t E, H, \partial_t H \in \Hw_{-\nu}(\R,\Hs)$.
  Moreover, the additional assumptions on the regularity in time of the data imply that Theorem~\ref{thm:Maxwell_ExpStab_FirstOrder} can be applied once more, with $\partial_t\varPhi, \partial_t\ph, \partial_t\varPsi$ instead of $\varPhi, \ph, \varPsi$, yielding $\partial_tE, \partial_tH, \partial_t^2E, \partial_t^2H \in \Hw_{-\nu}(\R,\Hs)$.
  From \eqref{eq:EH_H1_estimates} we obtain the following estimates
  \begin{align*}
    \lVert E \rVert_{-\nu,0} &\lesssim \lVert \partial_t(\varPhi - \ph) + \Curl\mu^{-1}\varPsi \rVert_{-\nu,0}
    + \lVert \Pi_0\partial_t^{-1}(\varPhi - \ph) \rVert_{-\nu,0}\\
    \lVert H \rVert_{-\nu,0} &\lesssim \lVert \partial_t(\varPhi - \ph) + \Curl\mu^{-1}\varPsi \rVert_{-\nu,0}
    + \lVert \varPhi - \ph \rVert_{-\nu,0} + \lVert \Pi_1\partial_t^{-1}\varPsi \rVert_{-\nu,0}\\
    \lVert \partial_tE \rVert_{-\nu,0} &\lesssim \lVert \partial_t(\varPhi - \ph) + \Curl\mu^{-1}\varPsi \rVert_{-\nu,0}
    + \lVert \varPhi - \ph \rVert_{-\nu,0}\\
    \lVert \partial_tH \rVert_{-\nu,0} &\lesssim \lVert \partial_t(\varPhi - \ph) + \Curl\mu^{-1}\varPsi \rVert_{-\nu,0}
    + \lVert \varPsi \rVert_{-\nu,0}\\
    \lVert \partial_t^2E \rVert_{-\nu,0} &\lesssim \lVert \partial_t^2(\varPhi - \ph) + \Curl\mu^{-1}\partial_t\varPsi \rVert_{-\nu,0}
    + \lVert \partial_t(\varPhi - \ph) \rVert_{-\nu,0}\\
    \lVert \partial_t^2H \rVert_{-\nu,0} &\lesssim \lVert \partial_t^2(\varPhi - \ph) + \Curl\mu^{-1}\partial_t\varPsi \rVert_{-\nu,0}
    + \lVert \partial_t\varPsi \rVert_{-\nu,0}.
  \end{align*}
  Collecting all the above terms, we observe that
  \begin{equation}
    \sum_{j=0}^2\bigl( \lVert \partial_t^jE \rVert_{-\nu,0} + \lVert \partial_t^jH \rVert_{-\nu,0}\bigr)
    \lesssim \lVert \varPhi \rVert_{\mathcal{V}_{-\nu}} + \lVert \varPsi \rVert_{\mathcal{W}_{-\nu}}
    + \sum_{j=-1}^2\lVert \partial_t^{j}\ph \rVert_{-\nu,0}.
    \label{eq:H2_data_base_estimate}
  \end{equation}
  Next, using the identities in the Maxwell system and Propositions \ref{prop:spatial_regularity} and \ref{prop:spatial_regularity_jump}, we have
  \begin{align*}
    \lVert \Curl_0E \rVert_{-\nu,1} &= \lVert -\partial_t\mu H + \varPsi \rVert_{-\nu,1}\\
    &\lesssim \lVert \partial_tH \rVert_{-\nu,0} + \lVert \Curl\partial_tH \rVert_{-\nu,0} + \lVert \Div\partial_t\mu H \rVert_{-\nu,0} + \lVert \varPsi \rVert_{-\nu,1}\\
    &\lesssim \lVert \partial_tH \rVert_{-\nu,0} + \lVert \epsilon(\partial_t)\partial_t^2E - \partial_t(\varPhi - \ph) \rVert_{-\nu,0} + \lVert \Div\varPsi \rVert_{-\nu,0} + \lVert \varPsi \rVert_{-\nu,1}\\
    &\lesssim \lVert \partial_tH \rVert_{-\nu,0} + \lVert \partial_t^2E \rVert_{-\nu,0} + \lVert \partial_t\varPhi \rVert_{-\nu,0} + \lVert \partial_t\ph \rVert_{-\nu,0} + \lVert \varPsi \rVert_{-\nu,1}\\
    \lVert \Div \partial_t\epsilon(\partial_t)E \rVert_{-\nu,1}
    &\lesssim \lVert \Div(\varPhi - \ph) \rVert_{-\nu,1}\\
    &\lesssim \lVert \Div\varPhi \rVert_{-\nu,1} + \lVert \ph \rVert_{-\nu,2}\\
    \lVert \Curl H \rVert_{-\nu,1}
    &= \lVert \partial_t\epsilon(\partial_t)E - \varPhi + \ph \rVert_{-\nu,1}\\
    &\lesssim \lVert \epsilon(\partial_t)\partial_tE \rVert_{-\nu,1} + \lVert \varPhi \rVert_{-\nu,1} + \lVert \ph \rVert_{-\nu,1}\\
    &\lesssim \lVert \partial_tE \rVert_{-\nu,1} + \lVert \varPhi \rVert_{-\nu,1} + \lVert \ph \rVert_{-\nu,1}\\
    &\lesssim \lVert \partial_tE \rVert_{-\nu,0} + \lVert \Curl_0 \partial_tE \rVert_{-\nu,0} + \lVert \Div \epsilon(\partial_t)\partial_tE \rVert_{-\nu,0} + \lVert \varPhi \rVert_{-\nu,1} + \lVert \ph \rVert_{-\nu,1}\\
    &\lesssim \lVert \partial_tE \rVert_{-\nu,0} + \lVert -\partial_t^2\mu H + \partial_t\varPsi \rVert_{-\nu,0} + \lVert \Div(\varPhi - \ph) \rVert_{-\nu,0} + \lVert \varPhi \rVert_{-\nu,1} + \lVert \ph \rVert_{-\nu,1}\\
    &\lesssim \lVert \partial_tE \rVert_{-\nu,0} + \lVert \partial_t^2H \rVert_{-\nu,0} + \lVert \partial_t\varPsi \rVert_{-\nu,0} + \lVert \varPhi \rVert_{-\nu,1} + \lVert \ph \rVert_{-\nu,1}\\
    \lVert \Div\mu H \rVert_{-\nu,1}
    &\lesssim \lVert \Div\partial_t^{-1}\varPsi \rVert_{-\nu,1}.
  \end{align*}
  Together with \eqref{eq:H2_data_base_estimate}, we see that all terms can be controlled by the right-hand side of \eqref{eq:ExpStabH2_estimate}. The claim now follows with \eqref{eq:H2_epsmu_estimate}.
\end{proof}

Now we can formulate an exponential stability result in the nonlinear $H^{2}$-setting.
We employ a fixed-point argument in $L^2_{-\nu}(\R,\spatH{2})$.

\begin{thm}[Small solutions of the nonlinear system in $H^{2}$]
  \label{thm:NonlinFixH2}
  Let $\Omega$, $\epsilon$, $\mu$, $\tilde{\nu}$ be given as in Theorem~\ref{thm:Maxwell_ExpStab_H2} and impose the regularity conditions of Theorem~\ref{thm:Maxwell_spatial_regularity} on the data $\varPhi,\varPsi$.
  Suppose now that the map $P_\mathrm{nl}$ is such that for $0 < \nu < \tilde{\nu}$ and for all $j\in\{0,1,2,3\}$,
  \begin{equation*}
    \begin{aligned}
      \partial_t^jP_\mathrm{nl} &\colon L^2_{-\nu}(\R,\spatH{2}) \to L^2_{-\nu}(\R,\spatH{2})
    \end{aligned}
  \end{equation*}
  are causal and satisfy the local Lipschitz estimates
  \begin{equation}
    \begin{aligned}
      \lVert \partial_t^jP_\mathrm{nl}(u) - \partial_t^jP_\mathrm{nl}(v) \rVert_{-\nu,2} &\le C \bigl( \lVert u \rVert_{-\nu,2} + \lVert v \rVert_{-\nu,2} \bigr)^\alpha \lVert u - v \rVert_{-\nu,2}
    \end{aligned}
    \label{eq:LocLipH2}
  \end{equation}
  for some $\alpha, C > 0$ and all $u,v \in \Hw_{-\nu}(\R,\spatH{2})$.
  Then there exist $\varepsilon_0,c_0 > 0$ such that for all $\varepsilon \in (0,\varepsilon_0]$, $\nu < \tilde{\nu}$, the following holds:
  If $\varPhi \in \mathcal{V}_{-\nu}$, $\varPsi \in \mathcal{W}_{-\nu}$ with
  $\lVert \varPhi \rVert_{\mathcal{V}_{-\nu}},\, \lVert \varPsi \rVert_{\mathcal{W}_{-\nu}} < c_0\varepsilon$, then the nonlinear Maxwell system
  \begin{equation*}
    \left( \partial_t\!
      \begin{pmatrix}
        \epsilon(\partial_t) & 0 \\ 0 & \mu
      \end{pmatrix} +
      \begin{pmatrix}
        0 & -\Curl \\ \Curl_0 & 0
      \end{pmatrix}
    \right)\! \binom{E}{H}
    =
    \begin{pmatrix}{\varPhi - \partial_tP_{\mathrm{nl}}(E)}\\{\varPsi}\end{pmatrix}
  \end{equation*}
  admits a unique solution $E,H \in L^2_{-\nu}(\R,\spatH{2})$ with $\lVert E \rVert_{-\nu,2}, \lVert H \rVert_{-\nu,2} \le \varepsilon$.
\end{thm}

\begin{proof}
  Define the solution map $S$ by
  \begin{equation*}
    S(u) \coleq \left(\overline{\partial_t\!
      \begin{pmatrix}
        \epsilon(\partial_t) & 0 \\ 0 & \mu
      \end{pmatrix} +
      \begin{pmatrix}
        0 & -\Curl \\ \Curl_0 & 0
      \end{pmatrix}}
    \right)^{-1}
    \begin{pmatrix} \varPhi - \partial_t P_{\mathrm{nl}}(u) \\ \varPsi \end{pmatrix},
  \end{equation*}
  then by Theorem~\ref{thm:Maxwell_ExpStab_H2}, for $\nu < \tilde{\nu}$ and $u, v \in \Hw_{-\nu}(\R,\spatH{2})$,
  \begin{equation*}
    \begin{aligned}
      \lVert S(u) \rVert_{-\nu,2}
      &\lesssim \lVert \varPhi \rVert_{\mathcal{V}_{-\nu}} + \lVert \varPsi \rVert_{\mathcal{W}_{-\nu}}
      + \sum_{j=0}^3\lVert \partial_t^jP_\mathrm{nl}(u) \rVert_{-\nu,2}\\
      &\lesssim \lVert \varPhi \rVert_{\mathcal{V}_{-\nu}} + \lVert \varPsi \rVert_{\mathcal{W}_{-\nu}} + \lVert u \rVert_{-\nu,2}^{1+\alpha}
    \end{aligned}
  \end{equation*}
  and
  \begin{equation*}
    \begin{aligned}
      \lVert S(u) - S(v) \rVert_{-\nu,2}
      &\lesssim
      \sum_{j=0}^3\lVert \partial_t^jP_{\mathrm{nl}}(u) - \partial_t^jP_{\mathrm{nl}}(v) \rVert_{-\nu,2}\\
      &\lesssim \bigl( \lVert u \rVert_{-\nu,2} + \lVert v \rVert_{-\nu,2} \bigr)^\alpha \lVert u - v \rVert_{-\nu,2}.
    \end{aligned}
  \end{equation*}
  Thus, with the projection $\pi_1(E,H) := E$, the map $\pi_1S$ is a contraction on the ball $B_\varepsilon^{\nu}(\R,\mathcal{H}^2) \coleq \{ u \in L^2_{-\nu}(\R,\mathcal{H}^2) : \lVert u \rVert_{-\nu,2} \le \varepsilon \}$ for small $\varepsilon,\nu > 0$, if $\lVert \varPhi \rVert_{\mathcal{V}_{-\nu}}, \lVert \varPsi \rVert_{\mathcal{W}_{-\nu}} < c_0\varepsilon$ for small $c_0$.
\end{proof}

The above result applies in particular to the fully nonlocal nonlinearities in Example \ref{ex:NonlocalExpStab} (or the local $\spatH{2}$-version in Example~\ref{ex:quad_nonlocal_H2}) after imposing suitable regularity on the spatial kernel $\Lambda$ as well as the temporal kernel $\kappa$.
We provide here an additional class of admissible nonlinearities.

\begin{ex}
  \label{ex:LocLipH2}
  Suppose each $F \in \{\partial_t^jP_\mathrm{nl}: j\in \{0,1,2,3\} \}$ is of the form
  \begin{equation*}
    F(E) = \int_\R \kappa (\mult - s)\,q(E(s))\dd s
  \end{equation*}
  with $\kappa \in \Hw_{-\nu_\kappa}(\R,\R^{3\times 3})$, $\supp\kappa \subseteq [0,\infty)$, and $\nu_\kappa > 0$, and $q\colon \spatH{2}\to \spatH{2}$ is such that
  \begin{equation}
    \lVert \snl(u) - \snl(v) \rVert_{\spatH{2}}
    \le C \bigl(\lVert u \rVert_{\spatH{2}} + \lVert v \rVert_{\spatH{2}}\bigr) \lVert u - v \rVert_{\spatH{2}}
    \label{eq:LocLipSpacH2}
  \end{equation}
  for all $u,v \in \spatH{2}$ (for instance if $q$ is given by a bilinear map $q^{(2)}\colon\R^3\times\R^3 \to \R^3$ via $q(u) = q^{(2)}(u,u)$).
  Then, by an analogous estimate to that in \eqref{eq:convolution_kappa_Lipschitz} we can show that $F$ fulfills
  \begin{equation*}
    \lVert F(u) - F(v) \rVert_{-\nu,2}
    \le \sqrt{2}C\lVert \kappa \rVert_{L^2_{-\nu_\kappa}} \bigl(\lVert u \rVert_{-\nu,2} + \lVert v \rVert_{-\nu,2} \bigr) \lVert u - v \rVert_{-\nu,2}
  \end{equation*}
  for $u, v \in \Hw_{-\nu}(\R,\spatH{2})$, for small $\nu > 0$.
  Indeed, take $\nu \in (0,\nu_\kappa]$ and let $u,v \in L^2_{-\nu}(\R,\spatH{2})$ with $\supp u, \supp v \subseteq [0, \infty)$.
  Then we have
  \begin{multline*}
    \int_\R \biggl\lVert \int_\R \kappa(t-s) \bigl( \snl(u(s)) - \snl(v(s)) \bigr) \dd s \biggr\rVert^2_{\spatH{2}} e^{2\nu t} \dd t \\
    \begin{aligned}
      &\le \int_\R \biggl( \int_\R \lVert \kappa(t-s)\rVert\, \bigl\lVert \snl(u(s)) - \snl(v(s)) \bigr\rVert_{\spatH{2}} \dd s \biggr)^2 e^{2\nu t} \dd t \\
      &\le   C^2 \int_\R \biggl( \int_\R \lVert\kappa(t-s)\rVert \left( \lVert u(s) \rVert_{\spatH{2}} + \lVert v(s) \rVert_{\spatH{2}} \right) e^{-\nu s} \cdot \lVert u(s) - v(s) \rVert_{\spatH{2}} e^{\nu s} \dd s \biggr)^2 e^{2\nu t} \dd t \\
      &\le   C^2 \lVert u - v \rVert^2_{-\nu,2} \int_\R \int_\R \lVert\kappa(t-s)\rVert^2 e^{2\nu_\kappa(t-s)} \underbrace{e^{2(\nu-\nu_\kappa)(t-s)}}_{\le 1} \dd t \left( \lVert u(s) \rVert_{\spatH{2}} + \lVert v(s) \rVert_{\spatH{2}} \right)^2 \dd s \\
      &\le 2 C^2 \lVert \kappa \rVert_{L^2_{-\nu_\kappa}}^2 \bigl( \lVert u \rVert^2_{0,2} + \lVert v \rVert^2_{0,2} \bigr)\, \lVert u - v \rVert^2_{-\nu,2} \\
      &\le 2 C^2 \lVert \kappa \rVert_{L^2_{-\nu_\kappa}}^2 \bigl( \lVert u \rVert_{-\nu,2} + \lVert v \rVert_{-\nu,2} \bigr)^2\, \lVert u - v \rVert^2_{-\nu,2}
    \end{aligned}
  \end{multline*}
  making repeated use of the Cauchy--Schwarz inequality and of $\lVert u \rVert_{0,2} \le \lVert u \rVert_{-\nu,2}$ for $\supp u \subseteq [0,\infty)$.

  It appears that this estimate cannot be generalized to the case when
  \begin{equation*}
    \lVert \snl(u) - \snl(v) \rVert_{\spatH{2}}
    \le C \bigl(\lVert u \rVert_{\spatH{2}} + \lVert v \rVert_{\spatH{2}}\bigr)^\alpha \lVert u - v \rVert_{\spatH{2}}
  \end{equation*}
  with $\alpha > 1$, like the $k$-linear ($k \ge 2$) map $q$ considered at the end of Example~\ref{ex:NonlocalNonlinearity}.
\end{ex}

\section{Closing remarks}
\subsection{Materials with conductivity}
The established stability results are based on Theorem \ref{thm:exp_stab_M_d}. Alternatively, one can impose a stricter material damping in the form of a uniformly positive conductivity term, cf.~\cite{LasieckaExpStab}.
We make use of the following criterion.

\begin{prop}[{\cite[Corollary 2.2.4]{Trostorff_habil}}]
  \label{prop:M_d_evo_sys_alternative}
  Let $M$ be a material law of the form $M(z) = M_0(z) + z^{-1}M_1(z)$ with $M_0, M_1 \colon \dom(M) \subseteq \C \to \BLO(\HX)$ analytic and bounded.
  Assume there exist $\nu_0, c > 0$ such that $\C_{\Re > -\nu_0}\setminus \dom(M)$ is discrete and
  \begin{equation*}
     \forall z \in \dom(M) \cap \C_{\Re > -\nu_0}:\quad \Re zM(z) \ge c.
  \end{equation*}
  Then there exists $d > 0$ such that the evolutionary problem \eqref{eq:M_d_evo_sys} is well-posed and exponentially stable.
\end{prop}

\noindent
Consider the following linear Maxwell system
\begin{equation*}
  \left( \partial_t
    \begin{pmatrix}
      \epsilon(\partial_t) & 0 \\ 0 & \mu
    \end{pmatrix} +
    \begin{pmatrix}
      \sigma & 0 \\ 0 & 0
    \end{pmatrix} +
    \begin{pmatrix}
      0 & -\Curl \\ \Curl_0 & 0
    \end{pmatrix}
  \right)
  \begin{pmatrix}
    E \\ H
  \end{pmatrix} =
  \begin{pmatrix}
    \varPhi \\ \varPsi
  \end{pmatrix}.
\end{equation*}
where the \emph{conductivity} $\sigma\in\mathcal{B}(\Hs)$ is selfadjoint and uniformly positive definite with $\sigma \ge c_\sigma > 0$.
For $\epsilon(\partial_t)$ we assume that the following conditions hold.
\begin{description}
  \item[(M4)\label{itm:MatCond_conductivity}]
    $\epsilon \colon \dom(\epsilon) \subseteq \C \to \BLO(\Hs)$ is holomorphic and uniformly bounded.
    Moreover, there exists $\nu_0, c > 0$ such that
    $\C_{\Re > -\nu_0}\setminus\dom(\epsilon)$ is discrete and
    \begin{equation*}
      \Re \epsilon(z) \ge c > 0,\quad
      \Re z\epsilon(z) + \sigma \ge c > 0
    \end{equation*}
    for all $z \in \C_{\Re > -\nu_0}\cap \dom(\epsilon)$.
\end{description}
Thus, \nameref{itm:MatCond_conductivity} assures that the material law $M$ given by
\begin{equation*}
  M(z) \coleq \epsilon(z) + z^{-1}\sigma
\end{equation*}
satisfies
\begin{equation}
  \Re z \ge -\nu_0 \implies
  \Re zM(z) \ge c.
  \label{eq:accr_perm_strict}
\end{equation}
Now we want to reiterate the proof of Theorem~\ref{thm:Maxwell_ExpStab_SecondOrder} with the material law $M(\partial_t)$ instead of $\epsilon(\partial_t)$ to obtain $E \in \Hw_{-\nu}(\R,\Hs)$ for some $\nu > 0$.
To this end, analogously to \eqref{eq:Mw_expstab2ord}, we derive the second-order system
\begin{equation*}
  (\partial_t^2M(\partial_t) + \Curl\mu^{-1}\Curl_0)E = g
\end{equation*}
and write it as
\begin{equation}
  \left(\partial_t^2
  \begin{pmatrix}
    M_{00}(\partial_t) & M_{01}(\partial_t)\\ M_{10}(\partial_t) & M_{11}(\partial_t)
  \end{pmatrix} +
  \begin{pmatrix}
    \iota_0^*\Curl\mu^{-1}\Curl_0\iota_0 & 0\\ 0 & 0
  \end{pmatrix}
  \right)
  \begin{pmatrix}
    E_0 \\ E_1
  \end{pmatrix} =
  \begin{pmatrix}
    g_0 \\ g_1
  \end{pmatrix},
  \label{eq:Mw_ord2_conductivity}
\end{equation}
where $M_{ij}(\partial_t) = \iota_i^*M(\partial_t)\iota_j$
(again $\iota_0 = \iota_{\ran(\Curl)}$, $\iota_1 = \iota_{\ker(\Curl_0)}$).
We now claim that the map
\begin{equation*}
  z \mapsto zM_{01}(z)(zM_{11}(z))^{-1} = (z\epsilon_{01}(z) + \sigma_{01})(z\epsilon_{11}(z) + \sigma_{11})^{-1}
\end{equation*}
(which is well-defined pointwise by $\Re zM(z) \ge c > 0$ and Lemma~\ref{lem:proj_positive})
is bounded on $\C_{\Re > -\nu_0}$.
To see this, write for $r > 0$
\begin{equation*}
  \begin{aligned}
    zM_{01}(z)(zM_{11}(z))^{-1}
    &= \epsilon_{01}(z)(\epsilon_{11}(z) + z^{-1}\sigma_{11})^{-1}1_{\{|z| > r\}}(z)\\
    &\quad + z\epsilon_{01}(z)(z\epsilon_{11}(z) + \sigma_{11})^{-1}1_{\{|z| \le r\}}(z)\\
    &\quad + \sigma_{01}(z\epsilon_{11}(z) + \sigma_{11})^{-1}
  \end{aligned}
\end{equation*}
and observe that, upon choosing $r$ large enough and using $\Re \epsilon(z) \ge c > 0$ and again Lemma~\ref{lem:proj_positive}, all factors are bounded on $\C_{\Re > -\nu_0}\cap \dom(\epsilon)$. Hence the claim follows by analytic continuation.
Consequently, $\partial_tM_{01}(\partial_t)(\partial_tM_{11}(\partial_t))^{-1}$ maps $\Hw_{-\nu}(\R,\mathscr{H}_1)$ into $\Hw_{-\nu}(\R,\mathscr{H}_0)$ for each $\nu < \nu_0$,
and \eqref{eq:Mw_ord2_conductivity} can be transformed into
\begin{multline*}
  \Biggl(\partial_t
  \begin{pmatrix}
    \partial_t\left( M_{00}(\partial_t) - \partial_tM_{01}(\partial_t)(\partial_tM_{11}(\partial_t))^{-1}M_{10}(\partial_t)\right) & 0\\ \partial_tM_{10}(\partial_t) & \partial_tM_{11}(\partial_t)
  \end{pmatrix}\\ +
  \begin{pmatrix}
    \iota_0^*\Curl\mu^{-1}\Curl_0\iota_0 & 0\\ 0 & 0
  \end{pmatrix}
  \Biggr)
  \begin{pmatrix}
    E_0 \\ E_1
  \end{pmatrix} =
  \begin{pmatrix}
    g_0 - \partial_tM_{01}(\partial_t)(\partial_tM_{11}(\partial_t))^{-1}g_1 \\ g_1
  \end{pmatrix}.
\end{multline*}
Now the system for $E_0$ in the first line satisfies the conditions for exponential stability of Proposition~\ref{prop:M_d_evo_sys_alternative}, and the rest follows analogously as in the proof of Theorem~\ref{thm:Maxwell_ExpStab_SecondOrder}.
We conclude that Theorem~\ref{thm:Maxwell_ExpStab_SecondOrder} holds for the Maxwell system with conductivity, with assumption \nameref{itm:MatCond_conductivity} replacing \nameref{itm:MatCond_accr1} and \nameref{itm:MatCond_accr2} if $g_0 - \partial_tM_{01}(\partial_t)(\partial_tM_{11}(\partial_t))^{-1}g_1 \in L^2_{-\nu}(\R,\mathscr{H}_0)$.
A similar conclusion can be drawn for Theorem~\ref{thm:Maxwell_ExpStab_FirstOrder}.

\subsection{Removal of the auxiliary terms \texorpdfstring{$\Pi_0\partial_t^{-1}\varPhi, \Pi_1\partial_t^{-1}\varPsi$}{Π0 ∂t−1 Φ, Π1 ∂t−1 Ψ}}
Under alternative (and in view of \eqref{eq:Maxwell_evo_inhomogeneity} rather natural) assumptions on the data, the dependence on the anti-derivatives $\Pi_0\partial_t^{-1}\varPhi, \Pi_1\partial_t^{-1}\varPsi$ in Theorem \ref{thm:Maxwell_spatial_regularity} and in the theorems in §\ref{sec:ExpStab_SecondOrder}, §\ref{sec:ExpStab_FirstOrder} can be removed.
Let $\Omega$ be a simply connected Lipschitz domain with connected complement.
Then the identities
\begin{equation*}
  \ker(\Div) = \ran(\Curl) = \ker(\Curl_0)^\bot,\quad \ker(\Div_0) = \ran(\Curl_0) = \ker(\Curl)^\bot
\end{equation*}
hold (see \cite[Lemma 5]{MilaniPicard_DecompNl} or \cite[Prop 6.2]{PaulyWaurick}). Consequently, suppose that $\varPhi \in \ker(\Div)$, $\varPsi \in \ker(\Div_0)$ are divergence-free.
Then the identities above imply $\Pi_0 \varPhi = 0$, $\Pi_1\varPsi = 0$ as well as $\Pi_0\partial_t^{-1}\varPhi = 0$, $\Pi_1\partial_t^{-1}\varPsi = 0$.
Thus, the conditions
\begin{equation*}
  \begin{aligned}
    \varPhi &\in \Hw_{-\nu}(\R,\ker(\Div))\\
    \varPsi &\in \Hw_{-\nu}(\R,\ker(\Div_0)) \\
    \partial_t\varPhi + \Curl\mu^{-1}\varPsi &\in \Hw_{-\nu}(\R,\Hs)
  \end{aligned}
\end{equation*}
are sufficient to obtain $E,H \in \Hw_{-\nu}(\R,\Hs)$ via Theorem \ref{thm:Maxwell_ExpStab_SecondOrder} or Theorem \ref{thm:Maxwell_ExpStab_FirstOrder}.
For Theorem \ref{thm:Maxwell_ExpStab_H2} and Theorem \ref{thm:NonlinFixH2}, the prerequisite can be changed from $\varPhi \in \mathcal{V}_{-\nu}$, $\varPsi \in \mathcal{W}_{-\nu}$ (cf.\ Definition \ref{defi:AdmissibleData_VW}) to 
\begin{equation*}
  \begin{aligned}
    \varPhi &\in H^2_{-\nu}(\R,\Hs) \cap \Hw_{-\nu}(\R,\spatH{1}\cap \ker(\Div))\\
    \varPsi &\in H^2_{-\nu}(\R,\Hs) \cap \Hw_{-\nu}(\R,\spatH{1}\cap \ker(\Div_0))\\
    \Curl\mu^{-1}\varPsi &\in H^1_{-\nu}(\R,\Hs).
  \end{aligned}
\end{equation*}

\appendix{}
\section{A class of permittivities of Drude--Lorentz-type and applications}
\label{sec:app1}
Here we check the accretivity conditions in Picard's theorem and conditions \nameref{itm:MatCond_accr1} and \nameref{itm:MatCond_accr2} for a particular class of scalar permittivities.

\subsection{Accretivity of the material law on a right half-plane}
The Drude--Lorentz model in its general form is given by
\begin{equation*}
  \tilde{\chi}_{\mathsf{DL}}(\omega) = \sum_{j=1}^n \frac{\alpha_j}{\omega_{0,j}^2 - \omega^2 - 2i\gamma_j\omega}, \quad \text{resp.}\quad
  \tilde{\chi}_{\mathsf{DL}}(iz) = \sum_{j=1}^n \frac{\alpha_j}{\omega_{0,j}^2 + z^2 + 2\gamma_jz},
\end{equation*}
where $n \in \N$, $\omega_{0,j} \ge 0$, $\alpha_j, \gamma_j > 0$.
Recall that $\tilde{\chi}(iz) = (2\pi)^{-1/2}\int \chi(t) e^{-zt}\dd t$ denotes the Laplace transform.
As such, we note that if $\omega_{0,j} > \gamma_j$, the above is the Laplace transform of a sum of exponentially damped sine functions:
\begin{equation*}
  \chi_{\mathsf{DL}}(t) = \theta(t) \sum_{j=1}^n a_j e^{-\gamma_j t}\sin(b_jt),\qquad a_j = \frac{\alpha_j}{b_j},\ b_j = (\omega_{0,j}^2 - \gamma_j^2)^{1/2}.
\end{equation*}
where $\theta$ is the Heaviside function. We take for simplicity $n = 1$ and consider the material law $\chi$ given by
\begin{equation}
  \chi(z) = \frac{\alpha}{\omega_0^2 + z^2 + 2\gamma z},
  \qquad \text{where}\quad \omega_0, \alpha, \gamma > 0.
  \label{eq:DL_generic}
\end{equation}
The zeros of the denominator are
\begin{align*}
  z =
  \begin{cases}
    -\gamma \pm i\sqrt{\omega_0^2 - \gamma^2}, & \text{if $\omega_0 > \gamma$} \\
    -\gamma \pm \sqrt{\gamma^2 - \omega_0^2}, & \text{if $\omega_0 \le \gamma$},
  \end{cases}
\end{align*}
and are all contained in $\C_{\Re < 0}$.
Writing $z = \nu + it$, we compute
\begin{align*}
  \Re {\bigl(\left(\nu + it \right){\chi}(\nu + it)\bigr)}
  &= \Re \frac{\alpha (\nu + it)}{\omega_0^2 + (\nu + it)^2 + 2\gamma(\nu + it)} \\
  &= \Re \frac{\alpha (\nu + it)}{(\omega_0^2 + \nu^2 - t^2 + 2\gamma \nu) + 2i(\nu t + \gamma t)} \\
  &= \frac{\alpha \nu (\omega_0^2 + \nu^2 - t^2 + 2\gamma \nu) + 2\alpha t(\nu t + \gamma t)}{\bigl( \omega_0^2 + \nu^2 - t^2 + 2\gamma \nu \bigr)^2 + 4\bigl( \nu t + \gamma t \bigr)^2} \\
  &= \frac{\alpha \nu (\omega_0^2 + \nu^2 + t^2 + 2\gamma \nu) + 2\alpha \gamma t^2}{\bigl( \omega_0^2 + \nu^2 - t^2 + 2\gamma \nu \bigr)^2 + 4\bigl( \nu t + \gamma t \bigr)^2}.
\end{align*}
With $M(z) \coleq \epsilon_0 + {\chi}(z)$, $\epsilon_0 > 0$, we thus have
\begin{align*}
  \Re {\bigl(\left( \nu + it \right) M(\nu + it)\bigr)}
  &= \epsilon_0\nu + \frac{\alpha \nu (\omega_0^2 + \nu^2 + t^2 + 2\gamma \nu) + 2\alpha \gamma t^2}{\bigl( \omega_0^2 + \nu^2 - t^2 + 2\gamma \nu \bigr)^2 + 4\bigl( \nu t + \gamma t \bigr)^2}.
\end{align*}
Because the second term on the right is positive, bounded, and vanishes as $|t| \to \infty$,
this shows 
$\Re zM(z) \ge c \Re z$ for some $c > 0$ if $\varrho_0 > 0$ and $z \in \C_{\Re > \varrho_0}$.
In particular, $M(\partial_t)$ fulfills the requirements of Theorem~\ref{thm:linear_solution_theory} and Proposition~\ref{prop:FixPointRefinement}.

In contrast, strict accretivity of $zM(z)$ on a half-plane $\C_{\Re > -\nu_0}$ (even outside of a disk $B[0,\delta]$, cf.\ assumption \nameref{itm:MatCond_accr1}) cannot hold for $\nu_0 > 0$.

\subsection{Analytic correction to the material law}
\label{sec:app_analytic_correction}

We have seen that the linear Maxwell system with the ``standard'' Drude--Lorentz susceptibility does not fulfill the criteria for exponential stability.
However, these are still global criteria that assure exponential decay of the solution for rather general right-hand sides $\varPhi,\varPsi \in \Hw_{-\nu}$.
If the Fourier--Laplace transform of the right-hand side is localized around a certain “frequency” $z = z_0 \in \C$, we can argue that the exact form of the solution operator plays no role for large $|z|$.

To be specific, fix $r \gg 1$ and consider the following modified material law (for the case $z_0 = 0$),
\begin{equation}
  M_r(z) \coleq \epsilon_0 + \frac{\alpha}{\omega_0^2 + z^2 + 2\gamma z} \left( 1 + \frac{z}{r} \right),
  \qquad \text{with}\quad \alpha, \gamma, \omega_0 > 0.
  \label{eq:DL_modified}
\end{equation}
Clearly, $M_r$ is bounded on a half-plane containing the imaginary axis.
First we show that \nameref{itm:MatCond_accr1} is satisfied if $r > \frac{\omega_0^2}{2\gamma}$ is large enough.

\begin{lem}
  Assume that $2\gamma r - \omega_0^2 > 0$.
  Then $M_r$ satisfies \nameref{itm:MatCond_accr1}, \ie, for all $\delta > 0$ there exist $\nu_1 > 0$ and $c > 0$ such that
\begin{equation}
  \forall z \in \C_{\Re > -\nu_1} \setminus B[0,\delta]: \quad  \Re zM_r(z) \ge c.
  \label{eq:Mr_accr}
\end{equation}
\end{lem}

\begin{proof}
  We denote $g(\nu,t) \coleq \Re\, (\nu + it) M_r(\nu + it) - \epsilon_0\nu$ and find that $g$ is explicitly given by
  \begin{equation*}
    g(\nu,t) = 
    \frac{\alpha}{r}\, \frac{\nu\omega_0^2 r + \left(2\gamma r + \omega_0^2 \right)\nu^2 + \left( r + 2\gamma \right) \nu^3 + \left( \left( r + 2\gamma \right) \nu + 2\nu^2 + 2\gamma r - \omega_0^2 \right) t^2 + \nu^4 + t^4}{\left( \omega_0^2 + \nu^2 - t^2 + 2\gamma\,\nu \right)^2 + \left( 2\,\nu + 2\gamma \right)^2 t^2}.
  \end{equation*}
  Define $\gamma_0 > 0$ by
  $\gamma_0 = \gamma$ if $\omega_0 > \gamma$ and $\gamma_0 = \gamma - (\gamma^2-\omega_0^2)^{1/2}$ if $\omega_0 \le \gamma$.
  Then $g$ is continuous on $\R_{> -\gamma_0} \times \R$, since the zeros of the denominator lie outside this set.
  Also, $g$ is positive on $\R_{\ge 0} \times (\R \setminus \{0\})$, since it is a sum of positive terms; here the assumption $2\gamma r - \omega_0^2 > 0$ has been used.
  Moreover, due to $\lim_{\nu \to \infty} g(\nu,t) = \lim_{|t| \to \infty} g(\nu,t) = \frac{\alpha}{r} > 0$, for all $\delta > 0$ the map
  \begin{equation*}
      \nu \mapsto \inf_{\{t:\, \nu^2 + t^2 > \delta^2\}} g(\nu,t) \qquad (\nu \ge 0)
  \end{equation*}
  is even bounded from below by a positive constant.
  By continuity, this remains true for $\nu < 0$ small enough, \ie, there are $\nu_0, c_0 > 0$ such that $g(\nu,t) \ge c_0$ for all $\nu > -\nu_0$.
  Now choose $0 < \nu_1 < \min\{\nu_0, c_0\epsilon_0^{-1}\}$, then
  \begin{equation*}
    \Re\,(\nu + it) M_r(\nu + it) = g(\nu,t) + \epsilon_0\nu \ge c_0 - \epsilon_0\nu_1 \eqcol c > 0
  \end{equation*}
    for all $\nu > -\nu_1$. 
\end{proof}

Let us now consider condition \nameref{itm:MatCond_accr2} for $M_r$, \ie,
\begin{equation}
  \exists c > 0\ \forall z\in \C_{\Re > -\nu_2}:\quad \Re M_r(z) \ge c
  \label{eq:Mr_inv}
\end{equation}
for some $\nu_2 \in (0, \nu_1]$.
A look at
\begin{equation*}
  \Re M_r(\nu + it) - \epsilon_0 =
  \frac{\alpha}{r} \frac{(2\gamma + \nu - r)t^2 + \omega_0^2r + (2\gamma r + \omega_0^2)\nu + (2\gamma + r)\nu^2 + \nu^3}{\bigl( \omega_0^2 + \nu^2 - t^2 + 2\nu\gamma \bigr)^2 + 4\bigl(\nu + \gamma \bigr)^2}
\end{equation*}
shows that this term is overall bounded and vanishes as $\left| \nu + it \right| \to \infty$.
By a suitable choice of parameters $\omega_0, \gamma, \alpha > 0$ (small $\alpha$, depending on $\gamma, \omega_0$), we can ensure $\Re M_r(\nu + it) \ge c$ for some $c > 0$ and small $\nu < 0$, providing \eqref{eq:Mr_inv}.

The modified material law \eqref{eq:DL_modified} thereby fulfills the conditions necessary for Theorem~\ref{thm:Maxwell_ExpStab_SecondOrder} to obtain exponential stability for the modified linearized Maxwell system
\begin{equation}
  \left( \partial_t
    \begin{pmatrix}
      M_r(\partial_t) & 0 \\ 0 & \mu
    \end{pmatrix}
    + 
    \begin{pmatrix}
      0 & -\Curl \\ \Curl_0 & 0
    \end{pmatrix}
  \right) \binom{E}{H} = \binom{\varPhi}{\varPsi}.
  \label{eq:Maxwell_bounded}
\end{equation}

\begin{rem}
  \begin{enumerate}
    \item[(i)]
      Multiplying ${\chi}(z)$ by the factor $(1 + \frac{z}{r})$ is the same as adding the following term to the convolution kernel $\chi_\mathsf{DL}(t)$:
      \begin{equation*}
        \frac{\alpha}{r}\theta(t)e^{-\gamma t} \bigl(\cos(bt) - \frac{\gamma}{b}\sin(bt) \bigr),\qquad b = (\omega_0^2 - \gamma^2)^{1/2}.
      \end{equation*}
      A similar modification was employed in \cite{SchneiderUecker2003} for nonlinear Maxwell equations in an optical fiber.
    \item[(ii)]
      In a similar fashion we could devise a modified material law of the form
      \begin{equation*}
        M_r(z) = \epsilon_0 + \sum_{j=1}^n \frac{\alpha_j}{\omega_{0,j}^2 + z^2 + 2\gamma_jz} \left( 1 + \frac{z-z_0}{r} \right),
      \end{equation*}
      which is localized around $z_0$. In this case the interplay between the parameters $\omega_{0,j}, \gamma_j, \alpha_j$, and $z_0$ becomes somewhat more delicate, but \eqref{eq:Mr_accr} and \eqref{eq:Mr_inv} are still satisfied for small $\alpha_j$.
  \end{enumerate}
\end{rem}

\subsection{Well-posedness and exponential stability of linear and nonlinear systems}
\label{sec:app_examples_overview}

In Tables~\ref{tab:LinearWPES} and \ref{tab:NonlinearWPES} we collect the various examples given throughout the paper into an overview of different Maxwell systems, sorted by the properties (defined below) proved in the paper.
Here \textsf{DL} denotes the Drude--Lorentz permittivity model in \eqref{eq:DL_generic} and \textsf{mod-DL} that in \eqref{eq:DL_modified}.
\medskip

\noindent In the linear case:
\begin{description}
  \item[(WP$_0$)] Well-posedness of the linear system in the spaces $L^2_\varrho(\R,\Hs)$ ($\varrho > \varrho_0$ for some $\varrho_0 > 0$)
    (via direct application of Picard's theorem, see Example~\ref{ex:LinearPermittivity})
  \item[(WP$_2$)] $\spatH{2}$-regularity of the linear system (\ie, well-posedness in the spaces $L^2_\varrho(\R,\spatH{2})$, $\varrho > \varrho_0 > 0$)
    (Theorem~\ref{thm:Maxwell_spatial_regularity}).
  \item[(ES$_0$)] Exponential stability of the linear system in $L^2_{-\nu}(\R,\Hs)$ ($0 < \nu < \nu_0$) (Theorem~\ref{thm:Maxwell_ExpStab_SecondOrder} resp.\ Theorem~\ref{thm:Maxwell_ExpStab_FirstOrder}).
  \item[(ES$_2$)] Exponential stability of the linear system in $L^2_{-\nu}(\R,\spatH{2})$ ($0 < \nu < \nu_0$) (Theorem~\ref{thm:Maxwell_ExpStab_H2}).
\end{description}
\medskip

\noindent In the nonlinear case:
\begin{description}
  \item[(nl-WP$_0$)] Well-posedness of the nonlinear system (with \textsf{DL} or \textsf{mod-DL}) in the spaces $L^2_\varrho(\R,\Hs)$ ($\varrho > \varrho_0 \ge 0$)
    (Picard's theorem and fixed-point argument)
  \item[(nl-WP$_{\mathrm{cut},0}$)] $\Hw_\varrho(\R,\Hs)$-solutions ($\varrho > \varrho_0 \ge 0$) of the nonlinear system with time-cutoff (with \textsf{DL} or \textsf{mod-DL}) (Proposition~\ref{prop:LocalWP_multilinear}).
  \item[(nl-WP$_{\mathrm{cut},2}$)] $\Hw_\varrho(\R,\spatH{2})$-solutions ($\varrho > \varrho_0 > 0$) of the nonlinear system with time-cutoff (with \textsf{DL} or \textsf{mod-DL})
    (Corollary~\ref{cor:H2_regularity_nonlinear_LocLip}).
  \item[(nl-ES$_0$)] Small solutions of the nonlinear system (only with \textsf{mod-DL}) in $L^2_{-\nu}(\R,\Hs)$ ($0 < \nu < \nu_0$) (Theorem~\ref{thm:NonlinFixL2} resp.\ Theorem~\ref{thm:Maxwell_ExpStab_FirstOrder} together with fixed-point argument).
  \item[(nl-ES$_2$)] Small solutions of the nonlinear system (only with \textsf{mod-DL}) in $L^2_{-\nu}(\R,\spatH{2})$ ($0 < \nu < \nu_0$) (Theorem~\ref{thm:NonlinFixH2}).
\end{description}

\begin{table}[H]
  \centering
  \begin{tabular}{r c c c}
    \toprule
    & $\Omega$ arbitrary & $\Omega$ bw-Lipschitz & $\Omega$ smooth \\\midrule
    \textsf{DL} \eqref{eq:DL_generic} & \textsf{(WP$_0$)} & \textsf{(WP$_0$)} & \textsf{(WP$_0$)}, \textsf{(WP$_2$)}\\\midrule
    \textsf{mod-DL} \eqref{eq:DL_modified} & \textsf{(WP$_0$)} & \textsf{(WP$_0$)}, \textsf{(ES$_0$)} & \textsf{(WP$_0$)}, \textsf{(ES$_0$)}, \textsf{(WP$_2$)}, \textsf{(ES$_2$)}\\\bottomrule
  \end{tabular}
  \caption{Well-posedness and exponential stability of the linear Maxwell system \eqref{eq:Mw1ord} depending on the type of permittivity and the domain. Here “bw-Lipschitz” stands for a bounded, weak Lipschitz domain, while “smooth” refers to a bounded domain satisfying the conditions of Proposition~\ref{prop:spatial_regularity}.}
  \label{tab:LinearWPES}
\end{table}
\begin{table}[H]
  \centering
  \begin{tabular}{>{\raggedleft}p{.35\textwidth} >{\centering}p{.15\textwidth} >{\centering}p{.15\textwidth} >{\centering\arraybackslash\hspace{0pt}}p{.24\textwidth}}
    \toprule
    & $\Omega$ arbitrary & $\Omega$ bw-Lipschitz & $\Omega$ smooth \\\midrule
    Saturable (Ex.~\ref{ex:SaturableNonlinearity1}) & \textsf{(nl-WP$_0$)} & \textsf{(nl-WP$_0$)} & \textsf{(nl-WP$_0$)}\\\midrule
    Fully nonlocal multilinear (Ex.~\ref{ex:NonlocalNonlinearity}, Ex.~\ref{ex:NonlocalExpStab}) & \textsf{(nl-WP$_{\mathrm{cut},0}$)} & \mbox{\textsf{(nl-WP$_{\mathrm{cut},0}$)}}, \mbox{\textsf{(nl-ES$_0$)}} & \mbox{\textsf{(nl-WP$_{\mathrm{cut},0}$)}}, \mbox{\textsf{(nl-ES$_0$)}}, \mbox{\textsf{(nl-WP$_{\mathrm{cut},2}$)}}, \mbox{\textsf{(nl-ES$_{\mathrm{cut},2}$)}}\\\midrule
    Multilinear algebraic spatial map~$q$ (Ex.~\ref{ex:quad_nonlocal_H2}) & & & \textsf{(nl-WP$_{\mathrm{cut},2}$)}, \textsf{(nl-ES$_2$)}\\\midrule
    Quadratic, single-variable $\kappa$ (Ex.~\ref{ex:LocLipH2}) & & & \textsf{(nl-ES$_2$)}\\\bottomrule
  \end{tabular}
  \caption{Well-posedness and exponential stability of the nonlinear Maxwell system \eqref{eq:Maxwell_nonlinear} depending on the type of nonlinearity and the domain.}
  \label{tab:NonlinearWPES}
\end{table}

\begin{rem}
  We recapitulate that the interface $\Gamma$ is arbitrary except for the higher regularity results in the last column in Tables~\ref{tab:LinearWPES} and~\ref{tab:NonlinearWPES}. Here a smoothness assumption on $\Gamma$ is required in Proposition~\ref{prop:spatial_regularity}.
\end{rem}

\section{Note on closedness of the range of \texorpdfstring{$\boldsymbol{\Curl}$}{curl}}

The following compactness result (see \cite{Picard_compactness}) and the Lemma below are instrumental for the proofs of Theorem~\ref{thm:Maxwell_ExpStab_SecondOrder} and Theorem~\ref{thm:Maxwell_ExpStab_FirstOrder}.
\begin{thm}[Picard--Weber--Weck selection theorem]
  \label{thm:PWW_selection_thm}
  Let $\Omega \subset \R^3$ be a bounded weak Lipschitz domain. Then the embeddings
  \begin{equation*}
    H(\Curl,\Omega)\cap H_0(\Div,\Omega) \hookrightarrow L^2(\Omega)^3\quad\text{and}\quad
    H_0(\Curl,\Omega)\cap H(\Div,\Omega) \hookrightarrow L^2(\Omega)^3
  \end{equation*}
  are compact.
\end{thm}

\begin{lem}
  \label{lem:RanCurl_closed}
  Let $\Omega \subset \R^3$ be a bounded weak Lipschitz domain. Then the ranges
  \begin{equation*}
    \ran(\Curl) = \{ \Curl u : u \in H(\Curl,\Omega) \} \quad\text{and}\quad
    \ran(\Curl_0) = \{ \Curl u : u \in H_0(\Curl,\Omega) \}
  \end{equation*}
  are closed subspaces of $L^2(\Omega)^3$.
\end{lem}

\begin{proof}
  In the following we establish the Poincaré-type estimates
  \begin{equation}
    \begin{aligned}
      \exists C > 0\ \forall \phi \in \dom(\Curl) \cap \ker(\Curl)^\bot :
      \quad \lVert \phi \rVert_{L^2} &\le C \lVert \Curl \phi \rVert_{L^2}\\
      \exists C > 0\ \forall \phi \in \dom(\Curl_0) \cap \ker(\Curl_0)^\bot :
      \quad \lVert \phi \rVert_{L^2} &\le C \lVert \Curl_0 \phi \rVert_{L^2}.
    \end{aligned}
    \label{eq:Poincare_curl}
  \end{equation}
  These estimates are sufficient (even equivalent, see \cite{TroWau14}) for the closedness of the ranges. Indeed, let $(\phi_n)_n$ be a sequence in $\dom(\Curl) = H(\Curl,\Omega)$ such that $\Curl\phi_n \to \psi$ for some $\psi \in L^2(\Omega)^3$.
  Decompose $\phi_n = \phi_{n,0} + \phi_{n,1}$, where $\phi_{n,0} \in \ker(\Curl)$ and $\phi_{n,1} \in \ker(\Curl)^{\bot}$.
  Since $\phi_{n,0}\in\ker(\Curl) \subseteq\dom(\Curl)$ and $\phi_n\in\dom(\Curl)$, it follows that $\phi_{n,1}\in\dom(\Curl)$.
  Using estimate \eqref{eq:Poincare_curl}, we infer for $n,m\in\N$:
  \begin{equation*}
    \lVert \phi_{n,1} - \phi_{m,1} \rVert_{L^2} \le C \lVert \Curl\phi_{n,1} - \Curl\phi_{m,1} \rVert_{L^2}
    = C\lVert \Curl\phi_n - \Curl\phi_m \rVert_{L^2},
  \end{equation*}
  which yields that $(\phi_{n,1})$ is a Cauchy sequence in $\dom(\Curl)$.
  Since $\Curl$ is a closed operator, there exist $\phi,\psi_1$ such that
  \begin{equation*}
    \phi_{n,1} \to \phi,\quad \Curl \phi_{n,1} \to \psi_1 = \Curl\phi\quad \text{in } L^2(\Omega)^3.
  \end{equation*}
  As $\Curl\phi_{n,1} = \Curl\phi_n$, it follows that $\psi = \psi_1 = \Curl\phi$, proving the closedness of $\ran(\Curl)$. The argument for $\ran(\Curl_0)$ is similar.

  Let us now show \eqref{eq:Poincare_curl}. Again, we only focus on $\Curl{}$.
  Assume that \eqref{eq:Poincare_curl} is not true.
  Then for each $n\in\N$ we find $\tilde{\phi}_n\in \dom(\Curl)\cap\ker(\Curl)^{\bot}$ such that $\lVert\tilde{\phi}_n\rVert_{L^2} > n \lVert \Curl\tilde{\phi}_n\rVert_{L^2}$. In particular, $\tilde{\phi}_n \ne 0$. For $n\in\N$ we let $\phi_n \coleq \lVert\tilde{\phi}_n\rVert_{L^2}^{-1}\tilde{\phi}_n$. Then
  \begin{equation*}
    \lVert\phi_n\rVert_{L^2} = 1\quad \text{and}\quad 1 > n\lVert \Curl\phi_n \rVert_{L^2}.
  \end{equation*}
  Thus, $(\phi_n)_n$ is a bounded sequence in $\dom(\Curl)\cap\ker(\Curl)^\bot \subseteq H(\Curl,\Omega)$. Without loss of generality (by possibly choosing a subsequence) we may assume that $(\phi_n)_n$ weakly converges in $\dom(\Curl) \cap \ker(\Curl)^{\bot}$, \ie,
  \begin{equation*}
    \phi_n\rightharpoonup\phi,\quad \Curl\phi_n \rightharpoonup \Curl\phi
  \end{equation*}
  for some $\phi \in \dom(\Curl)\cap \ker(\Curl)^{\bot}$.
  Next, since
  \begin{equation*}
    \ker(\Curl)^\bot = \ran(\Curl^*) = \ran(\Curl_0) \subseteq \ker(\Div_0) \subseteq \dom(\Div_0) = H_0(\Div,\Omega),
  \end{equation*}
  we infer by the Picard--Weber--Weck selection theorem that $\phi_n \to \phi$ in $L^2(\Omega)^3$.
  In particular, $1 = \lVert\phi_n\rVert_{L^2} \to \lVert\phi\rVert_{L^2} = 1$.
  Finally,
  \begin{equation*}
    \lVert\Curl\phi_n\rVert_{L^2} < \frac{1}{n} \to 0
  \end{equation*}
  and therefore
  \begin{equation*}
    \lVert\Curl\phi\rVert_{L^2} \le \liminf_{n\to\infty}{} \lVert \Curl\phi_n\rVert_{L^2} = 0.
  \end{equation*}
  Thus, $\phi\in\ker(\Curl)$. Since $\phi\in\ker(\Curl)^\bot$, we obtain $\phi = 0$, contradicting $\lVert\phi\rVert_{L^2} = 1$, which eventually proves \eqref{eq:Poincare_curl}.
\end{proof}

\section*{Acknowledgment}
Tomá\v{s} Dohnal acknowledges funding by the German Research Foundation (Deutsche Forschungsgemeinschaft DFG), project-ID DO1467/4-1.
Marcus Waurick acknowledges useful discussions with Rainer Picard.

\bibliography{main}{}

\begin{thebibliography}{BDPW22}

\bibitem[ABH99]{AmmariBaoHamdache}
Habib Ammari, Gang Bao, and Kamel Hamdache.
\newblock The effect of thin coatings on second harmonic generation.
\newblock {\em Electronic Journal of Differential Equations}, 1999(36):1–13,
  1999.

\bibitem[BD94]{BaoDobsonSHG}
Gang Bao and David~C. Dobson.
\newblock Second harmonic generation in nonlinear optical films.
\newblock {\em Journal of Mathematical Physics}, 35:1622–1633, 1994.

\bibitem[BDPW22]{BDPW22}
Malcolm Brown, Tom\'{a}\v{s} Dohnal, Michael Plum, and Ian Wood.
\newblock Spectrum of the {Maxwell} equations for a flat interface between
  homogeneous dispersive media.
\newblock \href{https://www.arXiv.org/abs/2206.02037}{arXiv:2206.02037}, 2022.

\bibitem[Boy08]{Boyd}
Robert~W. Boyd.
\newblock {\em Nonlinear Optics}.
\newblock Academic Press, 2008.

\bibitem[Bre11]{brezis}
Haïm Brezis.
\newblock {\em Functional Analysis, Sobolev Spaces and Partial Differential
  Equations}.
\newblock Springer, 2011.

\bibitem[BS22]{SchnaubeltBresch}
Christopher Bresch and Roland Schnaubelt.
\newblock Local wellposedness of {Maxwell} systems with scalar-type retarded
  material laws.
\newblock
  \url{https://www.math.kit.edu/iana3/~schnaubelt/media/maxwell-ret.pdf}, 2022.

\bibitem[DL90a]{DautrayLions1}
Robert Dautray and Jacques-Louis Lions.
\newblock {\em Mathematical Analysis and Numerical Methods for Science and
  Technology. {V}ol. 1}.
\newblock Springer, 1990.

\bibitem[DL90b]{DautrayLions3}
Robert Dautray and Jacques-Louis Lions.
\newblock {\em Mathematical Analysis and Numerical Methods for Science and
  Technology. {V}ol. 3}.
\newblock Springer, 1990.

\bibitem[DST22]{DST}
Tom\'{a}\v{s} Dohnal, Roland Schnaubelt, and Daniel~P. Tietz.
\newblock Rigorous envelope approximation for interface wave-packets in
  {Maxwell's} equations with {2D} localization.
\newblock \href{https://www.arXiv.org/abs/2206.03154}{arXiv:2206.03154}, 2022.
\newblock Accepted to SIAM J. Math. Anal.

\bibitem[GR86]{GiraultRaviart}
Vivette Girault and Pierre-Arnaud Raviart.
\newblock {\em Finite Element Methods for {Navier-Stokes} Equations}.
\newblock Springer, 1986.

\bibitem[Lei86]{Leis}
Rolf Leis.
\newblock {\em Initial Boundary Value Problems in Mathematical Physics}.
\newblock Springer, 1986.

\bibitem[LPS19]{LasieckaExpStab}
Irena Lasiecka, Michael Pokojovy, and Roland Schnaubelt.
\newblock Exponential decay of quasilinear {Maxwell} equations with interior
  conductivity.
\newblock {\em Nonlinear Differential Equations and Applications NoDEA},
  26(51), 2019.

\bibitem[Mai07]{MaierPlasmonics}
Stefan~A. Maier.
\newblock {\em Plasmonics: Fundamentals and Applications}.
\newblock Springer US, 2007.

\bibitem[MP88]{MilaniPicard_DecompNl}
Albert Milani and Rainer Picard.
\newblock Decomposition theorems and their applications to non-linear electro-
  and magneto-static boundary value problems.
\newblock In {\em Partial Differential Equations and Calculus of Variations,
  Lect. Notes Math. 1357}, pages 317--340. Springer, 1988.

\bibitem[MP02]{MilaniPicard}
Albert Milani and Rainer Picard.
\newblock Evolution equations with constitutive laws and memory effects.
\newblock {\em Differential and Integral Equations}, 15(3):327--344, 2002.

\bibitem[Pic84]{Picard_compactness}
Rainer Picard.
\newblock An elementary proof for a compact imbedding result in generalized
  electromagnetic theory.
\newblock {\em Mathematische Zeitschrift}, 187:151--164, 1984.

\bibitem[Pic09]{PicardStructObserv}
Rainer Picard.
\newblock A structural observation for linear material laws in classical
  mathematical physics.
\newblock {\em Mathematical Methods in the Applied Sciences},
  32(14):1768--1803, 2009.

\bibitem[PM11]{PicardMcGhee}
Rainer Picard and Des McGhee.
\newblock {\em Partial Differential Equations: A Unified Hilbert Space
  Approach}.
\newblock De Gruyter Expositions in Mathematics 55. De Gruyter, 2011.

\bibitem[PTW17]{PTW_MaxReg}
Rainer Picard, Sascha Trostorff, and Marcus Waurick.
\newblock On maximal regularity for a class of evolutionary equations.
\newblock {\em Journal of Mathematical Analysis and Applications},
  449:1368--1381, 2017.

\bibitem[PW22]{PaulyWaurick}
Dirk Pauly and Marcus Waurick.
\newblock The index of some mixed order {Dirac} type operators and generalised
  {Dirichlet--Neumann} tensor fields.
\newblock {\em Mathematische Zeitschrift}, 301:1739--1819, 2022.

\bibitem[Rae88]{raether}
Heinz Raether.
\newblock {\em Surface Plasmons on Smooth and Rough Surfaces and on Gratings}.
\newblock Springer, 1988.

\bibitem[She89]{Shen1989}
Yuen-Ron Shen.
\newblock Optical second harmonic generation at interfaces.
\newblock {\em Annual Review of Physical Chemistry}, 40:327--50, 1989.

\bibitem[SS22]{SchnaubeltSpitz2022}
Roland Schnaubelt and Martin Spitz.
\newblock Local wellposedness of quasilinear {Maxwell} equations with
  conservative interface conditions.
\newblock {\em Communications in Mathematical Sciences}, 20(8):2265--2313,
  2022.

\bibitem[STW22]{STW22}
Christian Seifert, Sascha Trostorff, and Marcus Waurick.
\newblock {\em Evolutionary Equations}.
\newblock Birkhäuser, 2022.

\bibitem[SU03]{SchneiderUecker2003}
Guido Schneider and Hannes Uecker.
\newblock Existence and stability of modulating pulse solutions in {Maxwell}'s
  equations describing nonlinear optics.
\newblock {\em Zeitschrift für Angewandte Mathematik und Physik ZAMP},
  54:677--712, 2003.

\bibitem[SW17]{SuessWaurick}
André Süß and Marcus Waurick.
\newblock A solution theory for a general class of {SPDE}s.
\newblock {\em Stochastics and Partial Differential Equations: Analysis and
  Computations}, 5:278--318, 2017.

\bibitem[Tro13]{TrostorffLinExpStab}
Sascha Trostorff.
\newblock Exponential stability for linear evolutionary equations.
\newblock {\em Asymptotic Analysis}, 85:179--197, 2013.

\bibitem[Tro18]{Trostorff_habil}
Sascha Trostorff.
\newblock Exponential stability and initial value problems for evolutionary
  equations.
\newblock Habilitation, TU Dresden, 2018.

\bibitem[TW14]{TroWau14}
Sascha Trostorff and Marcus Waurick.
\newblock A note on elliptic type boundary value problems with maximalmonotone
  relations.
\newblock {\em Mathematische Nachrichten}, 287:1545–1558, 2014.

\bibitem[TW21]{TW_NonautoMaxReg}
Sascha Trostorff and Marcus Waurick.
\newblock Maximal regularity for non-autonomous evolutionary equations.
\newblock {\em Integral Equations and Operator Theory}, 93(30), 2021.

\bibitem[Web80]{weber80}
Christian Weber.
\newblock A local compactness theorem for {Maxwell's} equations.
\newblock {\em Mathematical Methods in the Applied Sciences}, 2:12--25, 1980.

\bibitem[Web81]{weber81}
Christian Weber.
\newblock Regularity theorems for {Maxwell's} equations.
\newblock {\em Mathematical Methods in the Applied Sciences}, 3:523--536, 1981.

\end{thebibliography}
\bibliographystyle{alpha}
\end{document}